\documentclass{amsart}
\usepackage{amsmath,amsthm,amssymb,amsfonts}

\numberwithin{equation}{section} 
\theoremstyle{plain}
\newtheorem{theo+}           {Theorem}      [section]
\newtheorem{prop+}  [theo+]  {Proposition}
\newtheorem{coro+}  [theo+]  {Corollary}
\newtheorem{lemm+}  [theo+]  {Lemma}
\newtheorem{deep+}  [theo+]  {Deep Result}
\newtheorem{fact+}  [theo+]  {Fact}
\newtheorem{defi+}  [theo+]  {Definition}

\theoremstyle{definition}
\newtheorem{exam+}  [theo+]  {Example}
\newtheorem{rema+}  [theo+]  {Remark}

\newenvironment{theorem}{\begin{theo+}}{\end{theo+}}
\newenvironment{proposition}{\begin{prop+}}{\end{prop+}}
\newenvironment{corollary}{\begin{coro+}}{\end{coro+}}
\newenvironment{lemma}{\begin{lemm+}}{\end{lemm+}}

\newenvironment{example}{\begin{exam+}}{\end{exam+}}
\newenvironment{remark}{\begin{rema+}}{\end{rema+}}
\newenvironment{definition}{\begin{defi+}}{\end{defi+}}  

\newcommand{\al}{\alpha}
\newcommand{\ba}{\beta}
\newcommand{\be}{\beta}       
\newcommand{\ga}{\gamma}  
\newcommand{\Ga}{\Gamma}
\newcommand{\de}{\delta}
\newcommand{\da}{\delta}
\newcommand{\De}{\Delta}
\newcommand{\Da}{\Delta}
\newcommand{\ep}{\varepsilon}
\newcommand{\la}{\lambda}
\newcommand{\ro}{\rho}
\newcommand{\om}{\omega}
\newcommand{\Om}{\Omega}
\newcommand{\si}{\sigma}
\newcommand{\ze}{\zeta}
\newcommand{\La}{\Lambda}
\newcommand{\gh}{\mathfrak h}
\newcommand{\End}{\operatorname{End}}
\newcommand{\id}{\operatorname{id}}
\newcommand{\Hom}{\operatorname{Hom}}
\newcommand{\mg}{M_{\gh^\ast}}
\newcommand{\op}{\operatorname{op}}
\newcommand{\Ima}{\operatorname{Im}}
\newcommand{\fr}{\mathcal F_R(\mathrm{SU}(2))}
\newcommand{\frl}{\mathcal F_R(\mathrm{SL}(2))}
\newcommand{\frm}{\mathcal F_R(\mathrm{M}(2))}
\newcommand{\ot}{\otimes}
\newcommand{\wh}{\widehat}
\newcommand{\wt}{\widetilde}
\newcommand{\qb}[2]{\genfrac{[}{]}{0pt}{}{#1}{#2}_{q^2}}

\newcommand{\mh}{M_{\gh^\ast}}
\newcommand{\Zp}{{\mathbb Z}_{\geq 0}}
\newcommand{\R}{\mathbb R}

\begin{document}
\baselineskip 18pt
\larger[2]
\title[Harmonic analysis on the $\mathrm{SU}(2)$ dynamical quantum group]
{Harmonic analysis on the $\mathrm{SU}(2)$ dynamical quantum group}
\author{Erik Koelink}
\address
{Technische Universiteit Delft\\ Faculteit Informatietechnologie en Systemen\\ 
Toegepaste Wiskundige Analyse\\ 
Postbus 5031\\ 2600 GA Delft\\ The Netherlands}
\email{koelink@twi.tudelft.nl}
\author{Hjalmar Rosengren}
\address
{Department of Mathematics\\ Chalmers University of Technology and G\"oteborg
 University\\SE-412~96 G\"oteborg, Sweden}
\email{hjalmar@math.chalmers.se}

\begin{abstract}
Dynamical quantum groups were recently introduced by Etingof and Varchenko 
as an algebraic framework for studying the dynamical Yang--Baxter equation, 
which is precisely the Yang--Baxter equation satisfied by $6j$-symbols.
We investigate one of the simplest examples,  generalizing the standard 
$\mathrm{SU}(2)$ quantum group. The matrix elements for its corepresentations 
are identified with Askey--Wilson polynomials, and the Haar measure with the
Askey--Wilson measure. The discrete orthogonality of the matrix elements 
yield the orthogonality of $q$-Racah polynomials  (or quantum $6j$-symbols).
  The Clebsch--Gordan
coefficients for  representations and corepresentations are also identified 
with $q$-Racah polynomials. This results in new algebraic proofs of the
 Biedenharn--Elliott identity satisfied by quantum $6j$-symbols.  
\end{abstract}

\maketitle        

\section{Introduction}  

Quantum groups first arose in the 1980's as an algebraic framework for 
studying $R$-matrices, which have their
origin in statistical mechanics. An $R$-matrix is a solution of the 
Yang--Baxter
 equation, which exists in several versions. While the simplest examples of 
 quantum groups are constructed from constant solutions of the Yang--Baxter 
 equation, the $R$-matrices of statistical mechanics usually depend on
 external parameters. For vertex models, these are known as spectral 
parameters,
 while for face models so called  dynamical parameters are present.    
        
The fundamental Faddeev--Reshetikhin--Sklyanin--Takhtajan (FRST) construction 
assigns a 
bialgebra (and in many cases a Hopf algebra) to any constant solution of the  
quantum Yang--Baxter equation. Generalizations of this construction to 
$R$-matrices with spectral parameters lead to Yangians, quantum affine algebras
and Sklyanin algebras, depending on whether the $R$-matrix is a rational, a
trigonometric or an elliptic function. In \cite{fv}, Felder and Varchenko 
gave a similar construction starting from an elliptic $R$-matrix involving 
both  spectral and dynamical parameters. Motivated by this example, Etingof and
Varchenko \cite{ev, ev2} have developped an algebraic framework for studying 
dynamical  $R$-matrices. The resulting ``dynamical quantum
groups'' are not Hopf algebras, but rather Hopf algebroids. 
The appearance of ``oids'' is reflected in
 the 
correspondence between quasiclassical limits of dynamical $R$-matrices and
Poisson structures on Lie groupoids discovered in \cite{ev0}.

In the present paper we study one of the simplest examples of dynamical 
quantum groups, constructed from a trigonometric dynamical $R$-matrix. In
particular, we are interested in the special functions related to its
representation theory. 
It turns out that fundamental objects such as matrix elements and 
Clebsch--Gordan coefficients can be identified with $q$-Racah polynomials 
(or quantum $6j$-symbols) and Askey--Wilson polynomials. This results in new 
algebraic proofs of
 the orthogonality relation and Biedenharn--Elliott
identity satisfied by quantum $6j$-symbols. Moreover, we can interpret the 
orthogonality measure of the  Askey--Wilson polynomials  as a Haar measure on 
the dynamical quantum group. To obtain these results we must extend the 
algebraic machinery 
introduced by Etingof and Varchenko in several ways, with  new definitions 
and new results.
  We  hope that the present case study will be useful when investigating 
Felder's 
elliptic quantum groups \cite{f,fv},
and relating their representation theory to the elliptic hypergeometric series
introduced in \cite{ft}. One would likewise expect connections between higher 
rank dynamical quantum groups and multivariable orthogonal polynomials.   

Let us recall the definition of the quantum dynamical Yang--Baxter (QDYB) 
equation,
 also known as the
 Gervais--Neveu--Felder equation. 
 Let $\gh$ be a finite-dimensional complex vector space, viewed as a 
commutative 
 Lie algebra, and 
$V=\bigoplus_{\al\in\gh^\ast}V_\al$  a
diagonalizable $\gh$-module. In the context of dynamical quantum groups, 
 $\gh$ will typically be a Cartan subalgebra of the corresponding Lie algebra. 
 The QDYB equation may be written as
$$R^{12}(\la-h^{(3)})R^{13}(\la)R^{23}(\la-h^{(1)})=R^{23}(\la)
R^{13}(\la-h^{(2)})R^{12}(\la).$$
This is an identity in the algebra of meromorphic functions $\gh^{\ast}
\rightarrow\End(V\otimes V\otimes V)$.
Here $R:\, \gh^\ast\rightarrow\End(V\otimes V)$ 
is a meromorphic function, $h$ indicates the action of $\gh$, and the upper
indices refer to the factors in the tensor product. For instance, 
$R^{12}(\la-h^{(3)})$ denotes the operator
$$R^{12}(\la-h^{(3)})(u\otimes v\otimes w)=(R(\la-\mu)(u\otimes v))\otimes
w,\quad w\in V_\mu.$$
A dynamical $R$-matrix is by definition a solution of the QDYB equation which 
is $\gh$-invariant, that is, $R:\, \gh^\ast\rightarrow\End_{\gh}(V\otimes V)$. 
For an introduction to the QDYB equation and its relation to other topics we 
refer to \cite{es}.

In the form given above, the QDYB equation first appeared in \cite{gn}. 
Felder \cite{f} pointed out its equivalence to the
star-triangle relation satisfied by the Boltzmann weights of face models.
It is also equivalent to one of the classical identities 
for the $6j$-symbols of quantum mechanics, 
 which reflects the symmetries of the $9j$-symbol \cite{ev2,n}. 
 In this context, the QDYB equation (for $\gh=\mathbb C$) goes back to Wigner's
   1940 paper \cite{w} (cf.~equation (26a) there). 

In the example that we will study $\gh$ is
one-dimensional, and may be viewed as a Cartan subalgebra of
$\mathfrak{sl}(2,\mathbb C)$. We identify $\gh=\gh^\ast=\mathbb C$
and take $V$ to be the two-dimensional $\gh$-module
$V=\mathbb C e_1\oplus\mathbb C e_{-1}$. In the basis
$e_1\otimes e_1$, $e_1\otimes e_{-1}$, $e_{-1}\otimes e_1$, $e_{-1}\otimes
e_{-1}$, the
dynamical $R$-matrix we will consider is given by  
\begin{equation}\label{r}R(\lambda)=\left(\begin{matrix}q&0&0&0\\
0&1&\frac{q^{-1}-q}{q^{2(\lambda+1)}-1}&0\\
0&\frac{q^{-1}-q}{q^{-2(\lambda+1)}-1}&\frac{(q^{2(\lambda+1)}-q^2)
(q^{2(\lambda+1)}-q^{-2})}{(q^{2(\lambda+1)}-1)^2}&0\\
0&0&0&q
 \end{matrix}\right).\end{equation}
 This is the $R$-matrix arising from $6j$-symbols of the quantum
 algebra $\mathcal
 U_q(\mathfrak{sl}(2))$, evaluated in the two-dimensional representation 
 \cite{es}. It can also be interpreted as the
 $R$-matrix for a quasi-Hopf algebra, 
 which can be obtained from $\mathcal U_q(\mathfrak{sl}(2))$ via a Drinfel'd 
twist \cite{b, bbb}. Babelon's construction of the twisting operator uses
 certain ``shifted boundaries" in the quantum algebra. By contrast to the
 interpretation in terms of $6j$-symbols, this construction collapses in
 the limit $q\rightarrow 1$. The shifted boundaries were rediscovered in 
\cite{r}, 
 where they appear as $q$-analogues of group
 elements.
 This gives a link between the QDYB
 equation and harmonic analysis with respect to twisted primitive elements of 
$\mathcal
 U_q(\mathrm{sl}(2))$.
  On the level of special functions, the results of this paper are largely
   parallel to those obtained via twisted primitive elements,
  cf.~\cite{gz,kaw,kf,kv,kz,nm,r}. However, the conceptual connection between 
these two
  approaches remains to be investigated.
 
We will now briefly summarize the contents of the paper. In \S 2 we  review 
the generalized FRST construction from \cite{ev}. We then describe  the 
dynamical quantum group $\frl$ which is obtained from the $R$-matrix \eqref{r}
 through this construction.  In \S 3 we introduce finite-dimensional 
corepresentations of $\frl$. The main result of this section is Theorem 
\ref{lmep}, where the matrix elements of our corepresentations are expressed 
in terms of Askey--Wilson polynomials.  In \S 4 we  introduce a family of 
infinite-dimensional representations of $\frl$, and use them to obtain the 
orthogonality of $q$-Racah polynomials as  discrete ortho\-gonality relations 
for the matrix elements.
In \S 5 and \S 6 we consider tensor product decompositions of 
corepresentations and representations, respectively. In both cases we obtain 
$q$-Racah polynomials as  Clebsch--Gordan coefficients. This gives new 
algebraic  proofs of the   Biedenharn--Elliott identity, which plays a 
fundamental role in quantum mechanics and is also a master identity from the 
viewpoint of special functions.  In \S 7 we show that there is a natural Haar 
functional 
on our algebra, and that it can be identified with the Askey--Wilson measure.

As was mentioned above, we have to extend the algebraic machinery of Etingof
 and Varchenko in several ways. 
Since some readers may   be mainly interested in these parts of the paper, we 
will indicate where they can be found. Our definition of antipode, Definition 
\ref{ad}, differs from the one given in \cite{ev}. With our modified 
definition, we can extend
some basic results for Hopf algebras to the present situation, cf.~Proposition
 \ref{al} and Lemma \ref{sal}. These are proved in Appendix 1.
In \S \ref{ssectSU2} we give a straight-forward definition of $\ast$-structure
 on an $\gh$-algebra. In \S
\ref{scorep} we introduce the concept of corepresentation of an 
$\gh$-bialgebroid. Instead of unitarity of corepresentations, we  speak of 
unitarizability, cf.~Definition \ref{und}. 

To discuss tensor products of corepresentations, we must introduce several 
new algebraic concepts, cf.~\S \ref{sscoalg}--\ref{sstc}. Recall that  if $A$ 
and $B$ are Hopf algebras, then so is $A\ot B$. This is not true for the 
$\gh$-Hopf algebroids which we study:  there is then one  kind of tensor 
product  (denoted $\wt\ot$ and introduced in \cite{ev}) which inherits the 
algebra structure and another one (denoted $\wh\ot$ and  introduced in \S 
\ref{sscoalg}) which  inherits the coalgebroid structure. The 
innocent-looking  Lemma \ref{l23} is a key result which relates these 
structures. Finally, from \S 7 
it is clear what the natural definition of  Haar functional on an 
$\gh$-bialgebroid should be.
      
{\bf Acknowledgement:} The work on this paper was 
mainly done while the second author was employed by
Technische Universiteit Delft.


\section{Preliminaries on dynamical quantum groups}
\subsection{Dynamical $R$-matrices and $\gh$-bialgebroids}
\label{sshb}
In this section we review some of the results of \cite{ev}. First we recall
the fundamental notions of $\gh$-algebra,
$\gh$-bialgebroid and $\gh$-Hopf algebroid. 
These structures are related to the more general Hopf algebroids introduced by 
Lu \cite{lu}. We then recall the generalized
FRST construction, associating an
$\gh$-bialgebroid to any dynamical $R$-matrix.

Throughout this section,  $\gh$ will be a finite-dimensional complex 
commutative Lie algebra and $\mg$ will denote the field of
meromorphic functions on the dual of $\gh$. 

An \emph{$\gh$-algebra} is a complex associative algebra $A$ with $1$, which is
 bigraded over $\gh^\ast$,
$A=\bigoplus_{\al,\be\in\gh^\ast}A_{\al\be},$
and equipped with two algebra embeddings $\mu_l$, $\mu_r:\,
\mg\rightarrow A_{00}$ (the left and right \emph{moment maps}),
such that 
\begin{equation}\label{dr}\mu_l(f)a=a\,\mu_l(T_\al f),
\qquad \mu_r(f)a=a\,\mu_r(T_\be f), \ \ \ a\in A_{\al\be},\ f\in \mg,
\end{equation}
where $T_\al$ denotes the automorphism $T_\al f(\la)=f(\la+\al)$ of $\mg$.
A morphism of $\gh$-algebras is an algebra homomorphism preserving the moment 
maps (and thus also the bigrading).

The \emph{matrix tensor product} $A\widetilde\otimes B$ of two $\gh$-algebras 
is the
$\gh^\ast$-bigraded vector space with
$$(A\widetilde \otimes B)_{\al\ba}=\textstyle\bigoplus_\ga\displaystyle
 (A_{\al\ga}\otimes_{\mg} B_{\ga\be}),$$
where    $\otimes_{\mg}$ denotes the usual tensor product modulo the
relations 
\begin{equation}\label{mtr}\mu_r^A(f)a\otimes b=a\otimes \mu_l^B(f) b,
\ \ \ a\in A,\ b\in B,\ f\in
\mg.\end{equation}
It follows that 
\begin{equation}\label{mtr2}a\mu_r^A(f)\otimes b=a\otimes
b\mu_l^B(f)\end{equation}
in $A\widetilde\ot B$.
The multiplication $(a\otimes b)(c\otimes d)=ac\otimes bd$ and the moment maps 
$$\mu_l^{A\widetilde\otimes B}(f)=\mu_l^A(f)\otimes 1,
\ \ \ \ \ \mu_r^{A\widetilde\otimes
B}(f)=1\otimes\mu_r^B(f)$$
make $A\widetilde\otimes B$ into an $\gh$-algebra.

We denote by $D_\gh$ the algebra of difference operators on $\mg$,
consisting of operators
$$\sum_if_i\,T_{\be_i},\ \ \ f_i\in \mg,\ \be_i\in\gh^\ast.$$
This is an $\gh$-algebra with the bigrading
defined by $f\,T_{-\be}\in (D_\gh)_{\ba\ba}$ and both moment maps equal
to  the natural
embedding. For any $\gh$-algebra $A$, there are canonical $\gh$-algebra
isomorphisms 
$A\simeq A\widetilde\otimes D_\gh\simeq D_\gh\widetilde\otimes A$, defined by
\begin{equation}\label{can}x\simeq x\otimes T_{-\be}\simeq
T_{-\al}\otimes x,\ \ \ x\in A_{\al\be}.\end{equation}
 Thus
the algebra $D_\gh$ plays the role of unit object in the category of
$\gh$-algebras.

An \emph{$\gh$-bialgebroid} is an $\gh$-algebra $A$ equipped with two 
$\gh$-algebra
homomorphisms, $\De:\, A\rightarrow A\widetilde\otimes A$ (the coproduct) and
 $\ep:\,
A\rightarrow D_{\gh}$ (the counit), such that 
$(\Da\otimes\id)\circ\De=(\id\otimes\De)\circ\De$ and, under the 
identifications
(\ref{can}),  $(\ep\otimes\id)\circ\De=(\id\otimes\,\ep)\circ\De=\id$.

We will also need the concept of an \emph{$\gh$-Hopf algebroid}. Our definition
differs slightly from the one given in \cite{ev}.
  
\begin{definition}\label{ad}
An $\gh$-Hopf algebroid is an $\gh$-bialgebroid $A$ equipped with a $\mathbb
C$-linear map $S:\,A\rightarrow A$, called the antipode, such that
\begin{gather}\label{sf}S(\mu_r(f)a)=S(a)\mu_l(f),\ \ \ 
S(a\mu_l(f))=\mu_r(f)S(a),
\ \ \ a\in A,\ f\in \mg,\\
\begin{split}m\circ(\id\otimes\, S)\circ\De(a)&
=\mu_l(\ep(a)1),\ \ \ a\in A,\\
\label{anti}m\circ(S\otimes \id)\circ\De(a)&=\mu_r(T_\al(\ep(a)1)),
\ \ \ a\in A_{\al\be},
\end{split}\end{gather}
where $m$ denotes multiplication and where $\ep(a)1$ is the 
result of applying the difference
operator $\ep(a)$ to the constant function $1\in\mg$.
\end{definition}

The conditions (\ref{sf})  guarantee that 
the left-hand sides of (\ref{anti}) are well-defined; cf.~Lemma \ref {cl} 
in Appendix 1. 
Note that since $\ep(A_{\al\be})=0$ for $\al\neq \be$, the translation 
$T_\al$ in
(\ref{anti}) can be replaced by $T_\be$. In \cite{ev}, the translation $T_\al$
is missing from (\ref{anti}), and (\ref{sf}) is replaced by the stronger
property that $S$ is an algebra antihomo\-morphism which interchanges the 
moment maps. This is a consequence of our definition.

\begin{proposition}\label{al}  
The antipode of an $\gh$-Hopf algebroid is unique. 
Moreover, it satisfies
\begin{gather}\label{c1}S(A_{\al\be})\subseteq A_{-\be,-\al},\\
\label{aprop}S(ab)=S(b)S(a),
\ \ \ \De\circ S=\sigma\circ(S\otimes S)\circ \De,\ \ \ S(1)=1,
\ \ \ \ep\circ S=S^{D_{\gh}}\circ\ep,\\
\label{sm}S(\mu_l(f))=\mu_r(f),\ \ \ S(\mu_r(f))=\mu_l(f),\end{gather}
where $\sigma$ is the flip $\sigma(a\otimes b)=b\otimes a$ and where 
$S^{D_{\gh}}$ is the algebra antihomomorphism of $D_{\gh}$
 defined by $S^{D_{\gh}}(f\,T_\al)=T_{-\al}\circ
f=(T_{-\al}f)T_{-\al}$.  

Moreover, if $A$ is generated as an algebra by a subset $X$, and $S$
is an algebra antihomomorphism of $A$ such that $S(\mu_l(f))=\mu_r(f)$,
$S(\mu_r(f))=\mu_l(f)$ and \emph{(\ref{anti})} holds for every 
$a\in X$ \emph{(}or, more
precisely,  for each component of $a$ with respect to the bigrading\emph{)}, 
then $S$ is an antipode. 
\end{proposition}
 
The proof will be given in Appendix 1. 
It is easy to check that $S^{D_{\gh}}$ is an antipode on $D_{\gh}$, where
$\De^{D_{\gh}}$ is the canonical isomorphism and $\ep^{D_{\gh}}$ the identity
map. If one used the definition of \cite{ev}, there would exist no antipode
on $D_{\gh}$.  Moreover, the last statement of the proposition would be
false. 

We now recall the generalized 
FRST construction. Let 
$V=\bigoplus_{\al\in\gh^\ast}V_\al$ be a finite-dimensional
diagonalizable $\gh$-module 
and $R:\, \gh^\ast \rightarrow \End_\gh(V\otimes V)$ a meromorphic function. 
To each such $R$ one associates an $\gh$-bialgebroid $A_R$. 
Though $R$ is not a priori required to satisfy the QDYB
equation, the construction is motivated by the case when it does.    

Pick a homogeneous basis $\{e_{x}\}_{x\in X}$ of $V$, where $X$ is an index 
set. Write $R_{xy}^{ab}$ for the matrix elements
$$R(\la)(e_a\otimes e_b)=\sum_{xy}R_{xy}^{ab}(\la)\,e_x\otimes e_y$$
of $R$, and define $\om:\, X\rightarrow \gh^\ast$ by $e_x\in V_{\om(x)}$.
The algebra $A_R$ is generated by elements 
$\{L_{xy}\}_{x,\, y\in X}$, 
together with two copies of $\mg$, embedded as subalgebras. 
We will write the elements of these
two copies as $f(\la)$, $f(\mu)$, respectively. The defining relations of $A_R$
are
$$f(\la)L_{xy}=L_{xy}f(\la+\om(x)),\ \ \ f(\mu)L_{xy}=L_{xy}f(\mu+\om(y)),
\ \ \ f(\la)g(\mu)=g(\mu)f(\la)$$
for $f$, $g\in\mg$, together with the $RLL$-relations
\begin{equation}\label{rll}\sum_{xy}R_{ac}^{xy}(\la)L_{xb}L_{yd}=
\sum_{xy}R_{xy}^{bd}(\mu)L_{cy}L_{ax}.\end{equation}
The bigrading on $A_R$ is defined by $L_{xy}\in A_{\om(x),\,\om(y)}$, 
$f(\la),\,
f(\mu)\in A_{00}$, and the moment maps by $\mu_l(f)=f(\la)$,
$\mu_r(f)=f(\mu)$. Note that  $R$  
must be $\gh$-invariant, or equivalently $R_{xy}^{ab}=0$ for
$\om(x)+\om(y)\neq\om(a)+\om(b)$, in order that the $RLL$-relations be
consistent with the grading. Finally one defines a coproduct and a counit on
$A_R$ by
\begin{gather*}\Da(L_{ab})=\sum_{x\in X}L_{ax}\otimes L_{xb},
\ \ \ \ \ \Da(f(\la))=f(\la)\otimes 1,\ \ \ \ \ \Da(f(\mu))=1\otimes f(\mu),\\
\ep(L_{ab})=\de_{ab}\,T_{-\om(a)},
\ \ \ \ \ \ep(f(\la))=\ep(f(\mu))=f.\end{gather*}
These definitions equip $A_R$ with the structure of an $\gh$-bialgebroid.

\subsection{The $\mathrm{SU}(2)$ dynamical quantum group}
\label{ssectSU2}

We will now write down in detail the results of the generalized FRST
construction when applied to
the dynamical $R$-matrix (\ref{r}), where we think of $q$ as a fixed 
number, $0<q<1$. We will denote the corresponding $\gh$-bialgebroid 
$A_R$ by $\frm$. It is
a dynamical analogue of the algebra of polynomials on the space of complex
$2\times 2$-matrices. The four $L$-generators will be denoted by  
$\al=L_{11}$, $\be=L_{1,-1}$,
$\ga=L_{-1,1}$, $\de=L_{-1,-1}$. We also introduce the auxiliary functions 
\begin{align*}F(\la)&=\frac{q^{2(\la+1)}-q^{-2}}{q^{2(\la+1)}-1},\\
G(\la)&=\frac{(q^{2(\la+1)}-q^2)(q^{2(\la+1)}-q^{-2})}{(q^{2(\la+1)}-1)^2},\\
H(\la,\mu)&=\frac{(q-q^{-1})(q^{2(\lambda+\mu+2)}-1)}
{(q^{2(\la+1)}-1)(q^{2(\mu+1)}-1)},\\
I(\la,\mu)&=\frac{(q-q^{-1})(q^{2(\mu+1)}-q^{2(\la+1)})}
{(q^{2(\la+1)}-1)(q^{2(\mu+1)}-1)}.\end{align*}     
The following lemma is useful when checking various
statements made below.

\begin{lemma}\label{hi} The functions $F$, $G$, $H$, $I$ satisfy the following 
relations:
\begin{align*}
q+q^{-1}&=qF(\la)+\frac{q^{-1}}{F(\la-1)}\\
 H(\la,\mu)&=qF(\la)-\frac {q^{-1}}{F(\mu-1)}
=qF(\mu)-\frac {q^{-1}}{F(\la-1)},\\
 I(\la,\mu)&=q(F(\la)-F(\mu))=q^{-1}\left(\frac 1{F(\mu-1)}-\frac
1{F(\la-1)}\right),\\
 G(\la)&=\frac{F(\la)}{F(\la-1)},\\
G(\mu)-G(\la)&=H(\la,\mu)I(\la,\mu).\end{align*}
\end{lemma}

We now write down the definition of $\frm$. As in the non-dynamical 
case, the sixteen $RLL$-relations (\ref{rll}) reduce to 
six independent relations.

\begin{definition}\label{deffrm}
The algebra $\frm$ is generated by the four generators $\al$, $\be$, $\ga$, 
$\de$, together with two copies of $\mg$, whose
elements we write as $f(\la)$, $f(\mu)$. The defining relations are
$$\al\be=qF(\mu-1)\be\al,\ \ \ \al\ga=qF(\la)\ga\al,\ \ \ \be\de=qF(\la)\de\be,
\ \ \ \ga\de=qF(\mu-1)\de\ga,$$
together with any two of the four relations
\begin{equation}\label{cross}\begin{split}
\al\de-\de\al&=H(\la,\mu)\ga\be,\ \ \ \ \ G(\mu)\al\de-G(\la)\de\al
=H(\la,\mu)\be\ga,\\
\be\ga-G(\mu)\ga\be&=I(\la,\mu)\de\al,\ \ \  \ \ \be\ga-G(\la)\ga\be
=I(\la,\mu)\al\de,\end{split}\end{equation}
and, for arbitrary $f,\,g\in \mg$, $f(\la)g(\mu)=g(\mu)f(\la)$,
\begin{equation}\label{f2}\begin{split}
f(\la)\al&=\al f(\la+1),\ \ \ \ \ f(\mu)\al=\al f(\mu+1),\\
f(\la)\be&=\be f(\la+1),\ \ \ \ \ f(\mu)\be=\be f(\mu-1),\\
f(\la)\ga&=\ga f(\la-1),\ \ \ \ \ f(\mu)\ga=\ga f(\mu+1),\\
f(\la)\de&=\de f(\la-1),\ \ \ \ \ f(\mu)\de=\de f(\mu-1).
\end{split}\end{equation}
The bigrading $\frm=\bigoplus_{m,n\in\mathbb Z,\,m+n\in 2\mathbb Z}
\mathcal F_{mn}$ is defined on the
generators by
$$\al\in \mathcal F_{11},\ \ \ \be\in \mathcal F_{1,-1},\ \ \ \ga\in 
\mathcal F_{-1,1},\ \ \ \de\in
\mathcal F_{-1,-1},\ \ \ f(\la),\,f(\mu)\in \mathcal F_{00}.$$
The coproduct $\De:\,\frm\rightarrow \frm\,\widetilde\otimes\,\frm$ and counit
$\ep:\,\frm\rightarrow D_{\gh}$ are algebra homomorphisms defined on the
generators by
\begin{gather*}
\begin{split}\De(\al)&=\al\otimes\al+\be\otimes\ga, \ \ \ \ \
\De(\be)=\al\otimes\be+\be\otimes\de,\\
\De(\ga)&=\ga\otimes\al+\da\otimes\ga,\ \ \ \ \
\De(\de)=\ga\otimes\be+\de\otimes\de,\end{split}\\
\De(f(\la))=f(\la)\otimes 1,\ \ \ \ \ \De(f(\mu))=1\otimes f(\mu),\\
\ep(\al)=T_{-1},\ \ \ \ \ \ep(\be)=\ep(\ga)=0,\ \ \ \ \ \ep(\de)=T_1,
\ \ \ \ \ 
\ep(f(\la))=\ep(f(\mu))=f.\end{gather*}
\end{definition}

That any two of the relations (\ref{cross}) imply the others follows from 
the last identity of Lemma \ref{hi}. By a straight-forward application of the
diamond lemma \cite{be}, any element in $\frm$ can be written uniquely as
a finite sum
$$\sum_{klmn}f_{klmn}(\la,\mu)\al^k\be^l\ga^m\de^n,$$
where $f_{klmn}\in\mg\ot\mg$ (and similarly for any other ordering of the
generators).

Next we describe the dynamical analogue of the determinant.

\begin{lemma}\label{lc}
The element
\begin{equation*}\begin{split}c&=\frac{F(\la)}{F(\mu)}\,\de\al-
\frac{q^{-1}}{F(\mu)}\,\be\ga=\al\de-qF(\la)\ga\be\\
&=\frac{F(\la-1)}{F(\mu-1)}\,\al\de-qF(\la-1)\,\be\ga=\de\al
-\frac{q^{-1}}{F(\mu-1)}\,\ga\be\end{split}\end{equation*}
is a central element of  $\frm$. Moreover, it satisfies
$\De(c)=c\otimes c$, $\ep(c)=1.$
\end{lemma}

The proof is straightforward. That the four expressions are equal follows from
Lemma \ref{hi} and \eqref{cross}. We write down the proof that $\al c=c\al$ 
in detail:
\begin{equation*}\begin{split}\al c&=\al\left(\de\al-\frac{q^{-1}}{F(\mu-1)}
\,\ga\be\right)=\al\de\al-\frac{q^{-1}}{F(\mu-2)}\,\al\ga\be\\
&=\al\de\al-\frac{F(\la)}{F(\mu-2)}\,\ga\al\be=\al\de\al-q\frac{F(\la)}
{F(\mu-2)}\,\ga F(\mu-1)\be\al\\
&=\al\de\al-qF(\la)\ga\be\al=\left(\al\de-qF(\la)\ga\be\right)\al=c\al.
\end{split}\end{equation*}
The element $c$ can also be obtained as a limit case of
the elliptic determinant in \cite{fv}.  

We can now introduce a dynamical analogue of the algebra of functions on the
group $\mathrm{SL}(2,\mathbb C)$.

\begin{definition}\label{deffrsl2}
 The algebra $\mathcal F_R(\mathrm{SL}(2))$ is the $\gh$-Hopf
 algebroid
obtained by adjoining the relation $c=1$ to $\frm$ and defining
the antipode by
\begin{gather*}S(\al)=\frac{F(\la)}{F(\mu)}\,\de,
\ \ \ S(\be)=-\frac{q^{-1}}{F(\mu)}\,\be,
\ \ \ S(\ga)=-qF(\la)\ga,\ \ \ S(\de)=\al,\\
S(f(\la))=f(\mu),\ \ \ S(f(\mu))=f(\la).\end{gather*}
\end{definition}

We need to check that the antipode axioms are satisfied. 
By Proposition \ref{al}, it suffices to show that $S$ reverses the defining
relations of $\frm$ and that (\ref{anti}) holds for the generators. This is a
straight-forward verification.
Note that, for the $L$-generators, (\ref{anti}) can be written compactly as
$$\left(\begin{matrix}S(\al)&S(\be)\\ S(\ga)&S(\de)\end{matrix}\right)
\left(\begin{matrix}\al&\be\\ \ga&\de\end{matrix}\right)=
\left(\begin{matrix}\al&\be\\ \ga&\de\end{matrix}\right)
\left(\begin{matrix}S(\al)&S(\be)\\ S(\ga)&S(\de)\end{matrix}\right)
=\left(\begin{matrix}1&0\\ 0&1\end{matrix}\right).$$
Here the diagonal relations correspond to the four expressions for $c$ given in
Lemma \ref{lc}, and the cross-diagonal relations are defining relations of 
$\frm$.

Another application of the diamond lemma gives the first part of the following
lemma. The second part is proved by a dimension count, cf.~\cite{kk} for the
non-dynamical case.

\begin{lemma}\label{bl}
The elements $\ga^k\be^l\al^m$, $k,\,l,\,m\geq 0$ and 
$\de^k\ga^l\be^m$, $k>0$, $l,\,m\geq 0$,
 form together a basis for $\frl$ as a module over $\mu_l(\mg)\mu_r(\mg)\simeq
 \mg\ot\mg$. The elements
 $(\ga^k\de^l\al^m\be^n)_{k+l+m+n=N}$ are for each $N$ linearly independent
 over $\mu_l(\mg)\mu_r(\mg)$. 
\end{lemma}

Next we will give a  $\ast$-structure to our algebra. In general, to 
introduce a
$\ast$-structure on an $\gh$-algebra $A$, we must assume that a conjugation
(or, equivalently, a real form) $\la\mapsto \bar \la$ has been chosen on 
$\gh^\ast$.  We can then define a $\ast$-structure on 
 $A$ to be a $\mathbb C$-antilinear and antimultiplicative 
involution $a\mapsto a^\ast$ on $A$ such that 
$\mu_l(f)^\ast=\mu_l(\bar f)$ and
$\mu_r(f)^\ast=\mu_r(\bar f)$, where $\bar f(\la)=\overline{f(\bar\la)}$. It
follows that $(A_{\al\be})^\ast=A_{-\bar\al,-\bar\be}$. A $\ast$-structure on 
an
$\gh$-bialgebroid is in addition required to satisfy
$$(\ast\otimes\ast)\circ\De=\De\circ\ast,
\ \ \ \ \ \ep\circ\ast=\ast^{D_{\gh}}\circ\ep,$$
 where $\ast^{D_{\gh}}$ is defined by
$(fT_\al)^\ast=T_{-\bar\al}\circ \bar f=(T_{-\bar\al}\bar f)T_{-\bar\al}$.

\begin{definition} \label{defFRSU2}
 The algebra $\mathcal F_R(\mathrm{SU}(2))$ is the
$\gh$-Hopf algebroid $\mathcal F_R(\mathrm{SL}(2))$ equipped with the
$\ast$-structure $f(\la)^\ast=\bar f(\la)$, $f(\mu)^\ast=\bar f(\mu)$,
$$\al^\ast=\de,\ \ \ \be^\ast=-q\ga,\ \ \ \ga^\ast=-q^{-1}\be,\ \ \
\de^\ast=\al.$$
\end{definition}

It is easy to check that the axioms for a $\ast$-structure are satisfied.
Moreover, since $S$ is invertible we may apply the 
following lemma. We indicate the proof in Appendix~1; in the case at
hand it is easy to verify directly.

\begin{lemma}\label{sal}
Let $A$ be an $\gh$-Hopf algebroid equipped with a $\ast$-structure, such that
the antipode $S$ is invertible. Then $S$ and $\ast$ are related by
$$S\circ\ast\circ S\circ\ast=\id.$$
\end{lemma}

\begin{remark}\label{flr}
When $\la\rightarrow-\infty$, the $R$-matrix \eqref{r} tends
to the $R$-matrix for the standard $\mathrm{SL}(2)$ quantum group.
Accordingly, the Hopf algebra $\mathcal F_q(\mathrm{SL}(2))$
can be obtained as a formal limit of $\frl$ when 
 $\la,\mu\rightarrow-\infty$. We will refer to this limit as the
 non-dynamical case.
 The formal limit of $\frl$
 when $\la,\mu\rightarrow\infty$ is  $\mathcal F_{q^{-1}}(\mathrm{SL}(2))$,
 which can also be viewed as $\mathcal F_q(\mathrm{SL}(2))$ with the opposite
 multiplication. We also need to consider the limit
 $\la\rightarrow -\infty$, $\mu\rightarrow\infty$, in which the defining
 relations of $\frl$ reduce to
 \begin{gather*}\be\al=q\al\be,\ \ \ \al\ga=q\ga\al,\ \ \ \be\de=q\de\be,
 \ \ \ \de\ga=q\ga\de,\\
 \al\de=\de\al,\ \ \ 1=\al\de-q\ga\be=q^2\al\de-q\be\ga.\end{gather*}
 Replacing 
 \begin{equation}\label{usl}(\al,\be,\ga,\de)\mapsto
 (\ga,\al,-q^{-1}\de,-q^{-1}\be)\end{equation}
  we again recover the algebra
 $\mathcal F_q(\mathrm{SL}(2))$.
 We also note that when $q\rightarrow 1$, the
 $R$-matrix \eqref{r} tends to the rational dynamical $R$-matrix
$$R'(\lambda)=\left(\begin{matrix}1&0&0&0\\
0&1&-\frac{1}{\lambda+1}&0\\
0&\frac{1}{\lambda+1}&\frac{\la(\la+2)}{(\la+1)^2}&0\\
0&0&0&1 \end{matrix}\right).$$
The corresponding
$\gh$-Hopf-algebroid $\mathcal F_{R'}(\mathrm{SL}(2))$ can be obtained as
the formal limit of $\mathcal F_R(\mathrm{SL}(2))$ when $q\rightarrow 1$.
This is a ``quantum group without $q$'' which has Racah polynomials 
(classical $6j$-symbols) as matrix elements for its corepresentations.
Yet another interesting limit is the
 one giving rise to a Poisson--Lie groupoid, cf.~\cite{ev0}. 
\end{remark}

\subsection{Some $q$-notation}\label{sec-qnot}

We will follow the standard notation of \cite{gr}, writing
\begin{gather*}\begin{split}(a;q)_k&=
\prod_{j=0}^{k-1}(1-aq^j),\ \ \ k\in\Zp,\\
(a_1,\dots,a_n;q)_k&=(a_1;q)_k\dotsm(a_n;q)_k,\\
\genfrac{[}{]}{0pt}{}{n}{k}_{q}&=\frac{(q;q)_n}{(q;q)_k(q;q)_{n-k}},\\
{}_{r+1}\phi_r\left[\begin{matrix}a_1,\dots,a_{r+1}\\
b_1,\dots,b_r\end{matrix};q,z\right]&=\sum_{k=0}^\infty
\frac{(a_1,\dots,a_{r+1};q)_k}{(q,b_1,\dots,b_{r};q)_k}\,z^k,\\
{}_{r+1}W_r(a;b_1,\dots,b_{r-2};q,z)&={}_{r+1}\phi_r\left[\begin{matrix}
a,q\sqrt{a},-q\sqrt{a},b_1,\dots,b_{r-2}\\
\sqrt{a},-\sqrt{a},aq/b_1,\dots,aq/b_{r-2}\end{matrix};q,z\right]\\
&=\sum_{k=0}^\infty\frac{1-aq^{2k}}{1-a}
\frac{(a,b_{1},\dots,b_{r-2};q)_k}{(q,aq/b_1,\dots,aq/b_{r-2};q)_k}\,z^k.
\end{split}\end{gather*}
Occasionally we will write
\begin{equation}\label{cps}
(a;q)_x=\prod_{j=0}^\infty\frac{1-aq^j}{1-aq^{x+j}},\ \ \  x\in\mathbb
C,\end{equation}
so that
\begin{equation}\label{nps}
(a;q)_{-n}=\frac 1{(aq^{-n};q)_n}.
\end{equation}
We will write
\begin{equation}\label{awm}h_k(\cos\theta,a;q)
=(ae^{i\theta},ae^{-i\theta};q)_k=\prod_{j=0}^{k-1}
\left(1-2aq^j\cos\theta+a^2q^{2j}\right);\end{equation}
this is a polynomial of degree $k$ in $\cos\theta$ which may be called the
 Askey--Wilson monomial.

We will mainly encounter terminating ${}_4\phi_3$- and
${}_8W_7$-series, for which we recall the transformation formulas 
\begin{align}
\nonumber &\quad {}_8 W_7(a;q^{-n},b,c,d,e;q,z)\\
\label{tr1}
&=\frac{(aq,aq/bc,{aq}/{bd},{aq}/{be};q)_n}
{(aq/b,aq/c,aq/d,aq/e;q)_n}\,b^n
\,{}_4\phi_3\left[\begin{matrix}
q^{-n},b,bq^{-n}/a,q/z\\bcq^{-n}/a,bdq^{-n}/a,beq^{-n}/a
\end{matrix};q,q\right],\\
\label{tr2}&=\frac{(aq,b,q/z;q)_n}
{(aq/c,aq/d,aq/e;q)_n}\,\left(-q^{-\frac 12(n+1)}z\right)^n
\,{}_4\phi_3\left[\begin{matrix}
q^{-n},aq/bc,aq/bd,aq/be\\q^{1-n}/b,aq/b,q^{-n}z
\end{matrix};q,q\right],
\end{align}
where $n\in\Zp$ and $z=a^2q^{n+2}/bcde$.  They are obtained by combining
 equations (III.15) and (III.18) in \cite{gr}; note also that the 
 ${}_4\phi_3$'s are obtained from each other by inverting the order of
 summation.

Finally we recall the $q$-Racah and Askey--Wilson polynomials, both 
introduced by Askey and Wilson \cite{askeyw,aw2}, cf. also \cite{gr}. The 
$q$-Racah polynomials are defined by
\begin{equation}\label{defqracah}
 R_n(\mu(x);a,b,c,d;q)= {}_4\phi_3
\left( \begin{array}{c}
q^{-n},abq^{n+1},q^{-x},cdq^{x+1}\\ aq,\ bdq,\ cq \end{array};q,q
\right),
\end{equation}
where it is assumed that one of the  quantities $aq$, $bdq$ or $cq$ equals 
$q^{-N}$ with $N\in\Zp$, $n\leq N$. This is a polynomial of degree $n$ in 
$\mu(x)=q^{-x}+cdq^{x+1}$.
For generic values of the parameters, $\{R_n\}_{n=0}^N$ is a system of 
orthogonal
polynomials, with the orthogonality measure supported on 
$\mu(\{0,1,\dots,N\})$. For later use we recall the symmetries
\begin{equation}\label{rso}R_n(\mu(x);a,b,c,d;q)=R_x(\mu(n);a,dc/a,c,ba/c;q),
\end{equation} 
which is obvious from the definition, and
\begin{align}\nonumber R_n(\mu(x);a,b,q^{-N-1},d;q)&=d^n
\frac{(bq,aq/d;q)_n}{(aq,bdq;q)_n}\\
\label{rs1}&\quad\times R_{N-x}(\mu(n);b,q^{-N-1}/bd,q^{-N-1},q^{N+1}ab;q)\\
\nonumber &=(dq^{x-N})^x\frac{(aq^{1+N-x}/d,bq^{1+N-x};q)_x}{(aq,bdq;q)_x}\\
\label{rs2}&\quad\times R_{N-n}(\mu(x);dq^{-N-1}/a,q^{-N-1}/bd,q^{-N-1},d;q),
\end{align}
which follows from \cite[(III.15)]{gr}.
The Askey--Wilson polynomials are defined by
\begin{equation}\label{defaw}p_n(\cos\theta;a,b,c,d;q)=\frac{(ab,ac,ad;q)_n}
{a^{n}}\, {}_4\phi_3
\left( \begin{array}{c}
q^{-n},abcdq^{n+1},ae^{i\theta},ae^{-i\theta}\\ ab, ac, ad \end{array};q,q
\right);\end{equation}
this is a polynomial of degree $n$ in $\cos\theta$ which is symmetric in the
four parameters $a$, $b$, $c$, $d$.

\section{Corepresentations}

\subsection{Corepresentations of $\gh$-bialgebroids}
\label{scorep}
In this section we will discuss 
corepresentations of  $\gh$-bialgebroids 
(or better $\gh$-coalgebroids, cf.~\S \ref{sscoalg}). 
As in the case of representations (cf.~\S \ref{ssecdynreps}), it
is natural to view corepresentation spaces as ``dynamical'' spaces. 
In \cite{ev} a category of so called $\gh$-vector spaces  is introduced, whose
objects are complex vector spaces but whose morphisms are $\mathbb C$-linear
maps $V\rightarrow
\mg\otimes W$. We choose to work instead with vector spaces over $\mg$, with
$\mg$-linear maps as morphisms. We must point out, however, that this is purely
a matter of taste, and that we could equivalently have used the
category introduced in \cite{ev}.
 
Thus we define an \emph{$\gh$-space} to be a vector space over $\mg$, 
which is also a diagonalizable $\gh$-module,
$V=\bigoplus_{\al\in\gh^\ast}V_\al$, with $\mg V_\al\subseteq V_\al$ 
for all $\al$. A morphism of $\gh$-spaces is an $\gh$-invariant
(that is, grade-preserving) $\mg$-linear map.

 If $A$ is an $\gh$-algebra
 and $V$ an $\gh$-space, we define
 $A\widetilde\otimes V=\bigoplus_{\al\be} A_{\al\be}\otimes_{\mg}V_\be,$
   where
 $\otimes_{\mg}$ denotes the usual tensor product modulo the relations
 \begin{equation}\label{awt}\mu_r^A(f)a\otimes v=a\otimes fv.\end{equation}  
 The grading
  $A_{\al\be}\otimes_{\mg} V_\be\subseteq(A\widetilde\otimes V)_\al$ and the
 extension of scalars 
 $f(a\otimes v)=\mu_l^A(f)a\otimes v$
  make $A\widetilde\otimes V$ into
 an $\gh$-space.  This definition is compatible with the matrix tensor
 product of $\gh$-algebras in the sense that $(A\widetilde\otimes B)
\widetilde\otimes
 V=A\widetilde\otimes(B\widetilde\otimes V)$ when $A$ and $B$ are 
$\gh$-algebras and $V$ an $\gh$-space.     

We can now define a (left) \emph{corepresentation} of an $\gh$-bialgebroid $A$
 on an  $\gh$-space
$V$ to be an $\gh$-space morphism $\pi:\,V\rightarrow A\widetilde\otimes V$ 
such that
\begin{equation}\label{cor}(\De\otimes\id)\circ\pi=(\id\otimes\,\pi)\circ\pi,
\ \ \ \ \ (\ep\otimes\id)\circ\pi=\id.\end{equation}
The first of these equalities is in the sense of the natural isomorphism
$(A\widetilde\otimes A)\widetilde\otimes V\simeq
A\widetilde\otimes(A\widetilde\otimes V)$, the
second one in terms of the isomorphism $D_{\gh}\widetilde\otimes V\simeq V$ 
defined
by $f\,T_{-\al}\otimes v\simeq fv$, $f\in\mg$, $v\in V_\al$. A morphism or
\emph{intertwiner} of corepresentations is an $\gh$-space morphism
$\phi:\,V_1\rightarrow V_2$ such that 
\begin{equation}\label{int}\pi_{2}\circ\phi=(\id\ot\,\phi)\circ \pi_{1},
\end{equation}
where $\pi_{i}:\,V_i\rightarrow A\wt\ot V_i$ are corepresentations (note that
$\id\ot\,\phi$ factors to
a map on $A\widetilde\ot V_1$). 

If we pick a homogeneous
basis $\{v_k\}_k$ of  $V$ (over $\mg$), 
$v_k\in V_{\om(k)}$, and introduce the matrix elements
$t_{kj}\in A$ by
\begin{equation}\label{mei}\pi(v_k)=\sum_j t_{kj}\otimes v_j,\end{equation}
which is possible in view of (\ref{awt}),
then (\ref{cor}) may be stated as
\begin{equation}\label{mer}\De(t_{kl})=\sum_j  t_{kj}\otimes t_{jl},\ \ \ \ \
\ep(t_{kl})=\de_{kl}\,T_{-\om(k)}.\end{equation}
This refers only to the complex vector space spanned by the chosen basis.
 Thus, as long as we consider a single corepresentation, the
 dynamical variables play no role. However,  if $V_1$ and $V_2$ are two 
 corepresentations with matrix
elements $t_{kj}^{1}$, $t_{kj}^2$, with respect to some bases $\{v_k^1\}_k$, 
$\{v_k^2\}_k$, and
$\phi:\,V_1\rightarrow V_2$ an intertwiner with matrix
elements $\phi_{kj}\in\mg$ defined by
\begin{equation}\label{ime}\phi(v_k^1)=\sum_j \phi_{kj} v_j^2,\end{equation}
then \eqref{int} may be written as
\begin{equation}\label{int2}\sum_j \phi_{jl}(\mu)t_{kj}^1
=\sum_j\phi_{kj}(\la)t_{jl}^2
\ \ \ \text{for all}\ k,\,l,\end{equation}
so the dynamical variables appear when considering intertwiners. 
For later use we observe that
if $A$ is an $\gh$-Hopf algebroid, then it follows from \eqref{anti} and 
\eqref{mer} that
\begin{equation}\label{mes}\de_{kl}=\sum_j  S(t_{kj})t_{jl}=
\sum_j  t_{kj}S(t_{jl}).\end{equation}

If $A$ is an $\gh$-bialgebroid, viewed as an
$\gh$-space with $A_\al=\bigoplus_\be A_{\al\be}$, $fv=\mu_l(f)v$, then the
coproduct $\De$ defines a corepresentation of $A$ on itself, the \emph{regular
corepresentation}. More generally, if $V$ is a subspace of $A$ with
$\De(V)\subseteq A\widetilde\otimes V$ and $\mu_l(\mg)V\subseteq V$,  then
$\De\big|_{V}$ defines a corepresentation of $A$ on
$V$.

\subsection{Corepresentations of $\frl$}

Let $V_N$ 
be the  subspace of $\frm$ span\-ned by
$\{\ga^{N-k}\al^k\}_{k=0}^N$ together with $\mu_l(\mg)$. 
It is easy to see, and will be made clear below,
that $\De(V_N)\subseteq \frm\,\widetilde\otimes V_N$, so that
   $\De\big|_{V_N}$ is a
corepresentation. In this section we will compute the matrix elements 
of these corepresentations. Our method follows that of \cite{ko} for the
non-dynamical case. 
The following lemma will be used 
repeatedly.

\begin{lemma}\label{cr}
The following relations hold in the algebra $\frm$:
\begin{align*}
\al^n\be^m&=q^{mn}\frac{(q^{-2(\mu+m)};q^2)_n}{(q^{-2\mu};q^2)_n}
\,\be^m\al^n,\\
\al^n\ga^m&=q^{mn}\frac{(q^{-2(\la+m+1)};q^2)_n}{(q^{-2(\la+1)};q^2)_n}
\,\ga^m\al^n,\\
\be^n\de^m&=q^{mn}\frac{(q^{-2(\la+m+1)};q^2)_n}{(q^{-2(\la+1)};q^2)_n}
\,\de^m\be^n,\\
\ga^n\de^m&=q^{mn}\frac{(q^{-2(\mu+m)};q^2)_n}{(q^{-2\mu};q^2)_n}
\,\de^m\ga^n.\\
\end{align*}
\end{lemma}

\begin{proof}
We prove the second relation. The other ones are derived similarly or by
observing that the four subalgebras generated by $\{\al,\be\}$, $\{\al,\ga\}$,
$\{\be, \de\}$ and $\{\ga,\de\}$ are all isomorphic.  
It is clear from the defining relations that
$$\al^n\ga^m=C_{mn}(\la)\ga^m\al^n$$
for some $C_{mn}\in\mg$.
Multiplying with $\al$ from the left gives
$$\al^{n+1}\ga^m=
\al\, C_{mn}(\la)\ga^m\al^n=C_{mn}(\la-1)\al\ga^m\al^n
=C_{mn}(\la-1)C_{m1}(\la)\ga^m\al^{n+1},$$
leading to the recursion relation
$$C_{m,n+1}(\la)=C_{mn}(\la-1)C_{m1}(\la),$$
while multiplying with $\ga$ from the right similarly leads to
$$C_{m+1,n}(\la)=C_{mn}(\la)C_{1n}(\la+m).$$
The coefficients $C_{mn}$ are determined by these two recursion relations 
together with the
initial conditions $C_{m0}=C_{0n}=1$, $C_{11}=qF$. It is easy to check that the
solution is indeed given by 
$C_{mn}(\la)=q^{mn}(q^{-2(\la+m+1)};q^2)_n/(q^{-2(\la+1)};q^2)_n$.
\end{proof}

We can now compute the matrix elements of our corepresentations. In what 
follows
it will be convenient to write 
\begin{equation}\label{ro}f(\la)=f(\la)\ot 1,\ \ \ \ \ f(\ro)=f(\mu)\ot 1
=1\ot f(\la) ,\ \ \ \ \ f(\mu)=1\ot f(\mu)\end{equation}
for the three dynamical variables present in a tensor product $A\widetilde\ot
A$.

\begin{proposition}\label{mep}
In the algebra $\frm$,
\begin{equation}\label{med}\De(\ga^{N-k}\al^k)=\sum_{j=0}^N t_{kj}^N\otimes 
\ga^{N-j}\al^j,\end{equation}
where the matrix elements $t_{kj}^N$ are given by
\begin{equation*}\begin{split}t_{kj}^N&=\sum_{l=\max(0,\,j+k-N)}^{\min(j,\,k)}
\qb{N-k}{j-l}\qb{k}{l}q^{j(j+2k-N)+l(3l-3k-3j+N)}\\
&\quad\times\frac{(q^{2(j-N-\mu-1)};q^2)_{j-l}}
{(q^{2(j+k-l-N-\mu-1)};q^2)_{j-l}}\,\ga^{j-l}\de^{N-k-j+l}\al^l\be^{k-l}.
\end{split}\end{equation*}
\end{proposition}

Note that $\ga^{N-k}\al^k\in \frm_{2k-N,N}$, which implies that
 $t_{kj}^N\in \frm_{2k-N,2j-N}$. This is of course in agreement with the 
 proposition. 

\begin{proof}
We will first prove the case $k=N$, which may be written as 
\begin{equation}\label{dak}\De(\al^k)=\sum_{l=0}^k\qb k l 
q^{l(l-k)}\al^l\be^{k-l}\otimes\ga^{k-l}\al^l.\end{equation}
It is clear from the defining relations that
$$\De(\al^k)=(\al\otimes\al+\be\otimes\ga)^k=\sum_{l=0}^k
C_{kl}(\ro)\al^l\be^{k-l}\otimes\ga^{k-l}\al^l$$
for some coefficients $C_{kl}\in\mg$.
To find a recursion formula for $C_{kl}$ we write
\begin{equation*}\begin{split} 
\De(\al^{k+1})&=(\al\ot\al+\be\ot\ga) \sum_{l}
C_{kl}(\ro)\al^l\be^{k-l}\otimes\ga^{k-l}\al^l\\
&=\sum_l C_{kl}(\ro-1)\al^{l+1}\be^{k-l}\otimes\al\ga^{k-l}\al^l
+C_{kl}(\ro+1)\be\al^l\be^{k-l}\otimes\ga^{k-l+1}\al^l\\
&=\sum_l C_{kl}(\ro-1)q^{k-l}\frac{1-q^{-2(\ro+k-l+1)}}{1-q^{-2(\ro+1)}}
\,\al^{l+1}\be^{k-l}\otimes\ga^{k-l}\al^{l+1}\\
&\quad+C_{kl}(\ro+1)q^{-l}\frac{1-q^{-2(\ro-l+1)}}{1-q^{-2(\ro+1)}}
\,\al^l\be^{k-l+1}\otimes\ga^{k-l+1}\al^l,\\
\end{split}\end{equation*} 
where we used Lemma \ref{cr} in the last step. This leads to the recursion
$$C_{k+1,l}(\ro)=C_{k,l-1}(\ro-1)q^{k-l+1}
\frac{1-q^{-2(\ro+k-l)}}{1-q^{-2(\ro+1)}}
+C_{kl}(\ro+1)q^{-l}\frac{1-q^{-2(\ro-l+1)}}{1-q^{-2(\ro+1)}}$$
for $l=0,\dots,k+1$,
where $C_{k,-1}=C_{k,k+1}=0$. We must check that this holds for the constant
functions $C_{kl}(\ro)=\qb k l q^{l(l-k)}$. After simplifications one arrives 
at
\begin{equation}\label{dp}\begin{split}
\left(1-q^{-2(\ro+1)}\right)\qb{k+1}{l}&=\left(q^{2(k-l+1)}-q^{-2(\ro+1)}
\right)
\qb{k}{l-1}\\
&\quad +\left(1-q^{2l}q^{-2(\ro+1)}\right)\qb{k}{l},\end{split}\end{equation}
which is equivalent to the two different Pascal's triangle identities for the
$q$-binomial coefficients \cite{gr}:
\begin{align*}
\qb{k+1}{l}&=q^{2(k-l+1)}\qb{k}{l-1}+\qb{k}{l},\\
\qb{k+1}{l}&=\qb{k}{l-1}+q^{2l}\qb{k}{l}.
\end{align*}

The case $k=0$ of the proposition, which may be  written as
$$\De(\ga^{N-k})=\sum_{m=0}^{N-k}\qb {N-k}{m} q^{m(m-N+k)}\ga^m\de^{N-k-m}
\otimes\ga^{N-k-m}\al^m,$$
may be proved in the same way. Multiplying this
expression and (\ref{dak}) gives
\begin{equation*}\begin{split}\De(\ga^{N-k}\al^k)&=\sum_{m=0}^{N-k}\sum_{l=0}^k
\qb{N-k}{m}\qb{k}{l}q^{m(m-N+k)+l(l-k)}\\
&\quad\times\ga^m\de^{N-k-m}\al^l\be^{k-l}\otimes
\ga^{N-k-m}\al^m\ga^{k-l}\al^l\\
&=\sum_{m=0}^{N-k}\sum_{l=0}^k
\qb{N-k}{m}\qb{k}{l}q^{m(m-N+k)+l(l-k)+m(k-l)}\\
&\quad\times\frac{(q^{-2(\ro+N-l-m+1)};q^2)_m}
{(q^{-2(\ro+N-k-m+1)};q^2)_m}\,\ga^m\de^{N-k-m}\al^l\be^{k-l}\otimes
\ga^{N-l-m}\al^{l+m}\end{split}\end{equation*}
by Lemma \ref{cr}.
Putting $m=j-l$ and simplifying completes the proof.
\end{proof}

We will now consider $V_N$ 
as a corepresentation space of the quotient algebra $\frl$ of $\frm$. 
Using the determinant relation we will factor each matrix element
as a trivial part times a function involving only  commuting variables. 
The following lemma is the key to finding these factorizations.

\begin{lemma}\label{lx} The element $\Xi$ defined by
\begin{equation}\label{xi}\begin{split}\Xi&=q^{\la-\mu+1}+q^{\mu-\la-1}
-q^{-(\la+\mu+2)}(1-q^{2(\la+2)})(1-q^{2\mu})\ga\be\\
&=q^{\la-\mu-1}+q^{\mu-\la+1}
-q^{-(\la+\mu+2)}(1-q^{2(\la+1)})(1-q^{2(\mu+1)})\be\ga\end{split}
\end{equation}
is a central
element of $\frl$. Moreover, it satisfies $\ep(\Xi)=0$, $S(\Xi)=\Xi$ and
$\Xi^\ast=\Xi$ with respect to the $\fr$ $\ast$-structure.
\end{lemma}

Note that the algebra $\mathcal F_q(\mathrm{SL}(2))$ has trivial center.
Therefore, the existence of $\Xi$ is a purely dynamical phenomenon.
In  \S \ref{tpdr} we will see that $\Xi$ plays the role of Casimir
operator in the representation theory of $\fr$. 

The proof of Lemma \ref{lx} is straight-forward.  
 That the two expressions for $\Xi$ are equal follows from the
determinant relation. To prove centrality, we first check that
$\Xi^\ast=\Xi$. Then it is enough to prove that $\Xi$ commutes with $\al$,
$\be$,
and $f(\la)$, $f(\mu)$ for $f\in\mg$.
 To check that $\Xi$ commutes with $\be$  
one should write $\be\Xi$ using the first and $\Xi\be$ using the second
expression in \eqref{xi}. 

We will need the following relations 
in the subalgebra $\frl_{00}$ of $\frl$. 
This is a commutative algebra generated by $1$, $\Xi$, $\mu_l(\mg)$ and
$\mu_r(\mg)$. 

\begin{lemma}\label{ccr}
In the algebra $\frl_{00}$, the following relations hold:
\begin{align*}
\al^k\de^k&=\frac{1}
{(q^{-2(\la+1)},q^{-2\mu};q^2)_k}\,h_k(\textstyle\frac 12\displaystyle\,
\Xi,q^{-(\la+\mu+1)};q^2),\\
\de^k\al^k&=\frac{1}
{(q^{2(\la+2)},q^{2(\mu+1)};q^2)_k}\,h_k(\textstyle\frac 12\displaystyle\,
\Xi,q^{\la+\mu+3};q^2),\\
\be^k\ga^k&=\frac{(-q^{-1})^k}
{(q^{-2(\la+1)},q^{2(\mu+1)};q^2)_k}\,h_k(\textstyle\frac 12\displaystyle\,
\Xi,q^{\mu-\la+1};q^2),\\
\ga^k\be^k&=\frac{(-q)^k}
{(q^{2(\la+2)},q^{-2\mu};q^2)_k}\,h_k(\textstyle\frac 12\displaystyle\,
\Xi;q^{\la-\mu+1};q^2),
\end{align*}
where we use the notation \eqref{awm}.
\end{lemma}

 To prove Lemma \ref{ccr}, one
checks that
the case $k=1$ follows from the determinant relation and the definition of
$\Xi$. The general case then follows immediately by induction on $k$, using 
that $\Xi$ is central.

We can now prove the main result of this section.

\begin{theorem}\label{lmep}
In the algebra $\frl$, the matrix elements $t_{kj}^N$ are given by
\begin{equation*}t_{kj}^N=\begin{cases}
\ga^{j-k}\de^{N-k-j}P_k,\ \ \ &k\leq j,\ k+j\leq N,\\
\ga^{j-k}P_{N-j}\,\al^{k+j-N},\  \ \ &k\leq j,\ N\leq k+j,\\
\de^{N-k-j}P_{j}\,\be^{k-j},\ \ \ &j\leq k,\ k+j\leq N,\\
P_{N-k}\,\al^{k+j-N}\be^{k-j},\ \ \ &j\leq k,\ N\leq k+j,
\end{cases}\end{equation*}
where $P_n\in\frl_{00}$ can be written in terms of  Askey--Wilson polynomials, 
cf.~\eqref{defaw}, as
\begin{equation*}\begin{split}
P_n&=q^{n(\la-\mu+1+n+|N-k-j|)+j(j-N)}\frac{(q^2;q^2)_{N-n}}{(q^2;q^2)_{j}
(q^2;q^2)_{N-j}(q^{2(\la+2)},q^{2(-\mu+
|N-k-j|)};q^2)_n}\\
&\quad\times p_n(\textstyle\frac 12\,\displaystyle\Xi;q^{\mu-\la+1+2|j-k|},
q^{\la-\mu+1},q^{\la+\mu+3},q^{-\la-\mu-1+2|N-k-j|};q^2).
\end{split}\end{equation*} 
\end{theorem}

In \S \ref{hfs}  we will see that   the 
Schur-type orthogonality relations for matrix
elements correspond to the orthogonality of  Askey--Wilson polynomials. 
We also remark that in the alternative notation of \cite{nm},
\begin{equation}\label{nm}p_n^{(\al,\be)}(x;s,t|q)
=p_n(x;q^{\frac 12}t/s,q^{\frac 12+\al}s/t,
-q^{\frac 12}/st,-q^{\frac12+\be}st;q),\end{equation}
the Askey--Wilson polynomial of Theorem \ref{lmep} may be suggestively
 written as 
\begin{equation*} p_n^{(|j-k|,|N -k-j|)}
(\textstyle\frac 12\,\displaystyle\Xi;iq^{-\la-1},iq^{-\mu-1};q^2).
\end{equation*}

After the preparations
that we have made, the proof of Theorem \ref{lmep} is straightforward. 
For instance, in
the case $k\leq j$, $k+j\leq N$, Proposition \ref{mep} gives
$$t_{kj}^N=\sum_{l=0}^k C_l(\mu)\ga^{j-l}\de^{N-k-j+l}\al^l\be^{k-l},$$ 
where
$$C_l(\mu)=\qb{N-k}{j-l}\qb{k}{l}q^{j(j+2k-N)+l(3l-3k-3j+N)}
\frac{(q^{2(j-N-\mu-1)};q^2)_{j-l}}
{(q^{2(j+k-l-N-\mu-1)};q^2)_{j-l}}.$$
Using Lemma \ref{cr} repeatedly one can rewrite this expression as
$$t_{kj}^N=\ga^{j-k}\de^{N-k-j}\sum_{l=0}^k D_l(\la,\mu)
\ga^{k-l}\be^{k-l}\de^l\al^l,$$
where one computes
\begin{multline*}D_l(\la,\mu)=\qb{N-k}{j-l}\qb{k}{l}
q^{j(j-N)+k(N-k)+jk+l(3l-2k-2j)}\\
\times\frac{(q^{2(-\mu-j-1)};q^2)_{j-l}(q^{-2\mu};q^2)_{k-l}
(q^{2(-\mu-l)};q^2)_{l}
(q^{2(-\la+l-k-1)};q^2)_{k-l}}
{(q^{2(-\mu+k-l-j-1)};q^2)_{j-l}(q^{2(-\mu+N-k-j)};q^2)_{k-l}
(q^{2(-\mu+k-2l)};q^2)_{l}(q^{2(-\la-1-k)};q^2)_{k-l}}.\end{multline*}
Plugging in the expressions from Lemma \ref{ccr} and using
elementary identities for $q$-shifted factorials gives an expression like the 
one we are looking for, but with
\begin{equation*}\begin{split}
P_k&=\qb{N-k}{j}q^{k(2\mu+2+3j-N)+j(j-N)}
\frac{(q^{-2(\mu+1+j)},q^{\la-\mu+1}\xi,q^{\la-\mu+1}\xi^{-1};q^2)_{k}}
{(q^{-2(\mu+1)},q^{2(\la+2)},q^{2(\mu+1-N+j)};q^2)_{k}}\\
&\quad\times{}_8W_7\left(q^{2(\mu+1-k)};q^{-2k},q^{-2j},q^{2(\mu+1-N+j)},
q^{\la+\mu+3}\xi,q^{\la+\mu+3}\xi^{-1};q^2,q^{2(N-k-\la)}\right),
\end{split}\end{equation*}
where $\xi+\xi^{-1}=\Xi$ (here $\xi$ is  a formal quantity used to write the 
above expression in standard $q$-notation).
Applying  \eqref{tr2} 
gives the desired expression for $P_k$. 
In the remaining three cases, the theorem can be proved similarly, or derived
 using symmetries of the matrix elements, cf.~Remark \ref{sr} below.

We conclude this section by showing that analogues of the Peter--Weyl theorem 
and Schur's Lemma hold for $\frl$.

\begin{proposition}\label{pw}
The matrix elements $\{t_{kj}^N\}$ for $k,\,j,\,N\in\Zp$, $k,\,j\leq N$ form a
basis for $\frl$ as a module over
 $\mu_l(\mg)\mu_r(\mg)$.
\end{proposition}

\begin{proof}
Let $I$ denote the set of invertible elements in
$\mu_l(\mg)\mu_r(\mg)$. First we observe that the element $P_n$ of 
 Theorem \ref{lmep} is a polynomial 
in $\Xi$ of degree $n$ with the leading
coefficient in $I$. By 
Lemma \ref{ccr}, $\ga^k\be^k$ is a polynomial in $\Xi$ of degree $k$, with the
leading coefficient in $I$. Applying Gauss elimination, we can expand
$\Xi^k=\sum_l c_{l}\ga^l\be^l$, again with the leading coefficient $c_{k}\in
I$. Combining these facts and using Lemma \ref{cr}, we can for $k+j\leq N$
write $t_{kj}^N=\sum_{l=0}^{\min(j,k)}d_l\,\de^{N-k-j}\ga^{j-l}\be^{k-l}$
with $d_0\in I$. Again by Gauss elimination,  each element 
$\de^{l}\ga^m\be^n$ is in the $\mu_l(\mg)\mu_r(\mg)$-span  of 
 $\{t_{kj}^N\}_{k+j\leq N}$. Similarly, $\ga^l\be^m\al^n$ is in the 
 $\mu_l(\mg)\mu_r(\mg)$-span
of $\{t_{kj}^N\}_{k+j\geq N}$. Thus, by
 Lemma \ref{bl},   the matrix elements span $\frl$.  A dimension count 
completes the proof.
\end{proof}

\begin{remark}\label{lrd}
It is easy to check
that, for any $0\neq g\in\mg$ and $0\leq j\leq N$, the coproduct restricted to 
$\bigoplus_{k=0}^N \mu_r(g)\mu_l(\mg)t_{kj}^N$ is a corepresentation
equivalent to $V_N$,  
$\ga^{N-k}\al^k\mapsto
g(\mu)t_{kj}^N$ being  an intertwiner. 
Then Proposition \ref{pw} gives the decomposition
$$\frl\simeq \bigoplus_{N=0}^\infty  V_N\ot\mg^{N+1}$$
of the regular corepresentation. In contrast to the non-dynamical case,
 each $V_N$ occurs with infinite multiplicity.
\end{remark}

\begin{corollary}
The corepresentations $V_N$ are irreducible, that is, if $U\subseteq
V_N$ is an $\gh$-subspace which is invariant in the sense that
 $\De(U)\subseteq A\wt\ot U$, then $U=0$ or $U=V_N$.
\end{corollary}

\begin{proof}
That $U$ is an $\gh$-subspace means precisely that
$U=\sum_{k\in\La}\mu_l(\mg) \ga^{N-k}\al^k$ for some 
$\La\subseteq \{0,\dots, N\}$.
Then the invariance of $U$ means that $t_{kj}^N=0$ for $k\in\La$, $j\notin\La$.
By Proposition \ref{pw}, this is impossible unless $U=0$ or $U=V_N$.
\end{proof}

\begin{corollary}\label{sl}
If $\phi:V_M\rightarrow V_N$ is an intertwining map, then $\phi=0$ for 
$M\neq N$ and $\phi=Z\id$ for some $Z\in\mathbb C$ otherwise. 
\end{corollary}

\begin{proof}
First one checks that the kernel and the image of an
intertwiner are always invariant $\gh$-subspaces. It follows that an
intertwiner between irreducible corepresentations is either zero or bijective. 
Counting dimensions (over $\mg$) then gives the first statement.

For the second statement, we assume that $\phi:\,V_M\rightarrow V_M$. 
Since $\phi$ preserves the grading we have 
$\phi(\ga^{M-k}\al^k)=\phi_k \ga^{M-k}\al^k$ with $\phi_k\in\mg$. 
By \eqref{int2}, the intertwining property
means that
$\phi_l(\mu)t_{kl}^M=\phi_k(\la)t_{kl}^M$
for all $k$, $l$. By Proposition \ref{pw}, this implies that $\phi_k(\la)$ is
independent of $k$ and $\la$, which completes the proof.
\end{proof}

More generally,  if $V$ is a corepresentation of an $\gh$-bialgebroid 
such that its matrix elements are linearly independent over 
$\mu_l(\mg)\mu_r(\mg)$, then
it follows from \eqref{int2} that any intertwiner $V\rightarrow V$ is a 
complex  multiple of the identity.
We expect this to be true for all irreducible corepresentations  
of a large class of interesting dynamical quantum groups. 
The following example  shows that it is not true for general
$\gh$-Hopf algebroids. This suggests that one might adopt the Peter--Weyl
theorem as an axiom for dynamical quantum groups, which would give an approach
similar to that of \cite{dk} in the non-dynamical case.

\begin{example}
This example is modelled on the group of rotations of the plane,
to which it formally reduces when $\la=-1$.
Let $A=\mg[x,y]$ be the algebra of polynomials in two commuting variables over
the field $\mg$, $\gh^\ast=\mathbb C$. Define the coproduct and counit by
\begin{gather*}\De(x)=x\ot x+\la y\ot y,\ \ \ \ \ \De(y)=y\ot x+x\ot y,
\ \ \ \ \ \ep(x)=1,\ \ \ \  \ \ep(y)=0.\end{gather*}
The moment maps $\mu_l(f)=\mu_r(f)=f$ and the bigrading $A=A_{00}$ give $A$ the
structure of an $\gh$-bialgebroid. 
Let $V$ be the two-dimensional subspace consisting of homogeneous polynomials 
of degree $1$. Then $\De\big|_V$ is a corepresentation of $A$.
One easily checks that if $fx+gy$ spans a
one-dimensional invariant subspace, then $g(\la)^2=\la f(\la)^2$,
which has no meromorphic solutions. Therefore, $V$ is irreducible. On the other
hand,  the intertwiners from $V$ to $V$ are given by
$$C(x)=fx+\la g y,\ \ \ \ \ C(y)=gx+fy$$
for $f$, $g\in\mg$ arbitrary.  To get a similar 
example for $\gh$-Hopf algebroids,
adjoin the relation $x^2-\la y^2=1$ and
define the antipode by $S(x)=x$, $S(y)=-y$. 
\end{example}

\subsection{Unitarity of the corepresentations}
\label{uns}

Our next task is to show that our co\-representations are, in a certain sense,
unitarizable. 

\begin{definition}\label{und} 
A corepresentation of a $\ast$-$\gh$-Hopf
algebroid $A$ on an $\gh$-space $V$ is \emph{unitarizable} 
if there exists a basis of $V$ such that the corresponding 
 matrix elements, cf.~\eqref{mei}, satisfy
$$\Ga_k(\mu)S(t_{kj})^\ast=\Ga_j(\la)t_{jk}$$
for some $0\neq \Ga_k\in\mg$ with $\bar \Ga_k=\Ga_k$. 
\end{definition}

We will call the functions $\Ga_k$ \emph{normalizing functions} for $V$, with
respect to the given basis $\{v_k\}_k$.
If we formally introduce 
normalized basis vectors $e_k$ and matrix elements $\tilde
t_{kj}$ by
\begin{gather*}e_k=\sqrt{\Ga_k}\,v_k,\\
 \De(e_k)=\sum_{j} \tilde t_{kj}\otimes e_j,\end{gather*}
then, computing formally, 
$$\tilde t_{kj}=\sqrt{\frac{\Ga_k(\la)}{\Ga_j(\mu)}} 
\,t_{kj}$$
and the unitarizability criterion can be stated as the unitarity 
$$S(\tilde t_{kj})^\ast=\tilde t_{jk}.$$

\begin{proposition}\label{un}When considered as elements of $\fr$, the
matrix elements $t_{kj}^N$ satisfy
$$\qb{N}{k}\frac{(q^{2(k-N-\mu-1)};q^2)_k}{(q^{-2\mu};q^2)_k}\,
S(t_{kj}^N)^\ast=\qb{N}{j}\frac{(q^{2(j-N-\la-1)};q^2)_j}
{(q^{-2\la};q^2)_j}\,t_{jk}^N.$$
In particular, $V_N$ is a unitarizable corepresentation of $\fr$.
\end{proposition}

To prove Proposition \ref{un}, we first note that 
$$\si\circ((\ast\circ
S)\ot(\ast\circ S))\circ\De=\De\circ\ast\circ S,$$
where, as above, $\si(a\ot b)=b\ot a$. Thus, applying 
$\si\circ((\ast\circ S)\ot(\ast\circ S))$ to (\ref{med}) gives 
\begin{equation}\label{e1}\De(S(\ga^{N-k}\al^k)^\ast)
=\sum_{j=0}^N S(\ga^{N-j}\al^j)^\ast \ot S(t_{kj}^N)^\ast.\end{equation}
On the other hand, one has the identity
$$\De(t_{Nk}^N)=\sum_{j=0}^N t_{Nj}^N\ot t_{jk}^N,$$
which is a special case of (\ref{mer}).  By (\ref{dak}), it
may be written as
\begin{equation}\label{e2}\qb{N}{k}q^{k(k-N)}\De(\al^k\be^{N-k})=
\sum_{j=0}^N\qb{N}{j}q^{j(j-N)}\al^j\be^{N-j}
\ot t_{jk}^N.\end{equation}
 Comparing (\ref{e1}) and (\ref{e2}), we see that we can relate $t_{jk}^N$ and
 $S(t_{kj}^N)^\ast$ using the following lemma.
 
  \begin{lemma}\label{unl}
In the algebra $\fr$,
$$S(\ga^{N-k}\al^k)^\ast=C_k^N(\la,\mu)\al^k\be^{N-k},$$
where
$$C_k^N(\la,\mu)=\frac{1-q^{-2(\la+1)}}{1-q^{-2(\la+1)+2N}} 
\cdot\frac{q^{k(k-N)}(q^{-2\mu};q^2)_k}{(q^{2(k-N-\mu-1)};q^2)_k}.$$  
  \end{lemma}
  
  \begin{proof}
  By definition,
  $$S(\ga)^\ast=F(\la-1)\be,\ \ \ \ \ S(\al)^\ast
=\frac{F(\la-1)}{F(\mu-1)}\,\al,$$
  which gives
\begin{equation*}\begin{split}S(\ga^{N-k}\al^k)^\ast
&=(F(\la-1)\be)^{N-k}\left(\frac{F(\la-1)}{F(\mu-1)}\,\al\right)^k\\
&=F(\la-1)\dotsm F(\la-N+k)\be^{N-k}\frac{F(\la-1)\dotsm
 F(\la-k)}{F(\mu-1)\dotsm F(\mu-k)}\,\al^k\\
 &=\frac{F(\la-1)\dotsm F(\la-N)}{F(\mu-1+N-k)\dotsm
 F(\mu-N-2k)}\,\be^{N-k}\al^k\\
&=\frac{F(\la-1)\dotsm F(\la-N)}{F(\mu-1+N-k)\dotsm
F(\mu-N-2k)}\frac{q^{k(k-N)}(q^{-2\mu};q^2)_k}{(q^{2(k-N-\mu)};q^2)_k}\,
 \al^k\be^{N-k},\end{split}\end{equation*}
where we used Lemma \ref{cr} in the last step. Inserting
 $$F(\la-1)\dotsm F(\la-N)=\frac{1-q^{-2(\la+1)}}{1-q^{-2(\la+1)+2N}}$$
 and simplifying gives the desired expression.
   \end{proof}

Using this lemma we can combine (\ref{e1}) and (\ref{e2}) to the equality
\begin{equation}\label{e3}\begin{split}&\quad
\qb{N}{k}q^{k(k-N)}\sum_{j=0}^N C_j^N(\la,\ro)\al^j\be^{N-j}
\ot S(t_{kj}^N)^\ast\\
&=C_k^N(\la,\mu)\sum_{j=0}^N\qb{N}{j}q^{j(j-N)}\al^j\be^{N-j}
\ot t_{jk}^N,\end{split}
\end{equation}
with notation as in (\ref{ro}). We now observe that $C_k^N(\la,\mu)$ factors
as a part independent of $k$ and $\mu$ times a part independent of $\la$. 
We may therefore cancel the factors involving $\la$ on both sides of (\ref{e3})
and move the dynamical variables to the right side of the tensor product, 
which gives
\begin{equation*}\begin{split}&\quad
\sum_{j=0}^N \al^j\be^{N-j}
\ot\qb{N}{k}q^{k(k-N)+j(j-N)}
\frac{(q^{-2\la};q^2)_j}{(q^{2(j-N-\la-1)};q^2)_j} \, S(t_{kj}^N)^\ast\\
&=\sum_{j=0}^N\al^j\be^{N-j}
\ot\qb{N}{j}q^{k(k-N)+j(j-N)}
\frac{(q^{-2\mu};q^2)_k}{(q^{2(k-N-\mu-1)};q^2)_k}\, t_{jk}^N.
\end{split}\end{equation*}
We want to conclude that this identity holds termwise.
In view of  \eqref{mtr}, this follows if the family
$(\al^j\be^{N-j})_{j=0}^N$ is linearly independent 
over $\mu_r(\mg)$, which is indeed true by Lemma \ref{bl}.
This completes the proof of Proposition \ref{un}.

Note that, in the non-dynamical case, $\ast\circ S$ is the (antilinear) algebra
automorphism defined by $\be\leftrightarrow\ga$. For completeness, we also 
state
the symmetry of the matrix elements with respect to $\al\leftrightarrow\de$.

\begin{proposition} There is an algebra automorphism $\Phi$ of
$\frl$, defined on the generators by
$$\Phi:\,(\al,\be,\ga,\de,f(\la),g(\mu))\mapsto(\de,\be,\ga,\al,
 f(-2-\mu),g(-2-\la)).$$
Moreover, $\De\circ\Phi=\si\circ(\Phi\ot\Phi)\circ\De$ and
$$\qb{N}{k}q^{k(k-N)}\Phi(t_{kj}^N)=\qb{N}{j}q^{j(j-N)}t_{N-j,N-k}^N.$$ 
 \end{proposition}
 
 This can be proved similarly to Proposition \ref{un}, using instead of Lemma
 \ref{unl} the trivial identity $\Phi(\ga^{N-k}\al^k)=\ga^{N-k}\de^{k}$. Note
 that since $\Phi$ interchanges the formal limits
 $\la,\,\mu\rightarrow\pm\infty$, it reduces to the \emph{anti}-automorphism
 $\al\leftrightarrow\de$ in the  non-dynamical case; cf.~Remark \ref{flr}.

\begin{remark}\label{sr}
The symmetries $\Phi$ and $\Psi=\ast\circ S$ of the matrix elements permute the
four parameter domains of Theorem \ref{lmep}. In fact, the symmetries can be
obtained from the explicit expressions given there  (note  that 
$\Phi(\Xi)=\Psi(\Xi)=\Xi$). Conversely, after proving  Theorem \ref{lmep}
 in one of the four cases, we can use the symmetries to deduce the remaining
 three.
\end{remark}

Combining Proposition \ref{un} and \eqref{mes} one obtains the 
orthogonality relations
\begin{equation}\label{aor}\begin{split}\de_{kl}&=\sum_{j=0}^N(t_{jk}^N)^\ast 
\frac{\qb Nj}{\qb Nk}\frac{(q^{2(j-N-\la-1)};q^2)_j(q^{-2\mu};q^2)_k}{
(q^{-2\la};q^2)_j(q^{2(k-N-\mu-1)};q^2)_k}\,t_{jl}^N\\
&=\sum_{j=0}^N t_{kj}^N (t_{lj}^N)^\ast\frac{\qb Nl}{\qb Nj} 
\frac{(q^{2(l-N-\la-1)};q^2)_l(q^{-2\mu};q^2)_j}{
(q^{-2\la};q^2)_l(q^{2(j-N-\mu-1)};q^2)_j}\end{split}\end{equation}
for matrix elements. Our next goal is to find commutative versions of these
identities by evaluating them in a representation of $\fr$. 
In fact, they will yield the ortho\-gonality relations for $q$-Racah
 polynomials. Thus we must first
discuss dynamical representations of $\fr$.

\section{Discrete orthogonality of matrix elements}

\subsection{Dynamical representations of $\frl$}
\label{ssecdynreps}
We need the concept of dynamical representations from \cite{ev}. However, we
prefer to realize these representations on vector spaces over $\mg$. 
Recall from \S \ref{scorep} that an
  $\gh$-space $V$  is a
diagonalizable $\gh$-module $V=\bigoplus_{\al\in\gh^\ast} V_\al$, where 
$V$ and $V_\al$ are vector spaces over $\mg$.
 Let
 $(D_{\gh,V})_{\al\be}$ be the space of $\mathbb C$-linear operators $U$
on  $V$
such that $U(gv)=T_{-\be}(g)U(v)$ and $U(V_\ga)\subseteq V_{\ga+\be-\al}$ for 
all $g\in\mg$, $v\in V$ and $\ga\in\gh^\ast$. Then the space
  $D_{\gh,V}=\bigoplus_{\al,\be\in\gh^\ast}
(D_{\gh,V})_{\al\be}$ is an $\gh$-algebra with
the  moment maps $\mu_l,\mu_r\colon \mh\to 
(D_{\gh,V})_{00}$ given by
$$\mu_l(f)(v)=T_{-\al}(f)v,\ \ \ \ \ \mu_r(f)(v)=f v,\ \ \ v\in V_\al.$$ 
We define a \emph{dynamical representation} of an $\gh$-algebra $A$ on  $V$ to
be an $\gh$-algebra homomorphism $A\rightarrow D_{\gh,V}$.    
 An \emph{intertwiner} between two dynamical representations 
$\pi_i:\,A\rightarrow D_{\gh,V_i}$, $i=1,\,2$, is an $\gh$-space morphism
$\phi:\,V_1\rightarrow V_2$ with $\phi\circ\pi_1(a)=\pi_2(a)\circ\phi$ for all
$a\in A$. 
If $V_0$ is a complex subspace of $V$ with $\mg V_0=V$, then $V_0$ is a
dynamical representation in the sense of \cite{ev}, and, conversely, if 
$V_0$ is a dynamical representation in the  sense of \cite{ev}, then 
$\mg\otimes V_0$ is one in
our sense. 

\begin{proposition}\label{dynreps} 
Let, for $\gh=\mathbb C$ and $\om\in\mathbb C$ arbitrary, 
${\mathcal H}^\om$ be the $\gh$-space with basis 
$\{e_k\}_{k=0}^\infty$ and the weight decomposition 
${\mathcal H}^\om=\bigoplus_{k=0}^\infty {\mathcal
H}_{\om+2k}^\om$, ${\mathcal H}^\om_{\om+2k}=\mg e_k$. 
Then there is a dynamical representation $\pi^\om\colon \frl \to 
D_{\gh, {\mathcal H}^\om}$, defined on the
generators by
\begin{eqnarray*}
&&\pi^\om(\al)\, ge_k = A_k(T_{-1}g)e_k, \quad 
\pi^\om(\be)\, ge_k = B_k(T_1g)e_{k-1}, \\
&&\pi^\om(\ga)\, ge_k = -q^{-1}(T_{-1}g)e_{k+1}, \quad
\pi^\om(\de)\, ge_k = D_k(T_1g)e_k,\\
&&\pi^\om(\mu_l(f))\,ge_k=(T_{-\om-2k}f)ge_k, \quad
\pi^\om(\mu_r(f))\,ge_k=fge_k,
\end{eqnarray*}
where $g\in\mg$ is arbitrary and
\begin{align*}
A_k(\la) &=  q^{-k}\frac{1-q^{2(\la-\om-k+1)}}{1-q^{2(\la-\om-2k+1)}},
\\
B_k(\la) &=  \frac{(1-q^{2k})(1-q^{2(\om+k-1)})}{(1-q^{2(\la+1)})
(1-q^{2(\om+2k-\la-3)})}, \quad B_0(\la)=0, \\
D_k(\la) &= q^k \frac{1-q^{2(\la+1-k)}}{1-q^{2(\la+1)}}.
\end{align*}
\end{proposition}

\begin{proof} 
We first check that $\pi^\om$ preserves the bigrading and the moment maps, 
and consequently 
the commutation relations \eqref{f2}. Then it suffices 
to check the first four
defining relations of Definition \ref{deffrm}, two relations 
from  (\ref{cross})
and  the determinant relation
$c=1$. This is straight-forward; for instance, the relation 
$\al\be=qF(\mu-1)\be\al$ is equivalent to
$$A_{k-1}(\la)B_k(\la-1)=qF(\la-1)B_k(\la)A_k(\la+1).$$
\end{proof}

The dynamical representation $\pi^\om$ is irreducible, in an obvious sense, 
if and only if $\om\notin\mathbb Z_{\leq 0}$. To see this, suppose that 
$U\subseteq \mathcal H^\om$ is  
an $\gh$-subspace which is closed under $\pi^\om$. 
By definition, $U=\bigoplus_{k\in\La}\mg e_k$ for some 
$\La\subseteq\Zp$. Since $U$ is closed under $\pi^\om(\ga)$, we have
$e_k\in\La\Rightarrow e_{k+1}\in\La$. If $\om\notin\mathbb Z_{\leq 0}$,
so that $B_k\neq 0$ for $k\geq 1$, we have also $e_{k+1}\in\La\Rightarrow
e_{k}\in\La$  for $k\geq 0$, so we can deduce that $U=0$ or 
$U=\mathcal H^\om$. On the other hand, if
$\om\in\mathbb Z_{\leq 0}$, then the subspace $\bigoplus_{k\geq 1-\om}\mg e_k$
 is invariant. 

\begin{remark}\label{ndl}The functions $A_k$, $B_k$, and $D_k$ have finite
limits as
$(\om-\la,\la)\rightarrow\pm\infty$ (two cases).
Let us consider the case $(\om-\la,\la)\rightarrow\infty$. On the level of
the algebra, this corresponds to $\la\rightarrow -\infty$,
$\mu\rightarrow\infty$. In view of \eqref{usl}, we let $\pi^\infty(\al)$,
$\pi^\infty(\be)$, $\pi^\infty(\ga)$, and $\pi^\infty(\de)$ be the operators on
$\bigoplus_{k=0}^\infty \mathbb C e_k$ which are formally obtained as the 
limits of $\pi^\om(\be)$, $\pi^\om(-q\de)$, $\pi^\om(\al)$, and 
$\pi^\om(-q\ga)$,
respectively. Let us also put $f_k=(q^2;q^2)_k^{-1/2}e_k$. Then we recover the 
well-known $\ast$-representation of $\mathcal
F_q(\mathrm{SU}(2))$ on $\ell^2(\Zp)$, cf.~\cite{vs}, given by
\begin{align*}\pi^\infty(\al)f_k&=\sqrt{1-q^{2k}}f_{k-1},
\ \ \ \ \ \pi^\infty(\be)f_k=-q^{k+1}f_k,\\
\pi^\infty(\ga)f_k&=q^kf_k,\ \ \ \ \ \pi^\infty(\de)f_k
=\sqrt{1-q^{2(k+1)}}f_{k+1}.
\end{align*}
  \end{remark}

\begin{lemma}\label{pixi} 
The element 
$\Xi$  defined in \emph{Lemma \ref{lx}} 
acts in the dynamical representation $\pi^\om$ by
$$
\pi^\om(\Xi) = (q^{\om-1}+q^{1-\om})\, \id. 
$$
\end{lemma}

This corresponds nicely to $\Xi$ being central.
To prove the lemma, we  use the first expression of \eqref{xi}. This gives
\begin{equation*}\begin{split}
\pi^\om(\Xi)\, ge_k =& 
\left[q^{-\om-2k+1}+q^{\om+2k-1}\right. \\
&\left.-q^{-(2\la-\om-2k+2)}(1-q^{2(\la-\om-2k+2)})(1-q^{2\la})
(-q^{-1})B_k(\la-1)\right] g e_k, 
\end{split}\end{equation*}
which indeed simplifies to $(q^{\om-1}+q^{1-\om})ge_k$.

We will now show that, in a certain sense, $(\pi^\om)_{\om\in\R}$
are unitarizable representations of $\fr$. Note 
that the $\ast$-operator on $D_\gh$, cf.~\S \ref{ssectSU2}, is the formal 
adjoint with respect to the formal inner product
 $\langle f,g\rangle = 
\int_\R  f(\la)\overline{g(\la)}\, d\la$ on $\mh$.
Similarly, we want to find functions
$\Ga_k\in\mg$ such that 
\begin{equation}\label{fuy}\pi^\om(x^\ast)=\pi^\om(x)^\ast,\qquad
x\in\fr,\end{equation}
 where the $\ast$
on the right-hand side is the formal adjoint with respect to the formal 
pairing
\begin{equation*}\langle fe_k,ge_l\rangle= 
\de_{kl}\int_\R f(\la)\overline{g(\la)} 
\,\Ga_k(\la)\,d\la\end{equation*}
on $\mathcal H^\om$.

If $x\in\fr_{jk}$, so that
\begin{align*}\pi^\om(x)(ge_l)&=X_l(T_{-k}g)e_{l+\frac 12(k-j)},\\
\pi^\om(x^\ast)(ge_l)&=X_l^\ast(T_{k}g)e_{l+\frac 12(j-k)}\end{align*}
for some $X_l$, $X_l^\ast\in\mg$, then
\begin{align*}\langle \pi^\om(x^\ast)fe_l,ge_{l+\frac 12(j-k)}\rangle&=\int
X_l^\ast(\la)f(\la+k)\overline{g(\la)}\,\Ga_{l+\frac12(j-k)}(\la)\,d\la,\\
 \langle fe_l,\pi^\om(x) ge_{l+\frac 12(j-k)}\rangle&=\int
f(\la)\overline{ X_{l+\frac 12(j-k)}(\la)g(\la-k)}\Ga_{l}(\la)\,d\la\\
&=\int
f(\la+k)\overline{ X_{l+\frac 12(j-k)}(\la+k)g(\la)}\Ga_{l}(\la+k)\,d\la,
\end{align*}
so we require (recall that we write $\bar f(\la)=\overline{f(\bar\la)}$)
\begin{equation}\label{pu}X_l^\ast(\la)
\Ga_{l+\frac 12(j-k)}(\la)=
\bar X_{l+\frac 12(j-k)}(\la+k)\Ga_l(\la+k).\end{equation}
This gives a precise meaning to  \eqref{fuy}.

Thus we must find $\Ga_l$ so that (\ref{pu}) holds for all generators.
For $x=\mu_l(f)$, (\ref{pu}) holds provided that
$\om\in\R$.  For the right moment map (\ref{pu}) is satisfied.
For $x=\al$, $x^\ast=\de$ and for $x=\be$, $x^\ast=-q\ga$, (\ref{pu}) takes the
form
\begin{align}\label{pua}
\Ga_k(\la+1) A_k(\la+1)&= \Ga_k(\la) D_k(\la),\\
\label{pub}\Ga_{k-1}(\la-1)B_k(\la-1)&=\Ga_k(\la),\end{align}
where we used that
  $D_k=\bar D_k$ for $\om\in\R$. Taking
$k=0$ in \eqref{pua} shows that $\Ga_0$ has to be $1$-periodic. 
Iterating \eqref{pub} then gives
$$\Ga_k(\la) = \prod_{i=1}^k B_{k-i+1}(\la-i) \, \Ga_0(\la-k)
= \Ga_0(\la)
\frac{(q^2,q^{2\om};q^2)_k}{(q^{2(\la-k+1)},q^{2(\om-\la+k-1)};q^2)_k}.$$ 
Choosing $\Ga_0(\la)=1$  we can immediately
verify that \eqref{pua} holds. Thus we have proved the following proposition.

\begin{proposition}\label{pun} 
 For $\om\in \R$,
$\pi^\om$ is, in a sense made precise above,
 a unitary representation of $\fr$ with respect to the
formal pairing
$$\langle fe_k,ge_l\rangle= \de_{kl}\int_\R f(\la)\overline{g(\la)} 
\,\frac{(q^2,q^{2\om};q^2)_k}{(q^{2(\la-k+1)},q^{2(\om-\la+k-1)};q^2)_k}
\,d\la.$$
\end{proposition}

\subsection{Discrete orthogonality relations}
We will now obtain commutative versions of the orthogonality relations
\eqref{aor}, by evaluating them in a representation $\pi^\om$.

\begin{proposition}\label{tl}
One has
$$\pi^\om(t_{kj}^N)e_m=T_{kjm}\,e_{m+j-k},$$
where  $T_{kjm}=T_{kjm}^{\om N}\in\mg$ can be expressed in terms of 
$q$-Racah polynomials
\eqref{defqracah} as
\begin{equation*}\begin{split}
T_{kjm}^{\om N}(\la)&=(-1)^{j+k}q^{2k(\la+1)+m(N-k-j)+(k-j)(N-j+1)}
\qb{N}{j}\\
&\quad\times\frac{(q^{2(\la+1+k-j-m)};q^2)_{N-k-j}(q^{-2(m+j)},
q^{-2(\om-1+m+j)};q^2)_k}
{(q^{2(\la+1-j)};q^2)_{N-j}(q^{2(\la-\om+2+N-2j-2m)};q^2)_k}\\
&\quad\times
R_j(\mu(k);q^{-2(N+1)},q^{-2(\la+1)},q^{-2(m+j+1)},q^{2(\la-\om+1-m-j)};q^2)
\end{split}\end{equation*}
where we use the notation \eqref{nps} if $N<j+k$. 
\end{proposition}

\begin{proof}
We use the expressions for $t_{kj}^N$ given in Theorem \ref{lmep}.
First suppose that $k\leq j$, $k+j\leq N$, so that 
$$\pi^\om\left(t_{kj}^N\right)e_m=\pi^{\om}(\ga^{j-k}\de^{N-k-j}P_k)e_m.$$
By Proposition \ref{dynreps} and Lemma \ref{pixi}, the element $P_k$ acts on 
$e_m$ by multiplying with an
element of $\mg$ which is obtained from $P_k$ by replacing 
$\mu\mapsto\la$, $\la\mapsto \la-\om-2m$, $\Xi\mapsto q^{\om-1}+q^{1-\om}$.
Then $\pi^{\om}(\ga^{j-k}\de^{N-k-j})$ acts by replacing $\la\mapsto\la+N-2j$, 
multiplying with
$$(-q^{-1})^{j-k}\prod_{l=0}^{N-k-j-1}D_m(\la+l+k-j)
=(-q^{-1})^{j-k}q^{m(N-k-j)}
\frac{(q^{2(\la+1+k-j-m)};q^2)_{N-k-j}}{(q^{2(\la+1+k-j)};q^2)_{N-k-j}}$$
and shifting $e_m\mapsto e_{m+j-k}$. In conclusion,
$\pi^\om(t_{kj}^N)e_m=T\,e_{m+j-k}$, where
\begin{multline}\label{taw}
\begin{split}T(\la)&=(-q^{-1})^{j-k}q^{m(N-k-j)}
\frac{(q^{2(\la+1+k-j-m)};q^2)_{N-k-j}}{(q^{2(\la+1+k-j)};q^2)_{N-k-j}}\,
P_k\left(\begin{smallmatrix}\la\mapsto\la-\om+N-2j-2m\\
 \mu\mapsto\la+N-2j\hfill \\
\Xi\mapsto  q^{\om-1}+q^{1-\om}\hfill\end{smallmatrix}\right)\,\\
&=(-1)^jq^{k(2\la-\om+3+N+k-3j-3m)+(m-j)(N-j)-j}\\
&\quad\times\frac{(q^2;q^2)_{N-k}
(q^{2(\la+1+k-j-m)};q^2)_{N-k-j}}
{(q^2;q^2)_{j}(q^2,q^{2(\la+1-j)};q^2)_{N-j}(q^{2(\la-\om+2+N-2j-2m)}
;q^2)_{k}}\end{split}\\
\times p_k(\textstyle\frac
12\displaystyle(q^{\om-1}+q^{1-\om});q^{1+\om+2(m+j-k)},q^{1-\om-2m},q^{3-\om+2
(\la+N-2j-m)},q^{\om-1-2(\la+k-j-m)};q^2).\end{multline}
Writing the Askey--Wilson polynomial as a ${}_4\phi_3$ and inverting 
the order of summation, we
obtain the desired expression.

In the remaining three cases of Theorem \ref{lmep}, the lemma can be proved 
similarly. Alternatively, one can use Proposition \ref{mep}, which leads to a
more involved computation, but allows one to treat all four cases
simultaneously.
\end{proof}

We now suppose that $\om\in\mathbb R$, and write with abuse of notation
\begin{equation}\label{nf}\Ga_k^N(\la)=\qb{N}{k}\frac{(q^{2(k-N-\la-1)};q^2)_k}
{(q^{-2\la};q^2)_k},\end{equation}
$$\Ga_k^\om(\la)=\frac{(q^2,q^{2\om};q^2)_k}{(q^{2(\la-k+1)},
q^{2(\om-\la+k-1)};q^2)_k}$$
for the normalizing functions of the corepresentation $V_N$ and the
representation $\mathcal H^\om$; cf.~Propositions \ref{un} and \ref{pun}.
Since $t_{kj}^N\in\fr_{2k-N,2j-N}$, it follows from 
\eqref{pu} and Proposition \ref{tl} that
$$\pi^\om((t_{kj}^N)^\ast)g e_m=T_{kj,m+k-j}(\la+2j-N)
\frac{\Ga_m^\om(\la+2j-N)}{\Ga_{m+k-j}^\om(\la)}\,g(\la+2j-N)\,e_{m+k-j}.$$
The first identity of \eqref{aor} then gives
\begin{equation*}\begin{split}\de_{kl}\,e_m&=\sum_{j=0}^N \pi^\om\left(
(t_{jk}^N)^\ast\frac{\Ga_j^N(\la)}{\Ga_k^N(\mu)}\,t_{jl}^N\right)e_m\\
&=\sum_{j=0}^{\min(N,m+l)} \pi^\om((t_{jk}^N)^\ast)
\frac{\Ga_j^N(\la-\om-2m-2l+2j)}{\Ga_k^N(\la)}
\,T_{jlm}(\la)e_{m+l-j}\\
&=\sum_{j=0}^{\min(N,m+l)}
T_{jk,m+l-k}(\la+2k-N)\frac{\Ga_{m+l-j}^\om(\la+2k-N)}{\Ga_{m+l-k}^\om(\la)}\\
&\quad\quad\times\frac{\Ga_j^N(\la-\om-2m-2l+2j+2k-N)}
{\Ga_k^N(\la+2k-N)}\,
T_{jlm}(\la+2k-N)e_{m+l-k}.\end{split}\end{equation*}
Replacing $\la$ by $\la+N-2k$ and $m$ by $M-l$ we obtain
$$\de_{kl}=\sum_{j=0}^{\min(M,N)}
\frac{\Ga_{M-j}^\om(\la)\Ga_j^N(\la-\om-2M+2j)}
{\Ga_{M-k}^\om(\la+N-2k)\Ga_k^N(\la)}\,T_{jk,M-k}(\la)T_{jl,M-l}(\la)$$ 
for $0\leq k,\,l\leq\min(M,N)$.
Using Proposition \ref{tl} we can write this explicitly  as
\begin{multline}\label{rk}\sum_{j=0}^{\min(M,N)}\frac{1-q^{2(\la-\om+1-2M)+4j}}
{1-q^{2(\la-\om+1-2M)}}\frac{(q^{2(\la-\om-2M+1)},q^{-2N},q^{-2M},
q^{-2(\om-1+M)};q^2)_j}
{(q^2,q^{2(\la-\om+2+N-2M)},q^{2(\la-\om+2-M)},q^{2(\la+1-M)};q^2)_j}\\
\begin{split}&\quad\times  q^{2j(\la+N+1)}
R_k(\mu(j);q^{-2(N+1)},q^{-2(\la+1)},q^{-2(M+1)},q^{2(\la-\om+1-M)};q^2)
\\
&\quad\times R_l(\mu(j);q^{-2(N+1)},q^{-2(\la+1)},q^{-2(M+1)},
q^{2(\la-\om+1-M)};q^2) \\
&=\de_{kl}\,\frac{(q^{2(\om-\la-1+2M-N)},q^{-2(\la+N)};q^2)_N}
{(q^{2(\om-\la-1+M-N)},q^{-2(\la+N-M)};q^2)_N}\frac{1-q^{-2(\la+1+N)}}
{1-q^{-2(\la+1+N)+4k}}\\
&\quad\times\frac
{(q^2,q^{-2\la},q^{-2(\la+N-M)},q^{2(\om-\la-1+M-N)};q^2)_k}
{(q^{-2(\la+1+N)},q^{-2N},q^{-2M},q^{-2(\om-1+M)};q^2)_k}\,
q^{2k(\la-\om+1-2M)},
\end{split}\end{multline}
which is the  orthogonality of $q$-Racah polynomials.
In this case, $\{R_k\}_{k=0}^{\min(M,N)}$ is a complete system of
orthogonal polynomials on $\{\mu(j)\}_{j=0}^{\min(M,N)}$. 
Similarly, the second equation in \eqref{aor} gives the orthogonality of the
dual system
$$R_j(\mu(k);q^{-2(M+1)},q^{2(\la-\om+1-M)},q^{-2(N+1)},q^{-2(\la+1)};q^2).$$
In the limits $\la,\om-\la\rightarrow \pm\infty$, 
cf.~Remark \ref{ndl}, these relations reduce to 
the ortho\-gonality of quantum
$q$-Krawtchouk polynomials \cite{ko}.

\begin{remark}
We could have considered more general representations $\mathcal H^{\om,\ep}$,
$\om,\,\ep\in\R$, defined by the same formulas as in Proposition
\ref{dynreps}, but with the basis $\{e_k\}_{k\in\mathbb Z+\ep}$.
For $\ep=0$, $\mathcal H^\om$ occurs as the invariant subspace
spanned by $\{e_k\}_{k\in\Zp}$.
Suppose for simplicity that $\ep,\,\om+\ep\notin\mathbb Z$. Then
Proposition \ref{pun} is valid for $\mathcal H^{\om,\ep}$, with the
convention \eqref{cps}.
Working with these representations would lead to more general $q$-Racah
polynomials, with $M$ in \eqref{rk} replaced by a continuous parameter.
\end{remark}

\section{Clebsch--Gordan coefficients for corepresentations}

\subsection{Tensor products of $\gh$-coalgebroids}
\label{sscoalg}

Our next goal is to compute the Clebsch--Gordan coefficients of our
corepresentations. We first need to discuss some additional 
algebraic concepts, in particular $\gh$-coalgebroids and their tensor products.
 These concepts will also play a crucial role in Appendix 1. 

We define an \emph{$\gh$-prealgebra} $A$ to be a complex vector space,
 equipped with a
bigrading $A=\bigoplus_{\al,\be\in\gh^\ast}A_{\al\be}$ 
and two left actions $\mu_l,\,\mu_r:\,\mg\rightarrow\End_\mathbb C(A)$
which preserve the bigrading, such that the images of $\mu_l$
and $\mu_r$ commute.
We also introduce two right actions
of $\mg$ on $A$ by (\ref{dr}). 
A homomorphism of $\gh$-prealgebras is
a linear map which preserves the four $\mg$-actions, or equivalently 
  the two left actions and the bigrading.

If $A$ and $B$ are $\gh$-prealgebras we define their matrix tensor product
$A\widetilde\otimes B$ as in \S \ref{sshb}. We also define another kind of 
tensor
product $A\widehat\otimes B$ which is equal to the algebraic tensor product
 modulo
the relations
\begin{equation}\label{fh}a\mu_l^A(f)\otimes b=a\otimes\mu_l^B(f)b,
\ \ \ \ \  a\mu_r^A(f)\otimes b=a\otimes\mu_r^B(f)b.\end{equation}
The bigrading $A_{\al\be}\wh\otimes
A_{\ga\de}\subseteq(A\widehat\otimes B)_{\al+\ga,\,\be+\de}$
and the moment maps 
\begin{align*}\mu_l^{A\widehat\otimes B}(f)(a\otimes b)&=\mu_l^A(f)a\otimes
b,\\
\mu_r^{A\widehat\otimes B}(f)(a\otimes b)&=\mu_r^A(f) a\otimes b,\end{align*}
 make $A\widehat\otimes B$ an
$\gh$-prealgebra.
It follows from these definitions that
\begin{align*}(a\otimes b)\mu_l^{A\widehat\otimes B}(f)&
=a\otimes b\mu_l^B(f),\\
(a\otimes b)\mu_r^{A\widehat\otimes B}(f)&=a\otimes
b\mu_r^B(f).\end{align*}
We will need the following two lemmas, the first of which is trivial.

\begin{lemma}\label{pht}
If $\phi:\,A\rightarrow C$, $\psi:\,B\rightarrow D$ are homomorphisms
of $\gh$-prealgebras, then $\phi\otimes\psi$ factors to $\gh$-prealgebra
homomorphisms $A\wh\ot B\rightarrow C\wh\ot D$ and 
$A\wt\ot B\rightarrow C\wt\ot D$. 
 \end{lemma}

\begin{lemma}\label{l23} If $A$, $B$, $C$, $D$ are $\gh$-prealgebras, then
$$\si_{23}(a\ot b\ot c\ot d)=a\ot c\ot b\ot d$$
factors to an $\gh$-prealgebra homomorphism
\begin{equation}\label{s23}\si_{23}:\,(A\wt\ot B)\wh\ot 
(C\wt\ot D)\rightarrow (A\wh\ot C)\wt\ot(B\wh\ot
D).\end{equation}
\end{lemma}

\begin{proof}
There are two things to be checked: first, that $\si_{23}$ maps into the 
subspace
of $A\ot B\ot C\ot D$ which maps onto $(A\wh\ot C)\wt\ot(B\wh\ot D)$, second,
that $\si_{23}$ factors through the defining relations of $(A\wt\ot B)\wh\ot
(C\wt\ot D)$. 

For the first part, note that the left hand side of \eqref{s23} splits into
a sum of quotients of spaces of the form 
$A_{\al\be}\ot B_{\be\ga}\ot C_{\de\ep}\ot D_{\ep\ze}$. This component is 
mapped to
$A_{\al\be}\ot C_{\de\ep}\ot B_{\be\ga}\ot D_{\ep\ze}$ by $\si_{23}$, and
is then projected to $(A\wh\ot C)_{\al+\de,\be+\ep}\ot(B\wh\ot
D)_{\be+\ep,\ga+\ze}$. Here we may indeed replace $\ot$ by $\wt\ot$. 

For the second part, we write down the relations valid on both sides of
\eqref{s23} explicitly. Those on the left-hand side may be written as
$$\mu_ra= \mu_lb,\ \ \ \mu_rc=\mu_l d,
\ \ \ a\mu_l=\mu_lc,\ \ \ b\mu_r=\mu_r d,$$
where, for instance, the first identity is an abbreviation for
$$\mu_r^A(f)a\ot b\ot c \ot d=a\ot \mu_l^B(f)b\ot c\ot d,
\ \ \ a\in A,\,b\in B,\, c\in C,\,d\in D,\,f\in\mg,$$
and those on the right-hand side as
$$\mu_ra=\mu_lb,\ \ \ a\mu_l=\mu_lc,\ \ \ a\mu_r=\mu_rc,\ \ \ b\mu_l=\mu_ld,
\ \ \ b\mu_r=\mu_rd.$$
We must show that the second group of relations implies the first group. 
This is
clear except for the relation $\mu_rc=\mu_ld$. However, by \eqref{mtr2}
 we have also that $a\mu_r=b\mu_l$, and
 thus indeed $\mu_rc=a\mu_r=b\mu_l=\mu_ld$.
\end{proof}

We define an \emph{$\gh$-coalgebroid} to be an $\gh$-prealgebra equipped with a
coproduct
and a counit satisfying the same axioms as in the case of $\gh$-bialgebroids,
except that they are required to be $\gh$-prealgebra homomorphisms rather than
$\gh$-algebra homomorphisms. An $\gh$-coalgebroid homomorphism
$\phi:\,A\rightarrow B$ is an $\gh$-prealgebra homomorphism with
$(\phi\otimes\phi)\circ\De^A=\De^B\circ\phi$, $\ep^B\circ\phi=\ep^A$. An
$\gh$-android is a person who studies $\gh$-algebroids.

\begin{proposition}\label{cat}
If $A$ and $B$ are $\gh$-coalgebroids, then 
$\De^{A\wh\ot B}=\si_{23}\circ(\De^A\ot \De^B)$ and
$\ep^{A\widehat\otimes B}(a\otimes b)=\ep^A(a)\ep^B(b)$
(a composition of difference operators) define 
an $\gh$-coalgebroid structure on
$A\widehat\otimes B$.
\end{proposition}

More explicitly, if $\De^A(a)=\sum_i a_i'\otimes a_i''$ 
and $\De^B(b)=\sum_j b_j'\otimes b_j''$, then we define
\begin{equation}\label{de}\De^{A\widehat\otimes B}(a\otimes b)
=\sum_{ij} a_i'\otimes b_j'\otimes a_i''\otimes b_j''.\end{equation}
As for the proof, note that it follows from the previous two lemmas that
$\De^{A\wh\ot B}$ is a well-defined $\gh$-prealgebra homomorphism. The 
remaining details are  straight-forward.

If, in particular, $A$ is an $\gh$-bialgebroid, there is an $\gh$-algebra
structure on $A\widetilde\otimes A$ and an $\gh$-coalgebroid structure on
$A\widehat\otimes A$. 
One may check that the multiplication factors to an $\gh$-coalgebroid
 homomorphism
$A\widehat\otimes A\rightarrow A$, similarly as the coproduct is an 
$\gh$-algebra
homomorphism $A\rightarrow A\widetilde\otimes A$. In the case of bialgebras 
($\gh=0$) these structures combine to a bialgebra structure on
 $A\widetilde\otimes A=A\widehat\otimes A=A\otimes A$. By contrast, there is
apparently  no natural tensor product  on the class of 
 $\gh$-bialgebroids (that is, $\gh$-bialgebroids do not form a 
 monoidal category in a natural way).

\subsection{Tensor products of corepresentations}
\label{sstc}

We will now discuss tensor products of corepresentations in
general. In \S \ref{scorep} we defined corepresentations of $\gh$-bialgebroids;
 this definition extends mutatis mutandis to corepresentations of
$\gh$-coalgebroids. 

When $V$ and $W$ are $\gh$-spaces we denote by $V\wh\ot W$ their 
tensor product over $\mathbb C$ modulo the relations
$$fv\ot w=v\ot T_\al f w,\ \ \ \ \ v\in V_{\al}.$$
The grading $V_\al\wh\ot W_\be\subseteq (V\wh\ot W)_{\al+\be}$ and the action 
of scalars $f(v\ot w)=fv\ot w$ make $V\wh\ot W$ into an $\gh$-space. Note
 that if we view $\gh$-prealgebras as $\gh$-spaces by ``forgetting''
their left (or right) moment map and grading, then the
tensor product $A\wh\ot B$ of $\gh$-prealgebras introduced above reduces
to the one defined here.

\begin{proposition}\label{crt}
If $V$ and $W$ are corepresentation spaces of an $\gh$-coalgebroid $A$, 
then there is a corepresentation 
$$\pi_{V\wh\ot W}:\,V\wh\ot W\rightarrow A\wt\ot(V\wh\ot W)$$  
of $A$ on on $V\wh\ot W$ defined by
\begin{equation}\label{ct}\pi_{V\wh\ot W}=(m\ot\id)\circ\si_{23}\circ
(\pi_V\ot\pi_W).\end{equation}
If $A$ is an $\gh$-bialgebroid, then the multiplication 
$m:\,A\wh\ot A\rightarrow A$ is an intertwiner for
the regular corepresentation. 
\end{proposition}

To prove the first statement, note that it
follows from Lemmas \ref{pht} and \ref{l23}
that $\pi_{V\wh\ot W}$ is a well-defined $\gh$-space morphism. 
 The remaining details are exactly as for coalgebras. The second part of the
 proposition follows from the
 fact that   $\De$ is an algebra homomorphism, which can be expressed as
$$\De\circ m=(m\ot m)\circ \si_{23}\circ(\De\ot\De).$$
This means precisely that $m$ is an intertwiner.

If $\{v_k\}_k$, $\{w_k\}_k$ are bases (over $\mg$) of  $V$ and
$W$, and $t_{kj}^V$, $t_{kj}^W$ are matrix elements of $\pi_V$ and $\pi_W$ with
respect to these bases, then \eqref{ct} may be written as
$$\pi_{V\wh\ot W}(v_j\ot w_k)=\sum_{lm}t_{jl}^V \,t_{km}^W\ot v_l\ot w_m,$$
or more compactly as
\begin{equation}\label{tt}t_{jk,lm}^{V\wh\ot W}=t_{jl}^V\, t_{km}^W.
\end{equation}

Now suppose that we are given three corepresentations 
$U$, $V$, $W$ and an
intertwiner $C:\,U\wh\ot V\rightarrow W$. If we pick bases of the
corepresentation spaces, the matrix elements $C_{jk,l}\in\mg$ defined
by
\begin{equation}\label{cint}C(u_j\ot v_k)=\sum_l C_{jk,l}\,w_l\end{equation}
are Clebsch--Gordan coefficients.
Combining (\ref{int2}) and (\ref{tt}) gives the intertwining property in terms
of these coefficients:
\begin{equation}\label{cgc}\sum_{kl}C_{kl,p}(\mu)t_{mk}^U \,t_{nl}^V=\sum_j
C_{mn,j}(\la)t_{jp}^W\ \ \ \text{for all}\ m,\,n,\,p.\end{equation}

Let us now apply the coproduct $\De$ to \eqref{cgc}. By \eqref{mer}, we
obtain 
\begin{equation*}\begin{split}\sum_{klxy}C_{kl,p}(\mu)t_{mx}^U \,t_{ny}^V\ot
t_{xk}^U\,t_{yl}^V&=\sum_{jz}
C_{mn,j}(\la)t_{jz}^W\ot t_{zp}^W\\
&=\sum_{xyz}C_{xy,z}(\ro)t_{mx}^U\,t_{ny}^V\ot t_{zp}^W,
\end{split}\end{equation*}
where we again applied \eqref{cgc}  and used the notation of
\eqref{ro}. 
Assuming that the family
$(t_{mx}^U t_{ny}^V)_{xy}$ is linearly independent over $\mu_r(\mg)$,
it follows that
$$\sum_{kl}C_{kl,p}(\mu)
t_{xk}^U\,t_{yl}^V=\sum_{z}C_{xy,z}(\la) t_{zp}^W, $$
which is \eqref{cgc} with $(m,n)$ replaced by $(x,y)$. Thus,
to prove that the operator $C$ defined by \eqref{cint} is
intertwining, it suffices to check \eqref{cgc} for all $p$ and with $m$ and $n$
fixed such that $(t_{mk}^U t_{nl}^V)_{kl}$ is  independent over $\mu_r(\mg)$.

\subsection{Clebsch--Gordan coefficients for $\fr$}

We are now ready to compute the Clebsch--Gordan coefficients of our
corepresentations. The proof will be similar to the one in \cite{kk} for the 
non-dynamical case. 

In analogy with the classical case, we will prove that
\begin{equation}\label{cgd}V_M\wh\ot V_N\simeq 
\bigoplus_{s=0}^{\min(M,N)}V_{M+N-2s}
\end{equation}
(direct sums of corepresentations may be defined in a 
straight-forward way).
Therefore we assume that $C:\,V_M\wh\ot V_N\rightarrow V_{M+N-2s}$
is an intertwiner, where $0\leq
s\leq\min(M,N)$. We write the Clebsch--Gordan coefficients with respect to the
standard bases as
$$C(\ga^{M-j}\al^{j}\ot\ga^{N-k}\al^{k})=\sum_{l=0}^{M+N-2s}
C_{jk,l}^{MN,M+N-2s}
(\la)\ga^{M+N-2s-l}\al^{l}.$$
Since $C$ preserves the grading, one has
$C_{jk,l}^{MN,M+N-2s}=0$ unless $j+k=l+s$, so (\ref{cgc}) may be written as
\begin{equation}\label{cg2}\sum_{\substack{k+l=p+s\\0\leq k\leq
        M\\0\leq l\leq N}}C_{kl,p}^{MN,M+N-2s}(\mu)t_{mk}^M \,
t_{nl}^N=C_{mn,m+n-s}^{MN,M+N-2s}(\la)t_{m+n-s,p}^{M+N-2s}.\end{equation}

We now observe that, by Lemma \ref{bl}, the elements 
$$t_{0k}^M\,t_{Nl}^N=\qb{M}{k}\qb{N}{l}q^{k(k-M)+l(l-N)}\ga^k\de^{M-k}
\al^l\be^{N-l},\ \ \ 0\leq k\leq M,\ 0\leq l\leq N,$$
 are linearly independent over $\mu_r(\mg)$. 
 Thus, by the observation concluding \S \ref{sstc}, it
 suffices to find $C$ so that \eqref{cg2} holds for $m=0$ and $n=N$, that is,
  so that
\begin{multline}\label{cf}C_{0N,N-s}^{MN,M+N-2s}(\la)
t_{N-s,\,p}^{M+N-2s}
=\sum_{l=\max(0,\,p+s-M)}^{\min(N,\,p+s)}\qb{M}{p+s-l}\qb{N}{l}
 \\
 \times q^{(p+s-l)(p+s-l-M)+l(l-N)}C_{p+s-l,l,p}^{MN,M+N-2s}(\mu)
\ga^{p+s-l}\de^{M-p-s+l}\al^l\be^{N-l}.\end{multline}
Note that, by Proposition
\ref{crt}, we can for $s=0$ choose $C$ as the
multiplication $V_M\wh\ot
V_N\rightarrow V_{M+N}$. 
One may then check that (\ref{cf}) reduces to the expression for
matrix elements given in Proposition \ref{mep}. 
We will in fact compute the Clebsch--Gordan coefficients by
deducing an expression of the form (\ref{cf}) from Proposition
\ref{mep}. For this we need the following lemma.

\begin{lemma}
In the algebra $\frl$,
\begin{equation}\label{hd}1=\sum_{m=0}^s\qb{s}{m}(-1)^{s-m}q^{(2m+1)(m-s)}
\frac{(q^{-2\mu};q^2)_{s-m}}{(q^{-2(\mu+m+1)};q^2)_{s-m}}
\,\ga^{s-m}\de^m\al^m\be^{s-m}.\end{equation}
\end{lemma}

\begin{proof}
We sketch two proofs. For the first one, we make the Ansatz
$$c^s=\sum_{m=0}^s B_{sm}(\la,\mu)\ga^{s-m}\de^m\al^m\be^{s-m},$$
where $c$ is the dynamical determinant of Lemma \ref{lc}. Writing
$$c^{s+1}=\left(\de\al-\frac{q^{-1}}{F(\mu-1)}\,\ga\be\right)c^s=
\de c^s\al
-\frac{q^{-1}}{F(\mu-1)}\,\ga c^s\be$$
and using the commutation relations of Lemma \ref{cr}, one derives the
recursion relations
\begin{equation*}\begin{split}B_{s+1,m}(\la,\mu)&=q^{2(m-1-s)}\left(
\frac{1-q^{-2\mu+2s-2m}}{1-q^{-2\mu-2}}\right)^2B_{s,m-1}(\la+1,\mu+1)\\
&\quad
-\frac{q^{-1}}{F(\mu-1)}\,B_{sm}(\la+1,\mu-1),\end{split}\end{equation*}
for $m=0,\dots,s+1$, where $B_{s,-1}=B_{s,s+1}=0$. This leads to an identity 
equivalent to \eqref{dp}.

Alternatively, one may plug the expression for
$\de^m\al^m$ given in Lemma \ref{ccr} into the right-hand side of \eqref{hd}, 
commute $\ga^{s-m}$ across this
expression (using Lemma \ref{lx}) and then again apply Lemma
\ref{ccr} to the factor $\ga^{s-m}\be^{s-m}$. This results in an identity 
involving only commuting variables, which turns out to be the
terminating ${}_6W_5$ summation formula \cite{gr}:
$$1=\frac{(q^{1+\la-\mu}\xi,q^{1+\la-\mu}\xi^{-1};q^2)_s}
{(q^{-2(\mu+1)},q^{2(\la+2)};q^2)_s}\,{}_6W_5(q^{2(1+\mu-s)};
q^{\la+\mu+3}\xi,q^{\la+\mu+3}\xi^{-1},q^{-2s};q^2,q^{-2(\la+1)}),$$
with $\xi$  a formal quantity satisfying $\xi+\xi^{-1}=\Xi$. 
\end{proof}

To compute the Clebsch--Gordan coefficients we insert the right-hand side of
\eqref{hd} into the expressions for matrix elements given in Proposition
\ref{mep}.
Using Lemma \ref{cr} twice and pulling the dynamical variables to the left 
gives {\allowdisplaybreaks
\begin{multline*}t_{N-s,p}^{M+N-2s}=\sum_{l=\max(0,\,p+s-M)}^{\min(N-s,\,p)}
\qb{M-s}{p-l}\qb{N-s}{l}q^{p(p+N-M)+l(3l+M+s-2N-3p)}\\
\begin{split}&\quad\times\frac{(q^{2(p+2s-M-N-\mu-1)};q^2)_{p-l}}
{(q^{2(p+s-M-l-\mu-1)};q^2)_{p-l}}\,\ga^{p-l}\de^{M-p-s+l}
\left(\sum_{m=0}^s\qb{s}{m}(-1)^{s-m}q^{(2m+1)(m-s)}\right.\\
&\quad\left.\times
\frac{(q^{-2\mu};q^2)_{s-m}}{(q^{-2(\mu+m+1)};q^2)_{s-m}}
\,\ga^{s-m}\de^m\al^m\be^{s-m}
\right)\al^l\be^{N-s-l}\end{split}\\ 
\begin{split}&
=\sum_{l=\max(0,\,p+s-M)}^{\min(N-s,\,p)}\sum_{m=0}^s
\qb{M-s}{p-l}\qb{N-s}{l}\qb{s}{m}(-1)^{s-m}\\
&\quad\times q^{p(p+N-M)+l(3l+M+s-2N-3p)+(M-p-s+2l+2m+1)(m-s)}\\
&\quad\times \frac{(q^{2(p+2s-M-N-\mu-1)};q^2)_{p-l}
(q^{2(2p+s-M-l-\mu)},q^{2(p-l-\mu)};q^2)_{s-m}}
{(q^{2(p+s-M-l-\mu-1)};q^2)_{p-l}(q^{2(2p+s-M-2l-m-\mu-1)},
q^{2(2p+s-M-2l-\mu)};q^2)_{s-m}}\\
&\quad\times\ga^{p+s-l-m}\de^{M-p-s+l+m}\al^{l+m}\be^{N-l-m}.
\end{split}\end{multline*} }
Replacing $l$ by $l-m$, we see that
 \eqref{cf} holds with 
\begin{multline}\label{cgw}C_{jk,j+k-s}^{MN,M+N-2s}(\la)=
\sum_{m=\max(0,s-j,k+s-N)}^{\min(s,k,M-j)}
\frac{\qb{M-s}{j+m-s}\qb{N-s}{k-m}\qb{s}{m}}{\qb{M}{j}\qb{N}{k}}\,(-1)^{s+m}\\
\begin{split}&\times
q^{(N-k-s)(j-s)-s+m(1+2j+2N+3m-2k-4s)}\\
&\times\frac{(q^{2(j+k+s-M-N-\la-1)};q^2)_{j+m-s}
(q^{2(2j+k+m-s-M-\la)},q^{2(j+m-s-\la)};q^2)_{s-m}}
{(q^{2(j+m-M-\la-1)};q^2)_{j+m-s}
(q^{2(2j+2m-s-M-\la)},q^{2(2j+m-s-M-\la-1)};q^2)_{s-m}}
\end{split}\end{multline}
(in particular, the Clebsch--Gordan coefficient on the left-hand side of
\eqref{cf} reduces to $1$). For all values of the parameters, this is an
 ${}_8W_7$-sum; 
for instance, when $s\leq \min(j,N-k)$ we get
\begin{multline*}C_{jk,j+k-s}^{MN,M+N-2s}(\la)=(-1)^sq^{(j-s)(N-k-s)-s}
\frac{(q^{2};q^2)_{M-s}(q^{2};q^2)_{N-s}(q^{2};q^2)_{j}(q^{2};q^2)_{N-k}}
{(q^{2};q^2)_{M}(q^{2};q^2)_{N}(q^{2};q^2)_{j-s}(q^{2};q^2)_{N-k-s}}\\
\begin{split} &\times\frac{(q^{2(j+k+s-M-N-\la-1)};q^2)_{j-s}
(q^{2(2j+k-s-M-\la)},
 q^{2(j-s-\la)};q^2)_{s}}
{(q^{2(j-M-\la-1)};q^2)_{j}(q^{2(2j-s-M-\la)};q^2)_{s}}\,
{}_8W_7\left(q^{2(2j-s-M-\la-1)};\right.\\
&\quad\left. q^{-2s},q^{2(j-M-\la-1)},q^{-2k},q^{2(j-M)},
q^{2(2j+k-M-N-\la-1)};q^2,q^{2(2+M+N-s)}\right).
 \end{split}\end{multline*} 
Using \eqref{tr1} we can write this as 
\begin{equation}\label{cgp}\begin{split}C_{jk,j+k-s}^{MN,M+N-2s}(\la)&=
q^{(j-s)(N-k)+sj}\frac{(q^{-2(j+k)};q^2)_s(q^{2(j+k-M-N-\la-1)};q^2)_j}
{(q^{-2N};q^2)_s(q^{2(j-M-\la-1)};q^2)_j}\\
&\quad\times\,
{}_4\phi_3\left[\begin{matrix}q^{-2s},q^{2(s-M-N-1)},q^{-2j},q^{2(j-M-\la-1)}\\
q^{-2M},q^{-2(j+k)},q^{2(j+k-M-N-\la-1)}
\end{matrix};q^2,q^2\right].
\end{split}\end{equation}
It is easy to verify that \eqref{cgp} holds also without the assumption 
$s\leq \min(j,N-k)$.

For each $s\leq\min(M,N)$, we have constructed an intertwiner 
$C_s:\,V_M\wh\ot
V_N\rightarrow V_{M+N-2s}$, which is clearly surjective. 
Since both sides of \eqref{cgd} have the same dimension, 
 $\oplus_s C_s$ is a bijective intertwiner from the left-hand to the
right-hand side. Thus we have proved the following theorem. 

\begin{theorem}\label{cgpr}
The decomposition \eqref{cgd} holds as an equivalence of corepresentations.
The equation
\begin{equation}\label{cgf}C_s(\ga^{M-j}\al^j\ot\ga^{N-k}\al^k)
=C_{jk,j+k-s}^{MN,M+N-2s}(\la)\ga^{M+N-s-j-k}
\al^{j+k-s},\end{equation}
where the coefficients are given in \eqref{cgp}, defines an intertwiner
$C_s:\,V_M\wh\ot V_N\rightarrow V_{M+N-2s}$. In particular, 
the Clebsch--Gordan formula  \eqref{cg2} holds.
\end{theorem}

\subsection{Orthogonality of the Clebsch--Gordan coefficients}

We will now discuss the orthogonality of Clebsch--Gordan coefficients,
 which is a consequence of Schur's Lemma
(Corollary \ref{sl}) and the unitarizability of the corepresentations 
(Proposition
\ref{un}). We will see that it yields orthogonality relations for $q$-Racah 
polynomials.  We start with two general facts about unitarizable
corepresentations.

\begin{lemma} \label{ucl}
Let $V_1$ and $V_2$ be two unitarizable corepresentation spaces of a
$\ast$-$\gh$-Hopf algebroid $A$, 
and let $\Ga_k^i$, $i=1,\,2$,  be
 normalizing functions with respect to some
 bases $\{v_k^i\}_k$ of $V_i$. 
 
 If $\phi:\,V_1\rightarrow V_2$ is an intertwiner, and
  $\phi_{kj}\in\mg$  the matrix elements of
 $\phi$ with respect to the bases $\{v_{k}^i\}_k$, as in \eqref{ime}, then
  $\phi^\ast_{kj}=(\Ga_j^1/\Ga_k^2)\bar\phi_{jk}$
defines an intertwiner $\phi^\ast:\,V_2\rightarrow V_1$.

Moreover, the tensor product corepresentation $V_1\wh\ot V_2$ 
is unitarizable, with
$$\Ga_{jk}(\la)=\Ga_j^1(\la)\Ga_k^2(\la-\om(j))$$
 a normalizing function with respect to the basis $\{v_j^1\ot v_k^2\}_{jk}$,
where $\om(j)\in\gh^{\ast}$ is defined by $v_j^1\in (V_1)_{\om(j)}$.
\end{lemma}

\begin{proof}
For the first statement, we
 apply $\ast\circ S$ to \eqref{int2}. Using  the unitarizability we obtain
$$\sum_{j}\bar\phi_{jl}(\la)\frac{\Ga_j^1(\la)}{\Ga_k^1(\mu)}\,t_{jk}^1=\sum_j
\bar\phi_{kj}(\mu)\frac{\Ga_l^2(\la)}{\Ga_j^2(\mu)}\,t_{lj}^2
\ \ \ \text{for all}\ k,\,l,$$
which means precisely that $\phi^\ast$ is intertwining.

For the second statement, we apply 
 $\ast\circ S$ to \eqref{tt}. This gives indeed
$$S(t_{jk,lm}^{V_1\wh\ot V_2})^\ast=\frac{\Ga_l^1(\la)}{\Ga_j^1(\mu)}\,
t_{lj}^1\,
\frac{\Ga_m^2(\la)}{\Ga_k^2(\mu)}\,t_{mk}^2=\frac{\Ga_l^1(\la)
\Ga_m^2(\la-\om(l))}
{\Ga_j^1(\mu)\Ga_k^2(\mu-\om(j))}\,t_{lm,jk}^{V_1\wh\ot V_2},$$
where we used that $t_{lj}^1\in A_{\om(l),\om(j)}$. 
\end{proof}

In view of the second
 part of Lemma \ref{ucl}, we may apply the first part to the intertwiner
$C_s:\,V_M\wh\ot V_N\rightarrow V_{M+N-2s}$ defined in \eqref{cgf}.
 This yields an 
intertwiner $C_s^\ast:\,V_{M+N-2s}\rightarrow V_M\wh\ot V_N$, defined by
$$C^\ast_s(\ga^{M+N-2s-l}\al^l)=\sum_{\substack{j+k=l+s\\0\leq j\leq M\\0\leq
k\leq N}}
\frac{\Ga_j^MT_{M-2j}\Ga_k^N}{\Ga_{l}^{M+N-2s}}\,C_{jk,j+k-s}^{MN,M+N-2s}
\ga^{M-j}\al^j\ot\ga^{N-k}\al^k,$$
where the Clebsch--Gordan coefficients are given by \eqref{cgp} and the
normalizing functions  by \eqref{nf}.

Let us now consider the map $C_s C_t^\ast:\,V_{M+N-2t}\rightarrow V_{M+N-2s}$,
where $0\leq s,\,t\leq\min(M,N)$. By Corollary \ref{sl}, 
\begin{equation}\label{cca}C_sC_t^\ast=\de_{st}Z_s\id \end{equation}
for some $Z_s\in\mathbb C$. Applying this identity to $\ga^{M+N-t-L}\al^{L-t}$
gives the orthogonality relations for Clebsch--Gordan coefficients:
\begin{equation}\label{cgo}\de_{st}Z_s=\sum_{\substack{j+k=L\\0\leq j\leq M\\
0\leq k\leq N}}\frac{\Ga_j^M(\la)\Ga_k^N(\la+M-2j)}{\Ga_{L-t}^{M+N-2t}(\la)}\,
C_{jk,L-s}^{MN,M+N-2s}(\la)C_{jk,L-t}^{MN,M+N-2t}(\la),\end{equation}
where $0\leq s,\,t\leq \min(L,M,N,M+N-L)$.

To compute $Z_s$, we specialize \eqref{cgo} to the case $L=s=t$. 
Then the Clebsch--Gordan coefficients on the right-hand side simplify,
since   
the sum in \eqref{cgw} reduces to the term with $m=k$, giving 
\begin{equation}\label{csp}C_{jk,0}^{MN,M+N-2s}(\la)=(-1)^jq^{k(N-s)-j}
\frac{(q^2;q^2)_s(q^2;q^2)_{M-j}
(q^2;q^2)_{N-k}(q^{-2\la};q^2)_j}{(q^2;q^2)_M(q^2;q^2)_N
(q^{2(j-M-\la-1)};q^2)_j}
\end{equation}
for $j+k=s$.
Plugging in this
expression and simplifying, the right-hand side of \eqref{cgo}
 reduces to a ${}_6W_5$ sum:
\begin{align}\nonumber Z_s&=(-1)^sq^{-s(s+1)}
\frac{(q^2,q^{2(s-M-N-\la-1)};q^2)_s}{(q^{-2N},q^{-2(\la+M)};q^2)_s}\\
\nonumber &\quad\times
\,{}_6W_5(q^{-2(\la+1+M)};q^{-2s},q^{2(1+N-s)},q^{-2\la};q^2,q^{2(2s-M-N-1)})\\
\label{zs}&=(-1)^sq^{-s(s+1)}
\frac{(q^2,q^{2(s-M-N-1)};q^2)_s}{(q^{-2M},q^{-2N};q^2)_s},
\end{align}
where the second step is obtained either by applying \cite[(II.21)]{gr}
or by observing that, 
since we know a priori that $Z_s$ is independent of $\la$, we can 
 put $\la=0$, so that the ${}_6W_5$ reduces to
$1$.

 Using \eqref{cgp} one may check that \eqref{cgo} is the
orthogonality of the $q$-Racah polynomials
\begin{equation}\label{cgr}R_s(\mu(j);q^{-2(M+1)},q^{-2(N+1)},q^{-2(L+1)},
q^{2(L-M-\la-1)};q^2),\end{equation}
cf.~\eqref{defqracah}.
For this special case,  $\{R_s\}_{s=0}^{\min(L,M,N,M+N-L)}$
is, for generic $\la$, a complete system of polynomials orthogonal on
$\{\mu(j)\}_{j=\max(0,L-N)}^{\min(L,M)}$.

To obtain the dual orthogonality relations, we  observe that 
\begin{equation}\label{du}\id\big|_{V_M\wh\ot V_N}=\sum_{s=0}^{\min(M,N)}
\frac{1}{Z_s}\,C_s^{\ast}C_s.\end{equation}
In fact, it follows from \eqref{cca} that the restriction of both sides to
the image of $C_t^\ast$ are equal for $0\leq t\leq\min(M,N)$. 
By Theorem \ref{cgpr}, these images span $V_m\wh\ot V_N$.
Applying \eqref{du} to a tensor product $\ga^{M-j}\al^j\ot\ga^{N-L+j}\al^{L-j}$
 gives
\begin{equation}\label{dcgo}\begin{split}\de_{jk}&
=\sum_{s=0}^{\min(L,M,N,M+N-L)}
\frac{1}{Z_s}\frac{\Ga_j^M(\la)\Ga_{L-j}^N(\la+M-2j)}
{\Ga_{L-s}^{M+N-2s}(\la)}\\
&\quad\times C_{j,L-j,L-s}^{MN,M+N-2s}(\la)C_{k,L-k,L-s}^{MN,M+N-2s}(\la),
\end{split}\end{equation}
for $\max(0,L-n)\leq j,k\leq\min(L,M)$.
This is the orthogonality of the system
$$R_j(\mu(s);q^{-2(L+1)},
q^{2(L-M-\la-1)},q^{-2(M+1)},q^{-2(N+1)};q^2) $$
dual to \eqref{cgr}.
In the non-dynamical limit, \eqref{cgo} and \eqref{dcgo} are orthogonality 
relations for  $q$-Hahn and dual $q$-Hahn
polynomials, respectively; cf.~\cite{kk}. 

Another consequence of \eqref{du} is
the Clebsch--Gordan formula dual to \eqref{cg2}.

\begin{proposition}
The identity
\begin{equation}\label{lin}\begin{split}t_{kx}^M\,t_{ly}^N
&=\sum_s\frac{1}{Z_s}\frac{\Ga_x^M(\mu)\Ga_y^N(\mu+M-2x)}
{\Ga^{M+N-2s}_{x+y-s}(\mu)}\\
&\quad\quad\times C_{kl,k+l-s}^{MN,M+N-2s}(\la)
\,C_{xy,x+y-s}^{MN,M+N-2s}(\mu)\,t_{k+l-s,x+y-s}^{M+N-2s},
\end{split}\end{equation}
holds, where $0\leq k,x\leq M$, $0\leq l,y\leq N$,
 the sum runs over $s$ with 
$$0\leq s\leq\min(M,N,k+l,M+N-k-l,x+y,M+N-x-y),$$
and $Z_s$ is given by \eqref{zs}.
\end{proposition}

This is proved by applying
 \eqref{du} to $\pi^{V_M\wh\ot V_N}(\ga^{M-k}\al^k\ot
\ga^{N-l}\al^l)$, using the intertwining property of $C_s$, 
and identifying the coefficient of
$\ga^{M-x}\al^x\ot\ga^{N-y}\al^y$.
If we apply the counit $\ep$ to \eqref{lin}, we recover
\eqref{dcgo}.

One can obtain commutative versions of the Clebsch--Gordan formulas 
\eqref{cg2} and \eqref{lin} by evaluating them in a dynamical representation.
Namely, applying the representation $\pi^\om$ of Proposition \ref{dynreps} to 
\eqref{cg2} and acting with both sides on $e_j\in\mathcal H^\om$ gives the
identity
\begin{multline}\label{conv}
\sum_{k+l=p+s}C_{kl,p}^{MN,M+N-2s}(\la)\,
T_{mk,j+l-n}^{\om M}(\la)\,T_{nlj}^{\om N}(\la+M-2k)\\
=C_{mn,m+n-s}^{MN,M+N-2s}(\la-\om-2(j+p+s-m-n))
\,T_{m+n-s,pj}^{\om,M+N-2s}(\la),\end{multline}
where $T_{kjm}^{\om N}$ is the function  from Proposition \ref{tl},
while \eqref{lin}  similarly gives
\begin{multline}\label{p2}T_{kx,j+y-l}^{\om M}(\la)\,T_{lyj}^{\om N}(\la+M-2x)
=\sum_s\frac 1{Z_s}\frac{\Ga_x^M(\la)\Ga_y^N(\la+M-2x)}
{\Ga^{M+N-2s}_{x+y-s}(\la)}\\
\times C_{kl,k+l-s}^{MN,M+N-2s}(\la-\om-2(j+x+y-k-l))
\,C_{xy,x+y-s}^{MN,M+N-2s}(\la)\,T_{k+l-s,x+y-s,j}^{\om,M+N-2s}(\la).
\end{multline}
Using  Proposition \ref{tl} and \eqref{cgp} we may express these identities 
in terms of  $q$-Racah polynomials, as indicated in \eqref{rk} and
\eqref{cgr}. It turns out that in both cases we obtain instances of the
Biedenharn--Elliott identity, which will be 
discussed in more detail in \S \ref{prss}.
In the non-dynamical case, \eqref{conv} and \eqref{p2} are
essentially different identities, cf.~\cite{kv,kv2} for related results. 

\section{Clebsch--Gordan coefficients for representations}

\subsection{Tensor products of dynamical representations}
\label{tpdr}

In this section we will obtain the Clebsch--Gordan decomposition of 
the dynamical representations
introduced in \S \ref{ssecdynreps}. 
Since our definition of dynamical
representations differs slightly from the one in \cite{ev}, we must 
accordingly
modify the definition of tensor product representations.

When $V$ and $W$ are $\gh$-spaces, we denote by
 $V\bar\otimes W$  their tensor product over $\mathbb C$
  modulo the relations
$$fv\otimes w=v\otimes T_{-\be}fw,\quad w\in W_\be.$$
The grading  $V_\al\bar\ot W_\be\subseteq (V\bar\ot W)_{\al+\be}$ and the 
action of scalars $f(v\ot w)=v\ot fw$ make $V\bar\ot W$ into an $\gh$-space. 
This is
closely related to
 the tensor product $V\wh\ot W$ introduced in \S
\ref{sstc}; in fact, the flip map $v\otimes w\mapsto w\otimes v$ defines an
$\gh$-space isomorphism $V\bar \ot W\mapsto W\wh\ot V$.

Let $\pi_V:\,A\rightarrow D_{\gh,V}$ and $\pi_W:\,A\rightarrow D_{\gh,W}$ be
 two dynamical representations of an $\gh$-algebra $A$. One may check that  
the identity operator factors to an $\gh$-algebra
homo\-morphism $\Theta:\,D_{\gh,V}\wt\ot D_{\gh,W}\rightarrow 
D_{\gh,V\bar\ot W}$.
Then 
$$\pi_{V\bar \ot W}=\Theta\circ(\pi_V\otimes \pi_W)\circ\Delta $$
defines a dynamical representation of $A$ on $V\bar\ot W$. 

We will obtain the decomposition
\begin{equation}\label{trd}\mathcal H^{\om_1}\bar\ot
\mathcal H^{\om_2}\simeq\bigoplus_{s=0}^\infty
\mathcal H^{\om_1+\om_2+2s}\end{equation}
(direct sums of dynamical representations may be defined in a straight-forward
way).
In view of Lemma \ref{pixi}, we can achieve this by 
diagonalizing the 
 action of $\Xi$ in the tensor product representation 
 $\pi=\pi_{\mathcal H^{\om_1}\bar\ot\mathcal H^{\om_2}}$.
Note the resemblance of \eqref{trd} to the case of highest weight
(discrete series)
representations of  $\mathfrak{su}(1,1)$.  
 In this analogy, $\Xi$ corresponds to the Casimir operator and $\om_i$ to the
 highest weights.   
  
From (\ref{xi}) we see that we first have to consider
$$
\De(\ga\be) = \ga\al\ot \al\be + \ga\be\ot\al\de +
\de\al\ot\ga\be +\de\be\ot\ga\de, 
$$
where we used 
Definition \ref{deffrm}. Thus, $\pi(\ga\be)(e_{k_1}\otimes e_{k_2})$ can be
written as a sum of four terms, the first of which is
\begin{equation*}\begin{split}\pi^{\om_1}(\ga\al)e_{k_1}\ot
\pi^{\om_2}(\al\be)e_{k_2}&
=-q^{-1}(T_{-1}A_{k_1}^{\om_1})e_{k_1+1}\otimes
A_{k_2-1}^{\om_2}(T_{-1}B_{k_2}^{\om_2})e_{k_2-1}\\
&=-q^{-1}(T_{-(\om_2+2(k_2-1))-1}A_{k_1}^{\om_1})A_{k_2-1}^{\om_2}
(T_{-1}B_{k_2}^{\om_2})(e_{k_1+1}\ot e_{k_2-1}),
\end{split}\end{equation*}
where we have written $\om_i$ as superscripts to the 
functions
$A_k$ and $B_k$ from Proposition \ref{dynreps} to indicate
their dependence on $\om$. Computing the other three terms similarly
gives 
\begin{equation}\label{tpdragabe}
\begin{split}\pi(\ga\be)\, e_{k_1}\ot e_{k_2} 
&= a_{k_1,k_2}  (e_{k_1+1}\ot e_{k_2-1})
+ b_{k_1,k_2} (e_{k_1}\ot e_{k_2}) \\ 
&\quad + c_{k_1,k_2}(e_{k_1-1}\ot e_{k_2+1}), 
\end{split}
\end{equation}
with 
\begin{align*}\label{ttr}
-q a_{k_1,k_2}(\la)&= 
B^{\om_2}_{k_2}(\la-1)A^{\om_2}_{k_2-1}(\la)
A_{k_1}^{\om_1}(\la-\om_2-2k_2+1), \\
-q b_{k_1,k_2}(\la)&= 
D^{\om_2}_{k_2}(\la-1)A^{\om_2}_{k_2}(\la)
B_{k_1}^{\om_1}(\la-\om_2-2k_2-1)\\ &\quad +
B^{\om_2}_{k_2}(\la-1)
A_{k_1}^{\om_1}(\la-\om_2-2k_2+1)D_{k_1}^{\om_1}(\la-\om_2-2k_2), \\
-q c_{k_1,k_2}(\la)&=
D^{\om_2}_{k_2}(\la-1)
B_{k_1}^{\om_1}(\la-\om_2-2k_2-1)D_{k_1-1}^{\om_1}(\la-\om_2-2k_2-2).
\end{align*}
From (\ref{tpdragabe}) 
or from $\ga\be\in\frl_{00}$ it follows
that $\pi(\ga\be)$ preserves the weight spaces 
$({\mathcal H}^{\om_1}\bar\ot{\mathcal H}^{\om_2})_{\om_1+\om_2+2p}=
\bigoplus_{k_1+k_2=p} {\mathcal H}^{\om_1}_{\om_1+2k_1}\bar\ot 
{\mathcal H}^{\om_2}_{\om_1+2k_2}$, so that in order to diagonalize
$\pi(\ga\be)$ it suffices to diagonalize $\pi(\ga\be)$ in every
weight space. Fix $p\in\Zp$; then it follows from 
(\ref{tpdragabe}) that $\sum_{k=0}^p v_k (e_{p-k}\ot e_k)$, $v_k\in\mg$, is
an eigenvector of $\pi(\ga\be)$ in the weight space
 $({\mathcal H}^{\om_1}\bar\ot{\mathcal H}^{\om_2})_{\om_1+\om_2+2p}$
 with
eigenvalue $x$ if and only if the $v_k$'s satisfy 
\begin{equation}\label{ttrforpktobesolved}
x\, v_k = a_{p-k-1,k+1} v_{k+1} + 
b_{p-k,k} v_k  + 
c_{p-k+1,k-1} v_{k-1},\end{equation}
where $v_{-1}=v_{p+1}=0$. 

The three-term recurrence relation (\ref{ttrforpktobesolved}) 
can be solved in terms
of $q$-Racah polynomials. We recall that the polynomials
$R_n(\mu(x))=R_n(\mu(x);a,b,c,d;q)$, cf.~\eqref{defqracah}, satisfy the
recurrence
\begin{equation}\label{ttrqRacah}
\begin{split}
\bigl( \mu(x)-\mu(0)\bigr)\, R_n(\mu(x)) &= 
A_n \bigl( R_{n+1}(\mu(x))- R_n(\mu(x))\bigr)\\ &\quad  + 
C_n \bigl( R_{n-1}(\mu(x))- R_n(\mu(x))\bigr), 
\end{split}
\end{equation}
where
\begin{equation*}
\begin{split} 
&A_n=\frac{(1-abq^{n+1})(1-aq^{n+1})(1-bdq^{n+1})(1-cq^{n+1})}
{(1-abq^{2n+1})(1-abq^{2n+2})}, \\
&C_n=\frac{q(1-q^n)(1-bq^n)(c-abq^n)(d-aq^n)}
{(1-abq^{2n})(1-abq^{2n+1})},
\end{split}
\end{equation*}
which holds for $n,x\in\{0,1,\ldots,N\}$ if
$aq$, $bdq$ or $cq$ equals $q^{-N}$, $N\in\Zp$. 

Upon replacing in (\ref{ttrforpktobesolved}) 
$$
v_k = (-1)^k q^{2k(\om_1+\om_2)+k(3p-1)} 
\frac{(q^{-2\la},q^{-2(\om_1+p-1)}, q^{-2p};q^2)_k}
{(q^2,q^{2\om_2},q^{2(\om_1+\om_2+p-\la-1)};q^2)_k}\, 
R_k,
$$
we find after a straightforward calculation 
that the resulting three-term recurrence relation 
can be written as 
\begin{equation}\label{f}\begin{split}&\quad
q^{1-2(2p+\om_1+\om_2-1)} (1-q^{2\la})(1-q^{2(2p+\om_1+\om_2-\la-2)})
\, xR_k\\
& = a_k(R_{k+1}-R_k) + c_k(R_{k-1}-R_k),\end{split}\end{equation}
where
\begin{equation*}
\begin{split}a_k&=\frac{(1-q^{2(k+\om_2-\la-1)})(1-q^{2(k-\la)})
(1-q^{2(k-\om_1-p+1)})(1-q^{2(k-p)})}
{(1-q^{2(2k+\om_2-\la-1)})(1-q^{2(2k+\om_2-\la)})}, \\
c_k&=\frac{q^2(1-q^{2k})(1-q^{2(k+\om_2-1)})
(q^{-2(p+1)}-q^{2(k+\om_2-\la-2)})
(q^{-2(\om_1+\om_2+p-1)}-q^{2(k-\la-1)})}
{(1-q^{2(2k+\om_2-\la-2)})(1-q^{2(2k+\om_2-\la-1)})}.
\end{split}
\end{equation*}
The right-hand side of \eqref{f} is the right-hand side of
(\ref{ttrqRacah}) in base $q^2$ with $a\mapsto q^{-2(\la+1)}$, 
$b\mapsto q^{2(\om_2-1)}$, $c\mapsto q^{-2(p+1)}$ and
$d\mapsto q^{-2(\om_1+\om_2+p-1)}$, so that $N=p$. 
Comparing the left-hand sides, one finds that the eigenvalue $x$  is of the 
form 
\begin{equation*}
 x = -q^{2(\om_1+\om_2+2p)-3}\frac{(1-q^{-2z})(1-q^{2(z+1-\om_1-\om_2-2p)})}
{(1-q^{2\la})(1-q^{2(2p+\om_1+\om_2-\la-2)})}
\end{equation*}
for $z\in\{0,1,\ldots,p\}$. The corresponding eigenvalue
of $\pi(\Xi)$ can then be computed from 
(\ref{xi})  using the fact that we restrict 
to the weight space 
$({\mathcal H}^{\om_1}\bar\ot{\mathcal H}^{\om_2})_{\om_1+\om_2+2p}$:
\begin{multline*}
q^{1-\om_1-\om_2-2p}+q^{\om_1+\om_2+2p-1} 
- q^{-(2\la-\om_1-\om_2-2p+2)}(1-q^{2(\la-\om_1-\om_2-2p+2)})
(1-q^{2\la})x\\
=q^{\om_1+\om_2+2(p-z)-1}+ q^{1-\om_1-\om_2-2(p-z)}\hfill\end{multline*}
for $z\in\{0,1,\ldots,p\}$. Thus we have proved the following
proposition. 

\begin{proposition}\label{propeigvxitp} In the tensor
product representation 
$\pi=\pi_{\mathcal H^{\om_1}\bar\ot\mathcal H^{\om_2}}$
the element $\Xi$
has eigenvectors $v(y;p)\in
({\mathcal H}^{\om_1}\bar\ot{\mathcal H}^{\om_2})_{\om_1+\om_2+2p}$, 
$y\in\{0,1,\ldots,p\}$, with the eigenvalue 
 $q^{\om_1+\om_2+2y-1}+ q^{1-\om_1-\om_2-2y}$. 
The eigenvector $v(y;p)$
is given in terms 
of $q$-Racah polynomials \emph{(\ref{defqracah})} by 
\begin{equation*}
\begin{split}
v(y;p)&= \sum_{k=0}^p v_k(e_{p-k}\ot e_k), \\
v_k(\la)&= (-1)^k q^{2k(\om_1+\om_2)+k(3p-1)} 
\frac{(q^{-2\la},q^{-2(\om_1+p-1)}, q^{-2p};q^2)_k}
{(q^2,q^{2\om_2},q^{2(\om_1+\om_2+p-\la-1)};q^2)_k}\\ 
&\quad\times  R_k(\mu(p-y); q^{-2(\la+1)}, q^{2(\om_2-1)},
q^{-2(p+1)},q^{-2(\om_1+\om_2+p-1)};q^2).
\end{split}
\end{equation*}
\end{proposition}

Note that for $\om_1+\om_2\notin \mathbb Z_{\leq 0}$ the eigenvalues per
weight space are different and independent of $\la$, 
so that the eigenvectors are
linearly independent over $\mh$ 
and form a basis for the tensor product
representation space. 

From now on we
assume the genericity condition   $\om_1+\om_2\notin \mathbb Z_{\leq 0}$. 
We need to calculate 
$\pi(\ga) v(y;p)$. Since $\pi(\ga)$
commutes with $\pi(\Xi)$ and raises the degree by $2$, we have 
$\pi(\ga) v(y;p) = C v(y;p+1)$ for some $C\in\mg$. On the other hand, using  
 $\De(\ga)=\ga\ot\al +\de\ot\ga$  we find that $\pi(\ga)v(y;p)$ equals
\begin{equation*}-q^{-1}\sum_{k=0}^p v_k(\la-1)
\left(A^{\om_2}_k(\la) (e_{p-k+1}\ot  e_k) + 
D^{\om_1}_{p-k}(\la-\om_2-2k-2) (e_{p-k}\ot e_{k+1})\right).
\end{equation*}
Comparing the coefficient of $e_{p+1}\ot e_0$ yields $C(\la)=-q^{-1}$ or
\begin{equation}\label{gt}\pi(\ga) v(y;p) = -q^{-1} v(y;p+1).\end{equation}
 A similar computation gives 
\begin{equation}\label{at} \pi(\al) v(y;p) = q^{-p} 
\frac{1-q^{2(\la-\om_1-\om_2-p+1)}}
{1-q^{2(\la-\om_1-\om_2-2p+1)}} \,  v(y;p).\end{equation}

We can now construct an intertwiner $C:\,\mathcal
H^{\om_1+\om_2+2s}\rightarrow\mathcal H^{\om_1}\bar\ot\mathcal H^{\om_2}$. 
Since $C$ must preserve the grading and the eigenspaces of $\Xi$,
one has
$$C e_k=\phi_{k}\,v(s;s+k) $$ 
for some $\phi_{k}\in\mg$.
Then the intertwining property
\begin{equation}\label{itr}C\circ \pi_{\mathcal H^{\om_1+\om_2+2s}}(x)=
\pi_{\mathcal H^{\om_1}\bar\ot\mathcal H^{\om_2}}(x)\circ C,\end{equation}
is automatically satisfied for $x=\mu_l(f)$, $x=\mu_r(f)$ and $x=\Xi$, and thus
for $x\in\frl_{00}$. Using \eqref{at} and \eqref{gt} we write out \eqref{itr}
explicitly for $x=\al$ and $x=\ga$. This leads to the equations
\begin{equation*}\begin{split}
\phi_k(\la)&=q^{-s}
\frac{1-q^{2(\la-\om_1-\om_2-s-k+1)}}{1-q^{2(\la-\om_1-\om_2-2s-k+1)}}\,
\phi_k(\la-1)\\
\phi_{k-1}(\la)&=\phi_k(\la+1),
\end{split}\end{equation*}
which are solved by
$\phi_k(\la)=q^{s(k-\la)}(q^{2(\la+2-\om_1-\om_2-2s-k)};q^2)_s$ 
(the general solution is obtained by
multiplying each $\phi_k$ with a fixed $1$-periodic function).
 With this choice  of $\phi_k$, we know that \eqref{itr} holds for $x=\al$,
 $x=\ga$,  $x=\al\de$ and $x=\ga\be$. Since it
is clear from \eqref{gt} and \eqref{at} that $\pi(\ga)$ and $\pi(\al)$ are
injective, we can conclude that \eqref{itr} holds also for $x=\de$ and 
$x=\be$, 
and thus for any $x\in\frl$.
Since we have already observed that the eigenvectors $v(y;p)$ form an
$\mg$-basis of $\mathcal H^{\om_1}\bar\ot\mathcal H^{\om_2}$, 
the following theorem is now clear.

\begin{theorem}\label{decompdtp} Assuming that $\om_1+\om_2\notin 
\mathbb Z_{\leq
0}$, the decomposition \eqref{trd} holds as an equivalence of dynamical
representations. Moreover,
$$C e_k=q^{s(k-\la)}(q^{2(\la+2-\om_1-\om_2-2s-k)};q^2)_s\,v(s;s+k),
$$
with the notation of \emph{Proposition \ref{propeigvxitp}},
 defines an intertwiner   $C:\,\mathcal
H^{\om_1+\om_2+2s}
\rightarrow\mathcal H^{\om_1}\bar\ot\mathcal H^{\om_2}$.
\end{theorem}

We can now interpret Proposition \ref{propeigvxitp} as stating that the
$q$-Racah polynomials are Clebsch--Gordan coefficients for the representations
$\pi^\om$. In analogy with corepresentations, we will write 
$$Ce_k=\sum_{l+m=s+k}C_{k,lm}^{\om_1+\om_2+2s,\om_1\om_2}(e_l\ot e_m),$$
where, writing $L=l+m=k+s$,
\begin{equation}\label{rmk}\begin{split}C_{k,lm}^{\om_1+\om_2+2s,\om_1\om_2}
(\la)&=(-1)^m q^{s(k-\la)+ 2m(\om_1+\om_2)+m(3L-1)}\\ 
&\quad\times\frac{(q^{-2\la},q^{-2(\om_1+L-1)}, 
q^{-2L};q^2)_m(q^{2(\la+2-\om_1-\om_2-s-L)};q^2)_s}
{(q^2,q^{2\om_2},q^{2(\om_1+\om_2+L-\la-1)};q^2)_m}\\
&\quad\times R_m(\mu(k);q^{-2(\la+1)},q^{2(\om_2-1)},q^{-2(L+1)},
q^{-2(\om_1+\om_2+L-1)};q^2).\end{split}\end{equation}
 For later use we note the alternative expression
\begin{equation}\begin{split}\label{cga}&
C_{k,lm}^{\om_1+\om_2+2s,\om_1\om_2}(\la)=(-1)^s 
q^{s(\la+1-L)+ 2l(1-\om_1-\om_2-L)-lm} \qb{L}{m}\\
&\times\frac{(q^{-2\la},q^{-2(\om_1+L-1)};q^2)_m(q^{2\om_2};q^2)_L
(q^{2(\om_1+\om_2-\la-1+m+L)};q^2)_l}
{(q^{-2\la},q^{-2(\om_1+L-1)};q^2)_k(q^{2\om_2};q^2)_m(q^{2\om_2};q^2)_s}\\
&\times R_l(\mu(k);q^{2(\la-\om_1-\om_2+1-2L)},q^{2(\om_1-1)},q^{-2(L+1)},
q^{-2(\om_1+\om_2+L-1)};q^2),\end{split}\end{equation}
which follows by combining \eqref{rmk} and \eqref{rs2}. Finally we remark that 
if we instead use \eqref{rs1} we obtain
\begin{equation*}\begin{split} C_{k,lm}^{\om_1+\om_2+2s,\om_1\om_2}(\la)&=
q^{s(k-\la)-lm}
(q^{2(\la+2-\om_1-\om_2-s-L)};q^2)_s\qb{L}{m}\\
&\quad\times R_s(\mu(m);q^{2(\om_2-1)},q^{2(\om_1-1)},q^{-2(L+1)},
q^{2(L+\om_2-\la-1)};q^2).
\end{split}\end{equation*}
Comparing with \eqref{cgr}, we see that the Clebsch--Gordan coefficients for
$\mathcal H^{\om_1}\bar\ot\mathcal H^{\om_2}$ can be obtained from those of
$V_M\wh\ot V_N$ by the formal replacements $\om_1=-N$, $\om_2=-M$. 
This is analogous to the connection between highest weight representations of 
(Lie or quantum) $\mathrm{SU}(2)$ and $\mathrm{SU}(1,1)$.

\subsection{The pentagonal identity}
\label{prss}

The classical and quantum $6j$-symbols satisfy the pentagonal 
or Biedenharn--Elliott
identity. For the $6j$-symbols of $\mathcal U_q(\mathfrak{su}(2))$ \cite{kr}, 
this identity can be written in terms of $q$-Racah polynomials in many ways,
such as 
\begin{multline}\label{pi}\frac{(aq,q^{-m_1-m_2},bdq^{1-m_3};q)_{k_1}
(abq^{2k_1+2};q)_{k_2}(bcq^2;q)_{k_1+k_2}}{(bq,bcdq^2,q^{-m_1-m_2-m_3};q)_{k_1}
(cq;q)_{k_2}(aq;q)_{k_1+k_2}}\,(bq^{1+k_1})^{-k_2}\\
\begin{split}&\quad\times R_{k_1}(\mu(m_1);a,b,q^{-(m_1+m_2+1)},dq^{-m_3};q)\\
&\quad\times R_{k_2}(\mu(m_1+m_2-k_1);abq^{2k_1+1},c,q^{-(m_1+m_2+m_3-k_1+1)},
bdq^{1+k_1};q)\\
&=\sum_{l=0}^{\min(k_1+k_2,m_2+m_3)}\frac{(bcq,q^{-k_1-k_2},q^{-m_2-m_3},
abcq^{k_1+k_2+2},bq,bcdq^{2+m_1};q)_l}{(q,bcq^{k_1+k_2+2},q^{-m_1-m_2-m_3},
q^{-k_1-k_2}a^{-1},cq,bcdq^2;q)_l}\frac{1-bcq^{2l+1}}{1-bcq}\\
&\quad\times(abq^{1+m_1})^{-l}\, R_{k_1}(\mu(l);b,a,q^{-(k_1+k_2+1)},
bcq^{k_1+k_2+1};q)\\
&\quad\times R_l(\mu(m_2);b,c,q^{-(m_2+m_3+1)},bdq^{1+m_1};q)\\
&\quad\times
R_{k_1+k_2-l}(\mu(m_1);a,bcq^{2l+1},q^{-(m_1+m_2+m_3-l+1)},dq^{-l};q).
\end{split}\end{multline}
From the viewpoint of special functions, \eqref{pi} is 
a master identity which contains many classical results as limit
cases, including
 various convolution, linearization and addition formulas for
orthogonal polynomials. We will be concerned with the case $a=q^{-N-1}$, 
$N\in\Zp$. If $N<k_1+k_2$, there is then a singularity which must be removed 
by multiplying with $(aq;q)_{k_1+k_2}$ and interpreting
$$\frac{(aq;q)_{k_1+k_2}}{(q^{-k_1-k_2}a^{-1};q)_{l}}=(-1)^l
q^{-\binom{l}{2}+l(k_1+k_2-N-1)}(q^{-N};q)_{k_1+k_2-l},$$
so that the  summation can be restricted to 
$$\max(0,k_1+k_2-N)\leq l\leq \min(k_1+k_2,m_1+m_2).$$

From the interpretation in terms of $\mathcal U_q(\mathfrak{su}(2))$, 
one only obtains \eqref{pi}
 for discrete values of the $9$ free parameters (not counting $q$). 
It can be extended to continuous values of $a$, $b$, $c$, $d$ by working 
instead
with $\mathcal U_q(\mathfrak{su}(1,1))$.
The equations \eqref{conv} and \eqref{p2} are instances of \eqref{pi} 
with   $2$
continuous parameters. In this section we will obtain 
\eqref{pi} with $3$ continuous parameters
using our interpretation of $q$-Racah polynomials as Clebsch--Gordan 
coefficients of dynamical representations. We point out that
an extension of \eqref{pi} to the case of $9$ continuous parameters 
was obtained in \cite{kv2}, again using a quantum algebraic interpretation.
This involves not necessarily terminating ${}_8W_7$-series.  

To obtain the pentagonal identity, we corepresent the intertwining property
\eqref{itr} by applying it to $x=t_{kj}^N$ and acting on
 $e_m\in \mathcal H^{\om_1+\om_2+2s}$:
\begin{equation}\label{add}
C \pi^{\om_1+\om_2+2s}(t_{kj}^N)e_m=\pi^{\om}(t_{kj}^N)C e_m.
\end{equation}
The left-hand side of this identity is
$$C
T_{kjm}^{\om_1+\om_2+2s,N}e_{m+j-k}=\sum_{x+y=s+m+j-k}
C_{m+j-k,xy}^{\om_1+\om_2+2s,\om_1\om_2}T_{kjm}^{\om_1+\om_2+2s,N}
(e_x\ot
e_y),$$
where $T_{kjm}^{\om N}$ is the function  from Proposition \ref{tl}.
The right-hand side of \eqref{add} is
\begin{multline*}\sum_{l=0}^N\pi^{\om_1}(t_{kl}^N)\otimes\pi^{\om_2}(t_{lj}^N) 
\sum_{x+y=s+m} e_x\otimes C_{m,xy}^{\om_1+\om_2+2s,\om_1\om_2}e_y\\
\begin{split}
&=\sum_{l=0}^N\sum_{x+y=s+m}T_{klx}^{\om_1N}e_{x+l-k}\ot T_{ljy}^{\om_2N}
(T_{N-2j}C_{m,xy}^{\om_1+\om_2+2s,\om_1\om_2})e_{y+j-l}\\
&=\sum_{l=0}^N\sum_{x+y=s+m}(T_{-\om_2-2(y+j-l)}T_{klx}^{\om_1N}) 
T_{ljy}^{\om_2N}
(T_{N-2j}C_{m,xy}^{\om_1+\om_2+2s,\om_1\om_2})(e_{x+l-k}\ot e_{y+j-l}).
\end{split}\end{multline*}
Identifying the coefficient of $e_x\otimes e_y$, we obtain
\begin{multline}\label{add2}
C_{m+j-k,xy}^{\om_1+\om_2+2s,\om_1\om_2}(\la)T_{kjm}^{\om_1+\om_2+2s,N}(\la)\\
=\sum_{l=\max(0,j-y)}^{\min(N,x+k)}
C_{m,x+k-l,y+l-j}^{\om_1+\om_2+2s,\om_1\om_2}(\la+N-2j)
T_{kl,x+k-l}^{\om_1N}(\la-\om_2-2y) T_{lj,y+l-j}^{\om_2N}(\la)
,\end{multline}
where $k+x+y=s+m+j$. 

To identify \eqref{add2} as a special case of \eqref{pi}, we first plug in 
the expressions from Proposition \ref{tl} and \eqref{cga}. We then transform 
the $q$-Racah polynomials coming from $T_{kjm}$ and $T_{lj,y+l-j}$ using 
\eqref{rso} and  the one coming from $T_{kl,x+k-l}$ using \eqref{rs1}. 
Finally we 
replace $l$ by $x+k-l$ in the summation. As the patient reader can  verify, 
we obtain \eqref{pi} in base $q^2$ with 
\begin{multline*}(k_1,k_2,m_1,m_2,m_3,a,b,c,d)\\
\mapsto(k,x,j,m,s,q^{-2(N+1)},
q^{2(\la-\om_1-\om_2+N+1-2j-2m-2s)},q^{2(\om_1-1)},q^{-2(\la+1+N-j-m-s)}).
\end{multline*}

In the case of the group $\mathrm{SU}(2)$, 
\eqref{mer} for $\De(t_{kk}^{2k})$ (the spherical case)
is the classical addition formula for Legendre polynomials. In 
\cite{k2}, Koornwinder showed that for $\mathcal
F_q(\mathrm{SU}(2))$, one may obtain an addition formula for 
little $q$-Legendre polynomials by evaluating the corresponding identity
in a tensor product of infinite-dimensional representations. 
 This became the starting point for much work on
quantum groups and $q$-special functions, cf.~\cite{kf}.
In view of Remark \ref{ndl} one might expect that Koornwinder's formula 
is a limit case of the spherical case of 
\eqref{add2}, and thus of the Biedenharn--Elliott identity \eqref{pi}.
However, since the formal limit $\frl\rightarrow\mathcal F_q(\mathrm{SL}(2))$
involved here does not preserve the coproduct, this is not at all clear a 
priori.  Nevertheless, such a limit transition exists, as will be explained
 below.
We stress that the general case of \eqref{pi} had not yet appeared in the 
literature when \cite{k2} was published.

To rewrite \eqref{add2} in a form similar to the addition formula of \cite{k2},
we must express the  functions $T_{kjm}^{\om N}$ in terms of appropriate
 Askey--Wilson polynomials.  Namely, when $k\leq j$, $k+j\leq N$ we
 use the expression \eqref{taw} for $T=T_{kjm}^{\om N}$. When $k\leq j$, 
$k+j\geq N$ we write 
\begin{multline*}
T_{kjm}^{\om N}(\la)=(-1)^{N-k}q^{(N-j)(2\la-\om+3+N+k-4j-m)+k(1-m)-N}\\
\begin{split}&\times\frac{(q^{2(\la-\om+2+N-2j-m)};q^2)_{k+j-N}}
{(q^2,q^{2(\la+1-j)};q^2)_{N-j}(q^{2(\la-\om+2+N-2j-2m)};q^2)_{k}}\\
&\times
p_{N-j}(\textstyle\frac
{q^{\om-1}+q^{1-\om}}2\displaystyle;q^{1+\om+2(m+j-k)},q^{1-\om-2m},q^{3-\om+2
(\la+k-j-m)},q^{\om-1-2(\la+N-2j-m)};q^2).
\end{split}\end{multline*}
When $j\leq k$, $k+j\leq N$, we write 
\begin{multline*}
T_{kjm}^{\om N}(\la)=(-1)^{k}q^{k(2\la+3+N-j-m)-j(\om+1+N+3m)+mN}\\
\begin{split}&\times\frac{(q^{2(\la+1+k-j-m)};q^2)_{N-k-j}(q^{-2m},
q^{-2(m+\om-1)};q^2)_{k-j}}
{(q^2;q^2)_j(q^{2(\la+1-j)};q^2)_{N-j}(q^{2(\la-\om+2+N-2j-2m)};q^2)_{k}}\\
&\times
p_{j}(\textstyle\frac {q^{\om-1}+q^{1-\om}}
2\displaystyle;q^{1+\om+2m},q^{1-\om-2(m+j-k)},
q^{3-\om+2(\la+N-2j-m)},q^{\om-1-2(\la+k-j-m)};q^2),
\end{split}\end{multline*}
while, finally, if $j\leq k$, $N\leq k+j$ we write
\begin{multline*}
T_{kjm}^{\om
N}(\la)=(-1)^{N-j}q^{(N-j)(2\la-\om+2+N+k-4j-m)+(k-j)(\om+1+j-k+2m)-mk}\\
\begin{split}&\times\frac{(q^2;q^2)_k(q^{-2m},q^{-2(m+\om-1)};q^2)_{k-j}
(q^{2(\la-\om+2+N-2j-m)};q^2)_{k+j-N}}
{(q^2;q^2)_j(q^2,q^{2(\la+1-j)};q^2)_{N-j}(q^{2(\la-\om+2+N-2j-2m)};q^2)_{k}}\\
&\times
p_{N-k}(\textstyle\frac{q^{\om-1}+q^{1-\om}}
2\displaystyle;q^{1+\om+2m},q^{1-\om-2(m+j-k)},
q^{3-\om+2(\la+k-j-m)},q^{\om-1-2(\la+N-2j-m)};q^2).
\end{split}\end{multline*}
These expressions can be obtained from Proposition \ref{tl} using 
transformation formulas
from \cite{gr}, or directly from Theorem \ref{lmep}. 
If we insist on using these expressions, we must split the sum in
\eqref{add2} into five parts, according to whether $l$ is smaller or larger 
than
 $j,\,k,\,N-j$ and $N-k$. We will  write this out explicitly only in the 
spherical case, when the four splitting points agree. 

Thus we let $j=k$, $N=2k$ in \eqref{add2}. After rewriting the sum as
$$\sum_{l=\max(0,k-y)}^{\min(2k,x+k)}a_l=a_k+\sum_{l=1}^{\min(k,y)}a_{k-l}
+\sum_{l=1}^{\min(k,x)}a_{k+l},$$
we express the functions $T_{kjm}^{\om N}$   
as indicated above and  the Clebsch--Gordan coefficients using  \eqref{cga}.
This results in the identity
\begin{multline*}
 q^{-k} (q^2,q^{2(\la-\om_2+2-2y)},q^{2(\om_2-\la+2y)};q^2)_k\\
\begin{split}&\quad\times R_x(\mu(m);q^{2(\la-\om_1-\om_2+1-2L)},
q^{2(\om_1-1)},q^{-2(L+1)},q^{-2(\om_1+\om_2+L-1)};q^2)\\
&\quad\times p_k(\Om;q^{1+\om_1+\om_2+2L},q^{1-\om_1-\om_2-2L},
q^{3-\om_1-\om_2+2(\la-L)},q^{\om_1+\om_2-1-2(\la-L)};q^2)\\
 &=R_x(\mu(m);q^{2(\la-\om_1-\om_2+1-2L)},q^{2(\om_1-1)},q^{-2(L+1)},
q^{-2(\om_1+\om_2+L-1)};q^2)\\
&\quad\times p_{k}(\Om_1;q^{1+\om_1+2x},q^{1-\om_1-2x},
q^{3-\om_1+2(\la-\om_2-x-2y)},q^{\om_1-1-2(\la-\om_2-x-2y)};q^2)\\
&\quad\times p_{k}(\Om_2;q^{1+\om_2+2y},q^{1-\om_2-2y},q^{3-\om_2+2(\la-y)},
q^{\om_2-1-2(\la-y)};q^2)\\
 &+\sum_{l=1}^{\min(k,y)}\frac{(q^2;q^2)_{k+l}(q^2;q^2)_y}
{(q^2;q^2)_{k-l}(q^2;q^2)_{y-l}} \frac{1-q^{2(\la-\om_2+1-2y+2l)}}
{1-q^{2(\la-\om_2+1-2y)}} 
\,q^{l(\om_1+3\om_2-2\la-4+2k+2x+4y)} \\
&\quad\times\frac{(q^{2(1-\om_2-y)},q^{2(\la-\om_1-\om_2+2-x-2y)},
q^{2(\la-\om_2+1-x-2y)},q^{2(\la-\om_2+1-k-2y)};q^2)_l}
{(q^{2(\la-\om_2+2+k-2y)};q^2)_l}\\
&\quad\times R_{x+l}(\mu(m);q^{2(\la-\om_1-\om_2+1-2L)},q^{2(\om_1-1)},
q^{-2(L+1)},q^{-2(\om_1+\om_2+L-1)};q^2)\\
&\quad\times p_{k-l}(\Om_1;q^{1+\om_1+2(x+l)},q^{1-\om_1-2x},
q^{3-\om_1+2(\la-\om_2-x-2y+l)},q^{\om_1-1-2(\la-\om_2-x-2y)};q^2)\\
&\quad\times p_{k-l}(\Om_2;q^{1+\om_2+2y},q^{1-\om_2-2(y-l)},
q^{3-\om_2+2(\la-y+l)},q^{\om_2-1-2(\la-y)};q^2)\\
&+\sum_{l=1}^{\min(k,x)}\frac{(q^2;q^2)_{k+l}(q^2;q^2)_x}
{(q^2;q^2)_{k-l}(q^2;q^2)_{x-l}}\frac{1-q^{2(\om_2-\la-1+2y+2l)}}
{1-q^{2(\om_2-\la-1+2y)}}\,q^{l(\om_1-\om_2+2\la+2k-2y)}\\
&\quad\times\frac{(q^{-2(\la-y)},q^{2(1-\om_1-x)},q^{2(\om_2-\la-1+y)},
q^{2(\om_2-\la-1+2y-k)};q^2)_l}{(q^{2(\om_2-\la+k+2y)};q^2)_l} \\
&\quad\times R_{x-l}(\mu(m);q^{2(\la-\om_1-\om_2+1-2L)},q^{2(\om_1-1)},
q^{-2(L+1)},q^{-2(\om_1+\om_2+L-1)};q^2)\\
&\quad\times p_{k-l}(\Om_1;q^{1+\om_1+2x},q^{1-\om_1-2(x-l)},
q^{3-\om_1+2(\la-\om_2-x-2y)},q^{\om_1-1-2(\la-\om_2-x-2y-l)};q^2)\\
&\quad\times p_{k-l}(\Om_2;q^{1+\om_2+2(y+l)},q^{1-\om_2-2y},
q^{3-\om_2+2(\la-y)},q^{\om_2-1-2(\la-y-l)};q^2),
\end{split}\end{multline*}
where $L=x+y=s+m$ and we write
$$\Om_1=\frac{q^{\om_1-1}+q^{1-\om_1}}{2},\ \ \ \Om_2
=\frac{q^{\om_2-1}+q^{1-\om_2}}{2},\ \ \ \Om=
\frac{q^{\om_1+\om_2+2s-1}+q^{1-\om_1-\om_2-2s}}2.$$ 
This identity generalizes Koornwinder's addition formula to the level of 
Askey--Wilson polynomials, but is itself a special case of \eqref{pi}.
It is also a special case of an identity obtained in \cite{kf},
using the non-dynamical quantum group and twisted primitive elements.
Note that since $\mu(m)=2q^{1-2L-\om_1-\om_2}\Om$, we may view it as 
a linearization formula, expanding a product of two polynomials in $\mu(m)$
into $q$-Racah polynomials. 

To obtain Koornwinder's formula as a limit case
(cf.~also \cite[Remark 5.2]{kf}) , 
we fix $k$, $x$ and $m$. We then let 
 $y,s,L,\la,\om_1+\om_2+L,-\om_2-L\rightarrow\infty$ in such a way that
$L=x+y=s+m$ and  $\la-L\rightarrow z$ for some constant $z$. 
Using limit transitions from \cite{ks}, or directly from the definitions, 
one may check that both the $q$-Racah  and the Askey--Wilson polynomials 
tend to little $q$-Jacobi polynomials, defined by (cf.~\cite{gr})
$$p_n(x;a,b;q)={}_2\phi_1\left[\begin{matrix}
q^{-n},abq^{n+1}\\aq\end{matrix};q,qx\right]. $$
 In the limit, one obtains the identity
\begin{multline*}
  p_{x}(q^{2m};q^{2z},0;q^2)\, p_k(q^{2m};1,1;q^2)\\
\begin{split}&= p_{x}(q^{2m};q^{2z},0;q^2)
\,p_{k}(q^{2x};1,1;q^2)\,p_{k}(q^{2(x+z)};1,1;q^2)\\
 &\quad+\sum_{l=1}^{k}\frac{(q^2;q^2)_{k+l}(q^{2(1+z)};q^2)_{x+l}}
{(q^2;q^2)_{k-l}(q^{2(1+z)};q^2)_x(q^2;q^2)_l^2}\,q^{2l(x+l-k)} \\
&\qquad\times  p_{x+l}(q^{2m};q^{2z},0;q^2)\,p_{k-l}(q^{2x};q^{2l},q^{2l};q^2)
\,p_{k-l}(q^{2(x+z)};q^{2l},q^{2l};q^2)\\
&\quad+\sum_{l=1}^{\min(k,x)}\frac{(q^2;q^2)_{k+l}(q^2;q^2)_x}
{(q^2;q^2)_{k-l}(q^2;q^2)_{x-l}(q^2;q^2)_l^2}\,q^{2l(x+z-k+1)}\\
&\qquad\times p_{x-l}(q^{2m};q^{2z},0;q^2)\,
p_{k-l}(q^{2(x-l)};q^{2l},q^{2l};q^2)\,p_{k-l}(q^{2(x+z-l)};q^{2l},q^{2l};q^2).
\end{split}\end{multline*}
For $z\in \Zp$, this is Koornwinder's formula. 

\section{The Haar functional}
\label{hfs}

In this section we will show that there is a natural  Haar
functional on $\frl$, and that it
 can be identified with a special case of the Askey--Wilson measure.

To motivate our definition, note that   
 the Haar functional on a compact group can be obtained as
 the projection from the regular representation to the isotypic subspace
 containing the trivial representation. Since the trivial representation
  occurs with  multiplicity
 one, the range of the Haar functional can be identified with $\mathbb C$. 
 In the present case, cf.~Remark \ref{lrd}, the trivial representation occurs
 with infinite multiplicity, and the corresponding isotypic component is
  $\mu_l(\mg)\mu_r(\mg)\simeq\mg\ot\mg$. 
The projection onto this space is
\begin{equation}\label{h}h(f(\la)g(\mu)t_{kj}^N)=f(\la)g(\mu)\,\de_{0N},
\ \ \ 0\leq j,k\leq N.\end{equation}
By Proposition \ref{pw}, this defines a $\mathbb C$-linear map
$\frl\rightarrow\mu_l(\mg)\mu_r(\mg)$, which is an $\gh$-prealgebra
homomorphism, cf.~\S \ref{sscoalg}. We call this map the \emph{Haar functional}
on $\frl$.  

We define a \emph{left-invariant integral} on an $\gh$-bialgebroid $A$ to be an
$\gh$-prealgebra homomorphism $h:\,A\rightarrow\mu_l^A(\mg)\mu_r^A(\mg)
\subseteq A$ such that, under the identifications \eqref{can},
$$ h=(\id\ot\,\ep\circ h)\circ\De.$$
If this condition is replaced  with
$$ h=(\ep\circ h\ot\id)\circ\De$$
we speak of a \emph{right-invariant integral}.
If $a=\sum_i a_i'\ot a_i''$, and we write $\mu_l(f)=f(\la)$, $\mu_r(f)=f(\mu)$,
then left-invariance means that
$$h(a)(\la,\mu)=\sum_i h(a_i'')(\mu,\mu)\,a_i'$$
and right-invariance that
$$h(a)(\la,\mu)=\sum_{i}h(a_i')(\la,\la)\,a_i''.$$
These definitions are motivated by the following fact.

\begin{proposition}
The Haar functional \eqref{h}
  is both the unique  
left-invariant  and the unique right-invariant integral on $\frl$ such that
$h(1)=1$.
\end{proposition}

\begin{proof}
Let $h$ be a left-invariant integral on $\frl$. Then
 $h(f(\la)g(\mu)t_{kj}^N)= f(\la)g(\mu)h_{kj}^N(\la,\mu)$
for some $h_{kj}^N\in\mg\ot\mg$.
The left-invariance of $h$ means that 
$$h_{kj}^N(\la,\mu)=\sum_{l=0}^N h_{lj}^N(\mu,\mu) t_{kl}^N.$$
By Proposition \ref{pw}, this implies that $h_{kj}^N=0$ unless
$N=0$, in which case the normalizing condition shows that $h$ is given by
\eqref{h}.
For the right-invariance, the proof is similar.
\end{proof}

We will now obtain the Schur
orthogonality relations for matrix elements. These are most elegantly 
discussed 
in terms of contragredient
corepresentations, but to save space we give a direct proof using results
obtained above.

\begin{theorem}\label{sop}
In the notation above, the 
 following Schur orthogonality relations are valid:
\begin{equation*}\begin{split}h(t_{jk}^M(t_{lm}^N)^\ast)&=\de_{MN}
\,\de_{jl}\,\de_{km}\,q^{2(M-k)}
\frac{1-q^2}{1-q^{2(M+1)}}\\
&\quad\times\frac{(q^2,q^{-2\la};q^2)_j(q^2;q^2)_{M-j}
(q^{-2(\mu+1+M-k)};q^2)_{M-k}}
{(q^{-2(\la+M+1-j)};q^2)_j(q^2;q^2)_k
(q^2,q^{-2(\mu+M-2k)};q^2)_{M-k}}.\end{split}\end{equation*}
\end{theorem}

\begin{proof}
  Applying $h$ to \eqref{lin} gives
\begin{multline}\label{htt}h(t_{kx}^M t_{ly}^N)=\de_{MN}\,
\de_{k+l,M}\,\de_{x+y,M}\,
\frac{\Ga_x^M(\mu)\Ga_y^M(\mu+M-2x)}{Z_M\Ga_{0}^{0}(\mu)}\,C_{kl,0}^{MM,0}(\la)
C_{xy,0}^{MM,0}(\mu)\\
=\de_{MN}\,
\de_{k+l,M}\,\de_{x+y,M}(-q)^{l+y}
\frac{(q^2,q^{-2\la};q^2)_k(q^2;q^2)_l(q^{-2(\mu+1+y)};q^2)_y}
{(q^{-2(\la+1+l)};q^2)_k(q^2;q^2)_x
(q^2,q^{-2(\mu+y-x)};q^2)_y}
\frac{1-q^2}{1-q^{2(M+1)}},
\end{multline}
where we inserted the expressions
 \eqref{csp} for the Clebsch--Gordan coefficients. 
Next we observe that
\begin{equation}\label{tast}(t_{kj}^N)^\ast=(-q)^{k-j}t_{N-k,N-j}^N.
\end{equation}
This follows easily from Proposition \ref{mep}, and can also
 be proved without using explicit expressions for the matrix
elements, similarly to the proof of Proposition \ref{un}.
Combining \eqref{htt} and \eqref{tast} completes the proof.
\end{proof}

Note that
 $h$ vanishes outside  $\frl_{00}$, which can be identified with the algebra
of polynomials in $\Xi$ (cf.~Lemma \ref{lx}) 
over the meromorphic functions in $\la$ and $\mu$. 
It is natural to seek a family of measures $dm_{\la\mu}$ so that
\begin{equation}\label{hp}h(p(\Xi))=\int p(x)\,dm_{\la\mu}(x)\end{equation}
for any such polynomial $p$. 
By Proposition \ref{pw}, it is enough to do this for $t_{kk}^{2k}$, $k\in\Zp$. 
It follows from Theorem \ref{lmep} that
$$t_{kk}^{2k}=N_k\, p_k(\textstyle\frac 1 2\,\displaystyle\Xi),$$
with $N_k=q^{k(\la-\mu+1)}/(q^2,q^{-2\mu};q^2)_k$,
$p_k(x)=p_k(x;q^{\mu-\la+1},q^{\la-\mu+1},q^{\la+\mu+3},q^{-\la-\mu-1};q^2).$
For $\la,\,\mu\notin\mathbb Z$, the polynomials 
$\{ p_k\}_{k=0}^\infty$ form an orthogonal system 
with respect to a moment functional, 
which is given by integration with respect to
an explicitly known (not necessarily positive) measure,
cf.~\cite{aw2}. Under additional conditions
on $\la$ and $\mu$, e.g.~$\la,\mu\in (j,j+1)$ for
$j\in\mathbb Z$ or $\Ima\la=\Ima\mu=\pi/2\log q$, the 
measure is positive. Assuming that $\la,\,\mu\notin\mathbb Z$, we 
 let $dm_{\la\mu}$ be the orthogonality measure, rescaled and normalized  so 
that
$$C_k\,\de_{kl}=\int p_k(\textstyle\frac x2\displaystyle)\,
p_l(\textstyle\frac x2\displaystyle)\,dm_{\la\mu}(x)$$
with $C_0=1$. Then, in particular,
 \eqref{hp} is satisfied.
 Therefore, the Haar measure on $\frl$ 
 can be identified with the orthogonality measure for a 
two-parameter family of Askey--Wilson poly\-nomials. 
A similar interpretation is obtained by Koornwinder in the non-dynamical
case \cite{kz}, using twisted primitive elements.

Let us now combine \eqref{hp} with the case of Theorem \ref{sop} when
  $h$   is applied to an element of $\frl_{00}$. It will be no restriction 
to assume that we are in the first parameter domain of Theorem \ref{lmep}. 
Thus we  consider the element 
$t_{j,j+k}^{2j+k+m}(t_{l,l+k}^{2l+k+m})^\ast$, where 
 $j,\,k,\,l,\,m\in\Zp$. Using  Theorem \ref{lmep}, Lemma \ref{cr} and 
Lemma \ref{ccr}, we can write
\begin{multline*} t_{j,j+k}^{2j+k+m}(t_{l,l+k}^{2l+k+m})^\ast
=(-1)^mq^{(j+l)(\la-\mu+1+k)-m(2\mu+1+m)}\frac{(q^{-2\mu};q^2)_{k-m}}
{(q^{-2\mu};q^2)_{j+k}(q^{-2\mu};q^2)_{l+k}}\\
\begin{split}&\times\frac{(q^2;q^2)_{j+k+m}(q^2;q^2)_{l+k+m}}
{(q^2;q^2)_{j+k}(q^2;q^2)_{j+m}(q^2;q^2)_{l+k}(q^2;q^2)_{l+m}}
\frac{(q^{2(\la+2)};q^2)_{k+m}}
{(q^{2(\la+2)};q^2)_{j+k+m}(q^{2(\la+2)};q^2)_{l+k+m}}\\
&\times h_k(\textstyle\frac 12\,\Xi,q^{\la-\mu+1};q^2)
\,h_m(\textstyle\frac 12\,\Xi,q^{\la+\mu+3};q^2)
\, p_j^{(k,m)}(\textstyle\frac 12\,\Xi)\, 
p_l^{(k,m)}(\textstyle\frac 12\,\Xi),\end{split}\end{multline*}
where we use the notation \eqref{nps} if $k<m$ and write (cf.~\eqref{nm})
$$p_j^{(k,m)}(x)=p_j(x;q^{\mu-\la+1},q^{\la-\mu+1+2k},q^{\la+\mu+3+2m},
q^{-\la-\mu-1};q^2).$$
Combining this with Theorem \ref{sop} and \eqref{hp}, we get after 
simplifications
\begin{equation*}\begin{split}&\quad\int p_j^{(k,m)}(\textstyle
\frac x2\displaystyle)\, p_l^{(k,m)}(\textstyle\frac x2\displaystyle)
\,h_k(\textstyle\frac x2,q^{\la-\mu+1};q^2)\,
h_m(\textstyle\frac x2,q^{\la+\mu+3};q^2)\,dm_{\la\mu}(x)\\
&=\de_{jl}\,\frac{1-q^2}{1-q^{2(2j+k+m+1)}}
\frac{(q^2;q^2)_j(q^2;q^2)_{j+k}(q^2;q^2)_{j+m}}{(q^2;q^2)_{j+k+m}}\\
&\quad\times(q^{-2\la};q^2)_j(q^{2(\la+2)};q^2)_{j+k+m}
(q^{-2\mu};q^2)_{j+k}(q^{2(\mu+2)};q^2)_{j+m}.
\end{split}\end{equation*}
Thus we obtain the orthogonality of a four-parameter family of Askey--Wilson 
polynomials. The above expression agrees with the known explicit formulas 
for the ortho\-gonality measure and the norm.

\section*{Appendix 1. Properties of the antipode}
\setcounter{section}{1}
\renewcommand{\thesection}{\Alph{section}}
\setcounter{equation}{0}

In this appendix we prove Proposition \ref{al} and Lemma \ref{sal}, though we 
leave many details to the
reader. Recall that when proving the corresponding statements for Hopf algebras
 (cf.~\cite{k} for a detailed exposition) it is convenient to work with the 
convolution product
\begin{equation}\label{cp}\phi\star \psi=m^B\circ(\phi\otimes \psi)\circ\De^A
\end{equation}
  on $\Hom_\mathbb C(A,B)$, where $A$ is a coalgebra and $B$ an algebra.
We will use the analogous convolution when
$A$ is an $\gh$-coalgebroid (cf.~\S \ref{sscoalg}) and $B$ an $\gh$-algebra. 
 We will write $H_l=H_l(A,B)$ etc.~for the
spaces
\begin{align*}H_l&=\{\phi\in\Hom_\mathbb C(A,B);\  
\phi(\mu_l^A(f)a)=\mu_l^B(f)\phi(a)
\ \ \text{for all}\ a\in A,\ f\in \mg\},\\
H_r&=\{\phi\in\Hom_\mathbb C(A,B);\  
\phi(a\mu_r^A(f))=\phi(a)\mu_r^B(f)
\ \ \text{for all}\ a\in A,\ f\in \mg\},\\
H_l^{\op}&=\{\phi\in\Hom_\mathbb C(A,B);\  
\phi(a\mu_l^A(f))=\mu_r^B(f)\phi(a)
\ \ \text{for all}\ a\in A,\ f\in \mg\},\\
H_r^{\op}&=\{\phi\in\Hom_\mathbb C(A,B);\ 
\phi(\mu_r^A(f)a)=\phi(a)\mu_l^B(f)
\ \ \text{for all}\ a\in A,\ f\in \mg\}.
\end{align*}
Note that, for $\gh\neq 0$,  the convolution product (\ref{cp}) 
is not globally defined
on $\Hom_\mathbb C(A,B)\times \Hom_\mathbb C(A,B)$, since 
$m\circ(\phi\otimes\psi)$ need not factor
through relation (\ref{mtr}). A sufficient condition for $\phi\star\psi$ to be
well-defined is
$$\phi(\mu_r^A(f)a)\psi(b)=\phi(a)\psi(\mu_l^B(f)b),\ \ \ a\in A_{\al\be},\, 
b\in
A_{\be\ga},\,f\in\mg.$$
Using this condition one proves the following lemma. 

\begin{lemma}\label{cl} The convolution $\star$ is well-defined on 
$H_r^{\op}\times
H_l$ and on $H_r\times H_l^{\op}$. The associative law
$(\phi\star\psi)\star\chi=\phi\star(\psi\star\chi)$ holds whenever both sides
are well-defined. \end{lemma}

We now define $1_l=1_l^{(A,B)},\,1_r=1_r^{(A,B)}\in \Hom_\mathbb C(A,B)$ by
$$1_l(a)=\mu_r^B(T_\al(\ep^A(a)1)),\ \ \ \ \ 1_r(a)=\mu_l^B(\ep^A(a)1),
\ \ \ a\in A_{\al\be}.$$
These elements are functorial in the sense that if 
$\chi:\,A_1\rightarrow A_2$
is an $\gh$-coalgebroid homomorphism and $\om:\, B_1\rightarrow B_2$ an
$\gh$-prealgebra homomorphism, then
\begin{equation}\label{fu}1_x^{(A_1,B)}=1_x^{(A_2,B)}\circ\chi,
\ \ \ \ \  1_x^{(A,B_2)}=\om\circ 1_x^{(A,B_1)},\ \ \ x=l,\,r.\end{equation}
One easily checks that $1_l\in H_r\cap H_l^{\op}$, $1_r\in H_l\cap H_r^{\op}$. 
In particular, the following lemma is meaningful. Note that it depends 
crucially on the dynamical shift in the definition of $1_l$.

\begin{lemma}
The elements $1_l$ and $1_r$ satisfy
$$1_l\star\phi=\phi,\ \ \ \ \ \psi\star 1_r=\psi,\ \ \ \phi\in H_l^{\op},\
\psi\in H_r^{\op}.$$
\end{lemma}

\begin{proof}
We write out the proof for $1_l$, the case of $1_r$ being slightly easier.
It suffices to evaluate both sides on
$a\in A_{\al\be}$. Let us write $\De(a)=\sum_i a_i'\otimes a_i''$,
$\ep(a_i')=f_i T_{-\al}$. 
One then has
\begin{equation}\label{1sp}(1_l\star\phi)(a)=\sum_i 1_l(a_i')\phi(a_i'')=\sum_i
\mu_r(T_\al f_i)\phi(a_i'').\end{equation}
On the other hand, 
$$(\ep\otimes\id)\circ\De(a)=\sum_i\ep(a_i')
\otimes a_i''= \sum_i f_i T_{-\al}\otimes a_i''=\sum_i T_{-\al}\otimes 
\mu_l(f_i)a_i'',
$$
using (\ref{mtr}) in the last step.
By the counit axioms this implies
$$a=\sum_i\mu_l(f_i)a_i''=\sum_i a_i''\,\mu_l(T_\al f_i),$$
using (\ref{dr}) and the fact that $f_i=0$ unless $a_i'\in A_{\al\al}$, 
$a_i''\in A_{\al\be}$. Applying $\phi$ to this identity and using that 
$\phi\in H_l^{\op}$
gives 
$$\phi(a)=\sum_i \mu_r(T_\al f_i)\phi(\al_i'').$$
Comparing with (\ref{1sp}) we see that $1_l\star\phi=\phi$.
\end{proof}

We can now begin the proof of Proposition \ref{al}. 
The antipode axioms can be written as $S\in H_l^{\op}\cap
H_r^{\op}$, 
$S\star\id=1_l$, $\id\star\, S=1_r$. If $S$ and $T$ are two such maps, it
follows from the previous lemmas that
$$S=S\star 1_r=S\star(\id\star\, T)=(S\star\id)\star T=1_l\star T=T.$$
This proves  the uniqueness of the antipode.

It is easy to check that, for any $0\neq f\in \mg$,
the maps
$S_f(a)=\mu_l(f^{-1})S(a\mu_r(f))$ and $S'_f(a)=S(\mu_l(f)a)\mu_r(f^{-1})$ 
satisfy
the antipode axioms. By the uniqueness it follows that $S=S_f=S'_{f}$,
which means that 
\begin{equation}\label{c2}S(\mu_l(f)a)=S(a)\mu_r(f),\ \ \ \ \ S(a\mu_r(f))
=\mu_l(f)S(a),\ \ \ a\in A,\ f\in\mg.\end{equation}
Together with \eqref{sf}, this
 implies \eqref{c1}.

To show the first identity of (\ref{aprop}) one
defines $\phi(a\otimes b)=S(b)S(a)$, $\psi(a\otimes b)=S(ab)$  and checks
that they factor into maps $\phi\in H_l^{\op}(A\widehat\otimes A,A)$, 
$\psi\in H_r^{\op}(A\widehat\otimes A,A)$ satisfying 
\begin{equation}\label{mc}m\star\phi=1_r,\ \ \ \ \ \psi\star
m=1_l.\end{equation}
 One needs (\ref{c2}) to show that $\phi$ factors through
relation (\ref{fh}). It follows that 
$$\phi=1_l\star\phi=\psi\star m \star\phi=\psi\star 1_r=\psi.$$
Similarly one defines $\pi=\De\circ S$, 
$\ro=\sigma\circ(S\otimes S)\circ\De$ and checks that 
$\pi\in H_l^{\op}(A,A\widetilde\otimes A)$, 
$\ro\in H_r^{\op}(A,A\widetilde\otimes A)$, 
\begin{equation}\label{dc}\Da\star\pi=1_r,
\ \ \ \ \ \ro\star\Da=1_l,\end{equation}  
which  implies that
$\pi=\ro$. One needs (\ref{c1}) to show that $\si\circ(S\otimes S)$ factors 
through  (\ref{mtr}), so that $\ro$ is well-defined.

We write down the proof of (\ref{mc}) and (\ref{dc}) and leave the remaining
details to the
reader. In the case of (\ref{mc}), it suffices to evaluate both sides on 
$a\otimes b\in A_{\al\be}\otimes A_{\ga\de}$. 
With notation as in (\ref{de}), one then has
\begin{multline*}(m\star\phi)(a\otimes b)=\sum_{ij}m(a_i'\otimes
b_j')\phi(a_i''\otimes b_j'')=\sum_{ij}a_i'b_j'S(b_j'')S(a_i'')\\
\begin{split}&=\sum_i a_i'(\id\star\, S)(b)S(a_i'')=\sum_i a_i'\, 
\mu_l(\ep(b)1)
S(a''_i)=\mu_l(T_{-\al}\ep(b)1)(\id\star\, S)(a)\\
&=\mu_l(T_{-\al}\ep(b)1)\,\mu_l(\ep(a)1)=\mu_l(\ep(a)\ep(b)1)
=1_r^{(A\widehat\otimes
A,A)}(a\otimes b),\end{split}\end{multline*}
using in the penultimate step that $\ep(a)\in \mg T_{-\al}$.
Similarly,
$$(\psi\star m)(a\otimes b)=\sum_{ij}S(a'_ib'_j)a''_ib''_j=(S\star\id)(ab)
=1_l^{(A,A)}(ab)=1_l^{(A\widehat\otimes A,A)}(a\otimes b),
$$
 where we in the last step used (\ref{fu}) with $\chi=m$. 
The first equation of (\ref{dc}) follows from
$$(\De\star\pi)(a)=\sum_i\De(a_i')\De(S(a_i''))
=\De\left(\sum_i a_i'S(a_i'')\right)=\De(1_r^{(A,A)}(a))
=1_r^{(A,A\widetilde\otimes
A)}(a).$$
For the last equation we use notation such as 
$(\De\otimes\id)\circ\De(a)=\sum_i a_i^1\otimes a_i^2\otimes a_i^3$.  
Then (note that the symbols $a_i^j$ have different meaning from equation
to equation) 
\begin{multline*}(\ro\star\De)(a)=\sum_i(\sigma(S\otimes
S)\De)(a_i^1)\De(a_i^2)=\sum_i S(a_i^2)a_i^3\otimes S(a_i^1)a_i^4\\
\begin{split}&=\sum_i(S\star\id)(a_i^2)\otimes S(a_i^1)a_i^3=
\sum_i \mu_r(T_{\be_i}\ep(a^2_i)1)
\otimes S(a_i^1)a_i^3\\
&=\sum_i 1\otimes \mu_l(T_{\be_i}\ep(a^2_i)1)S(a^1_i)a^3_i
=\sum_i 1\otimes S(a^1_i)\mu_l(\ep(a^2_i)1)a^3_i\\
&=1\otimes(S\star 1_r\star\id)(a)=1\otimes(S\star\id)(a)
=1\otimes \mu_r(T_\al\ep(a)1)=1_l^{(A,A\widetilde\otimes
A)}(a),\end{split}\end{multline*}
where we assume that in a three-fold tensor product $a_i^1\in A_{\al\be_i}$,
$a_i^2\in A_{\be_i\ga_i}$. In the sixth equality we used that, by (\ref{c1}),
$S(a_i^1)\in A_{-\be_i,-\al}$. 

That $S(1)=1$ is the case $a=1$ of (\ref{anti}). Then \eqref{sm} is obtained as
the case $a=1$ of \eqref{sf}.
To prove that $\ep\circ
S=S^{D_{\gh}}\circ\ep$, we write
\begin{multline*}\ep(S(x))=\ep((1_l\star S)(x))=\ep\left(\sum_i
1_l^{(A,A)}(x_i')S(x_i'')\right)=\sum_i 1_l^{(A,D_{\gh})}(x_i')\,
\ep(S(x_i''))\\
\begin{split}&=\sum_i T_{\al}(\ep(x_i')1)\,\ep(S(x_i''))
=\sum_i T_{\al}(\ep(x_i')\ep(S(x_i''))1)T_\be
=T_\al(\ep((\id\star\, S)(x))1)T_{\be}\\
&=T_\al(\ep(1_r^{(A,A)}(x))1)T_{\be}=T_\al(1_r^{(A,D_{\gh})}(x)1)T_{\be}
=T_\al(\ep(x)1)T_\be=S^{D_{\gh}}(\ep(x)),
\end{split}\end{multline*}
where $x\in A_{\al\be}$. In the fifth equality we used that, by
(\ref{c1}), $\ep(S(x_i''))\in \mg T_{\be}$.

We now turn to the proof of the second part of Proposition \ref{al}. Let
$X\subseteq A$ and $S\in\End_\mathbb C(A)$ satisfy the conditions stated there.
It is clear that (\ref{sf}) holds. Since the relations (\ref{anti}) are linear
in $a$, it suffices to show that if they hold for $a$ and $b$, they hold for
$ab$. For the relation $\id\star\, S=1_r$, we write
$$(\id\star\,
S)(ab)=\sum_{ij}a_i'b_j'S(a_i''b_j'')=\sum_{ij}a_i'b_j'S(b_j'')S(a_i'').$$
From the proof that $m\star\phi=1_r$ given above we see that this equals
$1_r(ab)$ if the relation holds for $a$ and $b$. 
For  $S\star\id$, the proof is similar.   

Finally we turn to the proof of Lemma \ref{sal}. By the uniqueness of the 
antipode, it suffices to check that
$\tilde S=\ast\circ S^{-1}\circ\ast$ satisfies the antipode axioms.  
We write down 
the proof that $\tilde S$ satisfies the first equality of (\ref{anti}) and 
leave the remaining details to the reader. One has 
\begin{equation*}\begin{split}m\circ(\id\otimes\,\tilde
S)\circ\De&=m\circ(\ast\otimes\ast)\circ(\id\otimes\, S^{-1})\circ
(\ast\otimes\ast)\circ\De\\
&=\ast\circ m\circ\si\circ(\id\otimes\, S^{-1})\circ\De\circ\ast\\
&=\ast\circ m\circ\si\circ(\id\otimes\, S^{-1})\circ\si\circ(S\otimes
S)\circ\De\circ S^{-1}\circ\ast\\
 &=\ast\circ m\circ(\id\otimes S)\circ\De\circ
 S^{-1}\circ\ast,\end{split}\end{equation*}
where we used (\ref{aprop}) in the third step. Since $S$ satisfies 
(\ref{anti}),
it follows that
\begin{equation*}\begin{split}m\circ(\id\otimes\,\tilde
S)\circ\De(a)&=\mu_l(\ep(S^{-1}(a^\ast))1)^\ast=
\mu_l\left((S^{D_{\gh}})^{-1}\circ\ast^{D_{\gh}}\circ
\ep(a)1\right)^\ast\\
&=\mu_l(\overline{\ep(a)1})^\ast=\mu_l(\ep(a)1),\end{split}\end{equation*}
where we used (\ref{aprop}) and the $\ast$-structure axioms.

\end{document}